\documentclass[final]{siamltex}%
\usepackage{stmaryrd}
\usepackage{amsfonts}
\usepackage{sw20siam}
\usepackage{enumerate}
\usepackage{color}
\usepackage{amsmath}
\usepackage{amssymb}
\usepackage{diagrams}
\usepackage{graphicx}

\newcommand{\Rr}{\mathbb{R}}
\newcommand{\Cc}{\mathbb{C}}
\newcommand{\Dd}{\mathcal{D}}
\newcommand{\Ss}{\mathcal{S}}
\newcommand{\Ssbar}{\overline{\mathcal{S}}}

\newcommand{\Span}{\operatorname{span}}

\newcommand{\transp}{\top}

\newcommand{\trank}{\operatorname{rank}_{\otimes}}
\newcommand{\multirank}{\operatorname{rank}_{\boxplus}}
\newcommand{\brank}{\underline{\operatorname{rank}}_{\otimes}}

\newcommand{\GL}{\operatorname{GL}}
\newcommand{\Oo}{\operatorname{O}}
\newcommand{\xx}{\mathbf{x}}
\newcommand{\yy}{\mathbf{y}}
\newcommand{\zz}{\mathbf{z}}
\newcommand{\uu}{\mathbf{u}}
\newcommand{\vv}{\mathbf{v}}
\newcommand{\ww}{\mathbf{w}}
\newcommand{\ee}{\mathbf{e}}
\newcommand{\ff}{\mathbf{f}}

\newcommand{\app}{{\sc{approx}}}
\newcommand{\brap}{{\app}$(A,r)$}
\newcommand{\Det}{\operatorname{Det}}

\newcommand{\Bb}{\mathcal{B}}
\newcommand{\Bbe}{\mathcal{B}_{\epsilon}}

\begin{document}

\title{Tensor rank and the ill-posedness of the best low-rank
approximation problem}
\author{Vin de Silva\thanks{Department of Mathematics, Pomona
College, Claremont, CA 91711-4411.\newline
\textit{E-mail:} \texttt{vin.desilva@pomona.edu}}
\and Lek-Heng Lim\thanks{Corresponding author. Institute for Computational and
Mathematical Engineering, Stanford University, Stanford, CA 94305-9025.
\textit{E-mail:} \texttt{lekheng@stanford.edu}}}
\maketitle

\begin{abstract}
There has been continued interest in seeking a theorem describing
optimal low-rank approximations to tensors of order~$3$ or higher,
that parallels the Eckart--Young theorem for matrices.
In this paper, we argue that the naive approach to this problem is
doomed to failure because, unlike matrices, tensors of order~$3$ or
higher can fail to have best rank-$r$ approximations. The phenomenon
is much more widespread than one might suspect: examples of this
failure can be constructed over a wide range of dimensions, orders
and ranks, regardless of the choice of norm (or even Br\`{e}gman
divergence). Moreover, we show that in many instances these
counterexamples have positive volume: they cannot be regarded as
isolated phenomena.  In one extreme case, we exhibit a tensor space
in which \emph{no} rank-$3$ tensor has an optimal rank-$2$
approximation.
The notable exceptions to this misbehavior are rank-$1$ tensors and
order-$2$ tensors (i.e.\ matrices).

In a more positive spirit, we propose a natural way of overcoming
the ill-posedness of the low-rank approximation problem, by using
\emph{weak solutions} when true solutions do not exist. For this to
work, it is necessary to characterize the set of weak solutions, and
we do this  in the case of rank~$2$, order~$3$ (in arbitrary
dimensions). In our work we emphasize the importance of closely
studying concrete low-dimensional examples as a first step towards
more general results. To this end, we present a detailed analysis of
equivalence classes of $2 \times 2 \times 2$ tensors, and we develop
methods for extending results upwards to higher orders and
dimensions.

Finally, we link our work to existing studies of tensors from an
algebraic geometric point of view. The rank of a tensor can in
theory be given a semialgebraic description; in other words, can be
determined by a system of polynomial inequalities. We study some of
these polynomials in cases of interest to us; in particular we make
extensive use of the \emph{hyperdeterminant}~$\Delta$ on
$\Rr^{2\times 2 \times 2}$.
\end{abstract}

\keyphrases{numerical multilinear algebra, tensors, multidimensional
arrays, multiway arrays, tensor rank, tensor decompositions, low rank
tensor approximations, hyperdeterminants, Eckart--Young theorem, principal
component analysis, \textsc{parafac}, \textsc{candecomp}, Br\`{e}gman
divergence of tensors}

\AMSclass{14P10, 15A03, 15A21, 15A69, 15A72, 49M27, 62H25, 68P01}

\section{Introduction}

Given an order-$k$ tensor $A\in\Rr^{d_{1}\times\dots\times d_{k}}$, one
is often required to find a \textit{best rank-}$r$\textit{ approximation} to
$A$ --- in other words, determine vectors $\xx_{i}\in\Rr^{d_{1}%
},\yy_{i}\in\Rr^{d_{2}},\dots,$ $\mathbf{z}_{i}\in
\Rr^{d_{k}}$, $i=1,\dots,r$, that minimizes%
\[
\lVert A-\xx_{1}\otimes\yy_{1}\otimes\dots\otimes\mathbf{z}%
_{1}-\dots-\xx_{r}\otimes\yy_{r}\otimes\dots\otimes
\mathbf{z}_{r}\rVert
\]
or, in short,%
\begin{equation}
\operatorname*{argmin}\nolimits_{\operatorname*{rank}_{\otimes}(B)\leq
r}\lVert A-B\rVert.
\tag{\brap}
\label{brap}
\end{equation}
Here $\lVert\cdot\rVert$ denotes some choice of norm on $\Rr%
^{d_{1}\times\dots\times d_{k}}$. When $k=2$, the problem is completely
resolved for unitarily invariant norms on $\Rr^{m\times n}$ with the
\textit{Eckart--Young theorem} \cite{EY}, which states that if
\[
A=U\Sigma V=\sum\nolimits_{i=1}^{\operatorname*{rank}(A)}\sigma_{i}%
\uu_{i}\otimes\vv_{i},\quad\sigma_{i}\geq\sigma_{i+1},
\]
is the singular value decomposition of $A\in\Rr^{m\times n}$, then a
best rank-$r$ approximation is given by the first $r$ terms in the
above sum \cite{GV}. The best rank-$r$ approximation problem for
higher order tensors is a problem of central importance in the
statistical analysis of multiway data \cite{Bro, C, LMV1, LMV2, KR2,
Krus1, HL, LS, SBG1, SBG, VT1, VT2, VT3}.

It is therefore not surprising that there has been continued interest in
finding a satisfactory `singular value decomposition' and an `Eckart--Young
theorem'-like result for tensors of higher order. The view expressed in the
conclusion of \cite{K1} is representative of such efforts and we reproduce it here:

\begin{quotation}
``\textit{An Eckart--Young type of best rank-$r$ approximation
theorem for tensors continues to elude our investigations but can
perhaps eventually be attained by using a different norm or yet
other definitions of orthogonality and rank.}''
\end{quotation}

It will perhaps come as a surprise to the reader that the problem of
finding an `Eckart--Young type theorem' is ill-founded because of a
more fundamental difficulty: the best rank-$r$ approximation problem
{\brap} has no solution in general! This paper seeks to provide an
answer to this and several related questions.

\subsection{Summary}

Since this is a long paper, we present an `executive summary' of selected results, in this section and the next.
%
We begin with the five main objectives of this article:

\newcounter{goals}
\begin{enumerate}
\item \label{obj:manyr}
\textbf{{\brap} is ill-posed for many~$r$.}
We will show that, regardless of the choice of norm, the problem of
determining a best rank-$r$ approximation for an order-$k$ tensor in
$\Rr^{d_{1}\times\dots\times d_{k}}$ has no solution in general for
$r=2,\dots,\min\{d_{1},\dots,d_{k}\}$ and $k\geq3$. In other words,
the best low rank approximation problem for tensors is ill-posed for
all orders (higher than $2$), all norms, and many ranks.

\item \label{obj:manyA}
\textbf{{\brap} is ill-posed for many~$A$.}
We will show that the set of tensors that fail to have a best low
rank approximation has positive volume. In other words, such
failures are not rare --- if one randomly picks a tensor~$A$ in a suitable
tensor space, then there is a non-zero probability that~$A$ will fail to
have a best rank-$r$ approximation for some $r<\trank(A)$.

\item \label{obj:weak}
\textbf{Weak solutions to {\brap}.}
We will propose a natural way to overcome the ill-posedness of the
best rank-$r$ approximation problem with the introduction of `weak
solutions', which we explicitly characterize in the case $r=2$,
$k=3$.

\item \label{obj:semialg}
\textbf{Semialgebraic description of tensor rank.}
From the Tarski--Seidenberg theorem in model theory~\cite{Tar,Sei} we will deduce the following:
for any $d_{1},\dots,d_{k}$, there exists a finite
number of polynomial functions, $\Delta_{1},\dots,\Delta_{m}$,
defined on $\Rr^{d_{1}\times\dots\times d_{k}}$ such that the rank
of any $A\in\Rr^{d_{1}\times\dots\times d_{k}}$ is completely
determined by the signs of $\Delta_{1}(A),\dots,\Delta_{m}(A)$.
We work this out in the special case~$\Rr^{2\times2\times2}$.

\item \label{obj:reduction}
\textbf{Reduction.}
We will give techniques for reducing certain questions about tensors (orbits, invariants, limits) from high-dimensional tensor spaces to lower-dimensional tensor spaces.
%
For instance, if two tensors in
$\Rr^{c_{1}\times\dots\times c_{k}}$ lie in distinct
$\GL_{c_{1},\dots,c_{k}}(\Rr)$-orbits, then they lie in distinct
$\GL_{d_{1},\dots,d_{k}}(\Rr)$-orbits in
$\Rr^{d_{1}\times\dots\times d_{k}}$ for any $d_{i}\geq c_{i}$.

\setcounter{goals}{\value{enumi}}
\end{enumerate}

The first objective is formally stated and proved in
Theorem~\ref{thm:rank+1}. The two notable exceptions where {\brap}
has a solution are the cases $r=1$ (approximation by rank-$1$
tensors) and $k = 2$ ($A$ is a matrix). The standard way to prove
these assertions is to use brute force: show that the sets where the
approximators are to be found may be defined by polynomial
equations. We will provide alternative elementary proofs of these
results in Propositions \ref{prop:goodranka} and
\ref{prop:goodrankb} (see also Proposition \ref{prop:goodrankc}).

The second objective is proved in Theorem~\ref{thm:posvol2}, which
holds true on
 $\mathbb{R}^{d_1 \times d_2 \times d_3}$ for arbitrary $d_1,d_2,d_3 \ge 2$.
Stronger results can hold in specific cases: in Theorem~\ref{thm:norank2}, we will give an instance where \textit{every} rank-$r$ tensor fails to
have a best rank-$(r-1)$ approximator.

The third objective is primarily possible because of the following theorem, which asserts that the boundary of the set of rank-$2$ tensors can be explicitly parameterized. The proof, and a discussion of weak solutions, is given in Section~\ref{sect:char32}.

\begin{theorem}
\label{thm:main}Let $d_{1},d_{2},d_{3}\geq2$. Let $A_{n}\in
\Rr^{d_{1}\times d_{2}\times d_{3}}$ be a sequence of tensors with
$\trank(A_{n})\leq2$ and
\[
\lim_{n\rightarrow\infty}A_{n}=A,
\]
where the limit is taken in any norm topology. If the limiting tensor $A$ has
rank higher than $2$, then $\trank(A)$ must be exactly
$3$ and there exist pairs of linearly independent vectors
$\xx_{1},\yy_{1}\in\Rr^{d_{1}}$, $\xx_{2},\yy_{2}
\in\Rr^{d_{2}}$, $\xx_{3},\yy_{3}\in\Rr^{d_{3}}$
such that
\begin{equation}
A=\xx_{1}\otimes\xx_{2}\otimes\yy_{3}+
\xx_{1}\otimes\yy_{2}\otimes\xx_{3}+\yy_{1}\otimes \xx_{2}\otimes\xx_{3}.
\label{dSLtensor}
\end{equation}
Furthermore, the above result is not vacuous since
\[
{
A_{n}= n \left(  \xx_{1} + \frac{1}{n}\yy_{1}\right)  \otimes\left(
\xx_{2}+\frac{1}{n}\yy_{2}\right)  \otimes \left(  \xx_{3} + \frac{1}{n}\yy_{3}\right)
- n \xx_{1}\otimes\xx_{2}\otimes \xx_{3}}
\]
is an example of a sequence that converges to $A$.
\end{theorem}

A few conclusions can immediately be drawn from
Theorem~\ref{thm:main}: (i) the boundary points of all order-$3$
rank-$2$ tensors can be completely parameterized
by~\eqref{dSLtensor}; (ii) a sequence of order-$3$ rank-$2$ tensors
cannot `jump rank' by more than $1$; (iii) $A$ in \eqref{dSLtensor},
in particular, is an example of a tensor that has no best rank-$2$
approximation.

The formal statements and proofs of the fourth objective appear in
Section~\ref{sect:SA}. The fifth objective is exemplified by our approach
throughout the paper; some specific technical tools are discussed in
Sections \ref{sect:red} and~\ref{subsec:inject}.

On top of these five objectives, we pick up the following smaller results along the way.
Some of these results address frequently asked questions in tensor approximation.
They are discussed in Sections \ref{subsec:divergence}--\ref{subsec:derive} respectively.



\begin{enumerate}
\setcounter{enumi}{\value{goals}}

\item \label{obj:divergence}
\textbf{Divergence of coefficients.}
Whenever a low-rank sequence of tensors converges to a higher-rank tensor, some of the terms in the sequence must blow up. In examples of minimal rank, all the terms blow up.
%

\item \label{obj:nobound}
\textbf{Maximum rank.} For $k \geq 3$, the maximum rank of an order-$k$ tensor in $\Rr^{d_1 \times \dots \times d_k}$ (where $d_i \geq 2$) always exceeds $\min(d_1, \dots, d_k)$. In contrast, for matrices $\min(d_1,d_2)$ does bound the rank.
%

\item \label{obj:gaps}
\textbf{Tensor rank can leap large gaps.}
Conclusion~(ii) in the paragraph above does not generalize to rank~$r>2$.
We will show that a sequence of fixed rank tensors can converge to a limiting tensor of
arbitrarily higher rank.
%

\item \label{obj:bregman}
\textbf{Br\`{e}gman divergences do not help.} If we replace norm by any
continuous measure of `nearness' (including non-metric measures like Br\`{e}gman divergences), it does not change the ill-foundedness of {\brap}.
%

\item \label{obj:leibniz}
\textbf{Leibniz tensors.} We will construct a rich family of sequences of tensors with degenerate limits, labeled by partial derivative operators. The special case $L_3(1)$ is in fact our principal example \eqref{dSLtensor} throughout this paper.
%


\end{enumerate}

%

\subsection{Relation to prior work}

The existence of tensors that can fail to have a best {rank-$r$}
approximation is known to algebraic geometers as early as the 19th
century, albeit in a different language --- the locus of $r$th
secant planes to a Segre variety may not define a (closed) algebraic
variety. It is also known to computational complexity theorists as
the phenomenon underlying the concept of \textit{border rank}
\cite{Bi1, Bi2, BCS, Kn, Land} {and is related to
(but different from) what chemometricians and psychometricians call
`\textsc{candecomp}/\textsc{parafac} degeneracy' {\cite{KDS, KHL, Paa,
Steg1, Steg2}}. We do not claim to be the first to have found such
an example --- the honor belongs to Bini, Capovani, Lotti, and
Romani, who gave an explicit example of a sequence of rank-$5$
tensors converging to a rank-$6$ tensor in 1979 \cite{BCLR}. The
novelty of Theorem~\ref{thm:main} is not in demonstrating that a
tensor may be approximated arbitrarily well by tensors of strictly
lower rank but in \textit{characterizing all such tensors} in the
order-$3$ rank-$2$ case.

Having said this, we would like to point out that the ill-posedness
of the best rank-$r$ approximation problem for high-order tensors
is not at all well-known, as is evident from the paragraph quoted
earlier as well as other discussions in recent publications
\cite{KR1, KR2, K1, K2, ZG}. One likely reason is that in algebraic
geometry, computational complexity, chemometrics, and psychometrics,
the problem is neither stated in the form nor viewed in the light of
obtaining a best low-rank approximation with respect to a choice of
norm (we give several equivalent formulations of {\brap} in
Proposition~\ref{propEquiv}). As such, one goal of this paper will
be to debunk, once and for all, the question of finding best
low-rank approximations for tensors of order~$3$ or higher. As we
stated earlier (as our first and second objectives), our
contribution will be to show that such failures (i)~can and will
occur for tensors of \textit{any} order higher than~$2$, (ii)~that
they will occur for tensors of many different ranks, (iii)~that they
will occur regardless of the choice of norm, and (iv)~that they will
occur with non-zero probability. Formally, we have the following two
theorems (which will appear as Theorems~\ref{thm:rank+1} and
\ref{thm:posvol2} subsequently):

\begin{theorem}
\label{thm:nonclosed}Let $k\geq3$ and $d_{1},\dots,d_{k}\geq2$. For
any $s$ such that $2\leq s\leq\min\{d_{1},\dots,d_{k}\}$, there
exists $A\in \Rr^{d_{1}\times\dots\times d_{k}}$ with $\trank(A)=s$
such that $A$ has no best rank-$r$ approximation for some $r<s$. The
result is independent of the choice of norms.
\end{theorem}

\begin{theorem}
\label{thm:posvol}
If $d_{1}, d_{2}, d_{3} \geq 2$, then the set
\[
\{ A\in\Rr^{d_{1}\times d_{2}\times d_{3}} \mid \mbox{$A$ does not
have a best rank-$2$ approximation}\}
\]
has positive volume; indeed, it contains a nonempty open set.
\end{theorem}



A few features distinguish our work in this paper from existing
studies in algebraic geometry \cite{CGG1, CGG2, Land, LM, Z} and
algebraic computational complexity \cite{AL, AS, Bi1, Bi2, BCLR,
BLR, BCS, S1}: (i) we are interested in tensors over $\Rr$ as
opposed to tensors over $\mathbb{C}$ (it is well-known that the rank
of a tensor is dependent on the underlying field, {cf.\ \eqref{eq:fielddep} and }\cite{Berg}); (ii)
our interest is not limited to order-$3$ tensors (as is often the
case in algebraic computational complexity) --- we would like to prove
results that hold for tensors of any order $k\geq3$; (iii) since we
are interested in questions pertaining to approximations in the
norm, the Euclidean (norm-induced) topology will be more relevant
than the Zariski topology\footnote{Note that the Zariski topology on
$\Bbbk^{n}$ is defined for any field $\Bbbk$ (not just algebraically
closed ones). It is the weakest topology such that all polynomial
functions are continuous. In particularly, the closed sets are
precisely the zero sets of collections of polynomials.} on the
tensor product spaces --- {note in particular that
the claim that a set is not closed in the Euclidean topology is a
stronger statement than the corresponding claim in Zariski topology.}

Our work in this paper in general, and in Section \ref{sec:notsemic} in
particular, is related to studies of `\textsc{candecomp}/\textsc{parafac}
degeneracy' or `diverging \textsc{candecomp}/\textsc{parafac} components'
in psychometrics and chemometrics {\cite{KDS, KHL, Paa, Steg1, Steg2}}.
Diverging coefficients are a necessary consequence of the ill-posedness of
{\brap} (see Propositions~\ref{prop:unbounded} and \ref{prop:unbounded2}).
In fact, examples of `$k$-factor divergence' abound, for arbitrary~$k$ --- see Sections \ref{subsec:higher} and~\ref{subsec:derive} for various constructions.

{Section \ref{sect:condition} discusses how the non-existence of a best rank-$r$ approximation poses serious difficulties for multilinear statistical models based on such approximations. In particular, we will see: (i) why it is meaningless to ask for a `good' rank-$r$ approximation when a best rank-$r$ approximation does not exist; (ii) why even a small perturbation to a rank-$r$ tensor can result in a tensor that has no best rank-$r$ approximation; (iii) why the computational feasibility of finding a `good' rank-$r$ approximation is questionable.}

\subsection{Outline of the paper}

Section~\ref{sect:basics} introduces the basic algebra of tensors and $k$-way arrays. Section~\ref{sect:rankalgebra} defines tensor rank and gives some of its known (and unknown) algebraic properties. Section~\ref{sec:top} studies the topological properties of tensor rank and the phenomenon of rank-jumping. Section \ref{sect:char32} characterizes the problematic tensors in~$\Rr^{2\times2\times2}$, and discusses the implications for approximation problems. Section~\ref{sect:SA} gives a short exposition of the semialgebraic point of view. Section~\ref{sect:orbits} classifies tensors in~$\Rr^{2\times2\times2}$ by orbit type. The orbit structure of tensor spaces is studied from several different aspects. Section~\ref{sec:volume} is devoted to the result that failure of \app$(A,2)$ occurs on a set of positive volume.}

\section{Tensors\label{sect:basics}}

Even though tensors are well-studied objects in the standard graduate
mathematics curriculum \cite{AW, DF, Hun, L, R} and more specifically in
multilinear algebra \cite{B, Gr, M, N, Y}, a `tensor' continues to be viewed
as a mysterious object by outsiders. We feel that we should say a few words to
demystify the term.

In mathematics, the question `what is a vector?' has the simple answer
`a vector is an element of a vector space' --- in other words, a vector
is characterized by the axioms that define the algebraic operations on a
vector space. In physics, however, the question `what is a vector?'
often means `what kinds of physical quantities can be represented by
vectors?' The criterion has to do with the change of basis theorem: an
$n$-dimensional vector is an `object' that is represented by $n$ real
numbers once a basis is chosen \textit{only} if those real numbers
transform themselves as expected when one changes the basis. For exactly
the same reason, the meaning of a tensor is obscured by its more
restrictive use in physics. In physics (and also engineering), a tensor is
an `object' represented by a $k$-way array of real numbers that transforms
according to certain rules (cf.\ \eqref{tran}) under a change of basis. In
mathematics, these `transformation rules' are simply consequences of the
multilinearity of the tensor product and the change of basis theorem for
vectors. Nowadays, books written primarily for a physics audience
\cite{Ge, Mar} have increasingly adopted the mathematical definition, but a handful of recently published books continue to propagate the obsolete (and vague) definition.
%
To add to the confusion, `tensor' is frequently used to refer to
a tensor field (e.g.\ metric tensor, stress tensor, Riemann curvature
tensor).

For our purposes, an order-$k$ \textit{tensor} $\mathbf{A}$ is simply an element
of a \textit{tensor product} of $k$ real vector spaces, $V_{1}\otimes V_{2}%
\otimes\dots\otimes V_{k}$, as defined in any standard algebra textbook
\cite{AW, B, DF, Gr, Hun, L, M, N, R, Y}. Up to a choice of bases on
$V_{1},\dots,V_{k}$, such an element may be \textit{coordinatized}, i.e.\ represented as a $k$-way array~$A$ of
real numbers --- much as an element of an $n$-dimensional vector space may, up to a choice of basis, be represented by an $n$-tuple of numbers in $\Rr^{n}$. We will
let $\Rr^{d_{1}\times\dots\times d_{k}}$ denote the vector space of
$k$-way arrays of real numbers $A=\llbracket a_{j_{1}\cdots j_{k}}%
\rrbracket_{j_{1}=1,\dots,j_{k}=1}^{d_{1},\dots,d_{k}}$ with addition and
scalar multiplication defined coordinatewise:%
\begin{equation}\label{eq:vecsp}
\llbracket a_{j_{1}\cdots j_{k}}\rrbracket+\llbracket b_{j_{1}\cdots j_{k}%
}\rrbracket:=\llbracket a_{j_{1}\cdots j_{k}}+b_{j_{1}\cdots j_{k}%
}\rrbracket\quad\text{and}\quad\lambda\llbracket a_{j_{1}\cdots j_{k}%
}\rrbracket:=\llbracket\lambda a_{j_{1}\cdots j_{k}}\rrbracket.
\end{equation}
A $k$-way arrays of numbers (or $k$-array) is also sometimes referred to as a $k$-dimensional \textit{hypermatrix} \cite{GKZ1}.

It may be helpful to think of a $k$-array as a \textit{data structure}, convenient for representing or storing the coefficients of a tensor with respect to a set of bases. The tensor itself carries with it an \textit{algebraic structure}, by virtue of being an element of a tensor product of vector spaces. Once bases have been chosen for these vector spaces, we may view the order-$k$ tensor as a $k$-way array equipped with the algebraic operations defined in \eqref{eq:vecsp} and \eqref{eq:mmm}. Despite this correspondence, it is not wise to regard `tensor' as being synonymous with `array'.

{\it Notation.} We will denote elements of abstract tensor spaces in boldface upper-case letters; whereas $k$-arrays will be denoted in italic upper-case letters. Thus $\mathbf{A}$ is an abstract tensor, which may be represented by an array of numbers~$A$ with respect to a basis. We will use double brackets to enclose the entries of a $k$-array --- $A=\llbracket a_{j_{1}\cdots j_{k}}\rrbracket_{j_{1}=1,\dots,j_{k}=1}^{d_{1},\dots,d_{k}}$ --- and when there is no risk of confusion, we will leave out the range of the indices and simply write $A=\llbracket
a_{j_{1}\cdots j_{k}}\rrbracket$.



\subsection{Multilinear matrix multiplication\label{sect:MMM}}

Matrices can act on other matrices through two independent
multiplication operations: left-multiplication and
right-multiplication. Matrices act on order-$3$ tensors via
\emph{three} different multiplication operations. These can be
combined into a single formula. If $A=\llbracket
a_{ijk}\rrbracket\in\Rr^{d_{1}\times d_{2}\times d_{3}}$ and
$L=[\lambda_{pi}]\in\Rr^{c_{1}\times d_{1}}$, $M=[\mu_{qj}%
]\in\Rr^{c_{2}\times d_{2}}$, $N=[\nu_{rk}]\in\Rr^{c_{3}\times
d_{3}}$, then the array~$A$ may be transformed into an array $A^{\prime
}=\llbracket
a_{pqr}^{\prime}\rrbracket\in\Rr ^{c_{1}\times c_{2}\times c_{3}}$,
by the equation:
\begin{equation}
a_{pqr}^{\prime}= \sum\nolimits_{i,j,k=1}^{d_{1},d_{2},d_{3}} \lambda_{pi}%
\mu_{qj}\nu_{rk}a_{ijk}\, \label{tran}%
\end{equation}
We call this operation the \emph{multilinear multiplication} of~$A$ by
matrices $L,M,N$, which we write succinctly as
\[
A^{\prime}= (L,M,N)\cdot A.
\]
Informally, we are multiplying the $3$-way array $A$ on its three `sides' by the
matrices $L,M,N$ respectively.

\textit{Remark.} This notation is standard in mathematics --- the
elements of a product $G_{1}\times G_{2} \times G_{3}$ are generally
grouped in the form $(L,M,N)$, and when a set with some algebraic
structure~$G$ acts on another set~$X$, the result of $g\in G$ acting
on $x\in X$ is almost universally written $g\cdot x$ \cite{AW, B,
DF, Hun, L, R}. Here we are just looking at the case $G =
\Rr^{c_{1}\times d_{1}} \times\Rr^{c_{2}\times d_{2}}
\times\Rr^{c_{3}\times d_{3}} $ and $X=\Rr^{d_{1}\times d_{2}\times
d_{3}}$. This is consistent with notation adopted in earlier work
\cite{J1} but more recent publications such
as \cite{LMV1, LMV2} have used $A\times_{1}L^{\transp}\times_{2}M^{\transp}%
\times_{3}N^{\transp}$ in place of $(L,M,N)\cdot A$.

Multilinear matrix multiplication extends in a straightforward way to arrays
of arbitrary order: if $A=\llbracket a_{i_{1}\cdots i_{k}}\rrbracket\in
\Rr^{d_{1}\times\dots\times d_{k}}$ and $L_{1}=[\lambda_{ij}^{(1)}%
]\in\Rr^{c_{1}\times d_{1}},\dots,L_{k}=[\lambda_{ij}^{(k)}%
]\in\Rr^{c_{k}\times d_{k}}$, then $A^{\prime}=(L_{1},\dots,L_{k})\cdot
A$ is the array $A^{\prime}=\llbracket
a_{i_{1}\cdots i_{k}}^{\prime}\rrbracket\in\Rr^{c_{1}\times\dots\times
c_{k}}$ given by%
\begin{equation}\label{eq:mmm}
a_{i_{1}\cdots i_{k}}^{\prime}=\sum\nolimits_{i_{1},\dots,i_{k}=1}%
^{d_{1},\dots,d_{k}}\lambda_{i_{1}j_{1}}\cdots\lambda_{i_{k}j_{k}}%
a_{j_{1}\cdots j_{k}}\,.
\end{equation}

We will now see how a $3$-way array representing a tensor in $V_{1}\otimes
V_{2}\otimes V_{3}$ transforms under changes of bases of the vector
spaces~$V_{i}$. Suppose the $3$-way array $A=\llbracket a_{ijk}\rrbracket\in
\Rr^{d_{1}\times d_{2}\times d_{3}}$ represents an order-$3$ tensor
$\mathbf{A}\in V_{1}\otimes V_{2}\otimes V_{3}$ with respect to bases
$\mathcal{B}_{1}=\{\ee_{i}\mid i=1,\dots,d_{1}\}$, $\mathcal{B}%
_{2}=\{\ff_{j}\mid j=1,\dots,d_{2}\}$, $\mathcal{B}_{3}=\{\mathbf{g}%
_{k}\mid k=1,\dots,d_{3}\}$ on $V_{1},V_{2},V_{3}$, i.e.
\begin{equation}
\mathbf{A}= \sum\nolimits_{i,j,k=1}^{d_{1},d_{2},d_{3}} a_{ijk}\ee%
_{i}\otimes\ff_{j}\otimes\mathbf{g}_{k}. \label{basis}%
\end{equation}
Suppose we choose different bases, $\mathcal{B}_{1}^{\prime}=\{\ee%
_{i}^{\prime}\mid i=1,\dots,d_{1}\}$, $\mathcal{B}_{2}^{\prime}=\{\ff%
_{j}^{\prime}\mid j=1,\dots,d_{2}\}$, $\mathcal{B}_{3}^{\prime}=\{\mathbf{g}%
_{k}^{\prime}\mid k=1,\dots,d_{3}\}$ on $V_{1},V_{2},V_{3}$ where
\begin{equation}
\ee_{i}= \sum\nolimits_{p=1}^{d_{1}} \lambda_{ip}\ee_{p}%
^{\prime}, \quad\ff_{j}=\sum\nolimits_{q=1}^{d_{2}}\mu_{jq}%
\ff_{q}^{\prime}, \quad\mathbf{g}_{k}=\sum\nolimits_{r=1}^{d_{3}}%
\nu_{kr}\mathbf{g}_{r}^{\prime}, \label{cob}%
\end{equation}
and $L=[\lambda_{pi}]\in\Rr^{d_{1}\times d_{1}}$, $M=[\mu_{qj}%
]\in\Rr^{d_{2}\times d_{2}}$, $N=[\nu_{rk}]\in\Rr^{d_{3}\times
d_{3}}$ are the respective change-of-basis matrices. Substituting the
expressions for \eqref{cob} into~\eqref{basis}, we get
\[
\mathbf{A}=\sum\nolimits_{p,q,r=1}^{d_{1},d_{2},d_{3}} a_{pqr}^{\prime
}\ee_{p}^{\prime}\otimes\ff_{q}^{\prime}\otimes\mathbf{g}%
_{r}^{\prime}%
\]
where
\begin{equation}
a_{pqr}^{\prime}= \sum\nolimits_{i,j,k=1}^{d_{1},d_{2},d_{3}} \lambda_{pi}
\mu_{qj}\nu_{rk}a_{ijk}\,, \label{tran1}%
\end{equation}
or more simply $A^{\prime}= (L,M,N)\cdot A$, where the $3$-way array
$A^{\prime}=\llbracket a_{pqr}^{\prime}\rrbracket\in\Rr^{d_{1}\times
d_{2}\times d_{3}}$ represents $\mathbf{A}$ with respect to this new choice of
bases $\mathcal{B}_{1}^{\prime},\mathcal{B}_{2}^{\prime},\mathcal{B}%
_{3}^{\prime}$ .

All of this extends immediately to {order-$k$} tensors and
{$k$-way} arrays. Henceforth, when a choice of basis is implicit, we will
not distinguish between an order-$k$ tensor and the $k$-way array that
represents it.

The change-of-basis matrices $L,M,N$ in the discussion above are of course invertible; in other words they belong their respective \emph{general linear} groups. We write
$\GL_{d}(\Rr)$ for the group of nonsingular matrices
in~$\Rr^{d \times d}$. Thus $L \in\GL_{d_{1}}%
(\Rr)$, $M \in\GL_{d_{2}}(\Rr)$, $N \in
\GL_{d_{3}}(\Rr)$. In addition to general linear
transformations, it is natural to consider orthogonal transformations. We
write $\Oo_{d}(\Rr)$ for the subgroup of
$\GL_{d}(\Rr)$ of transformations which preserve the
Euclidean inner product. The following shorthand is helpful:
\begin{align*}
\GL_{d_{1},\dots,d_{k}}(\Rr)  &  :=\GL%
_{d_{1}}(\Rr)\times\dots\times\GL_{d_{k}}(\Rr)\\
\Oo_{d_{1},\dots,d_{k}}(\Rr)  &  :=\Oo%
_{d_{1}}(\Rr)\times\dots\times\Oo_{d_{k}}(\Rr)
\end{align*}
Then $\Oo_{d_{1},\dots,d_{k}}(\Rr) \leq\GL%
_{d_{1},\dots,d_{k}}(\Rr)$, and both groups act on~$\Rr%
^{d_{1}\times\dots\times d_{k}}$ via multilinear multiplication.

\begin{definition}
\label{def:equivalence} Two tensors
$A,A^{\prime}\in\Rr^{d_{1}\times\dots\times d_{k}}$ are said to be
$\GL$-equivalent (or simply `equivalent') if there exists
$(L_{1},\dots,L_{k})\in\GL_{d_{1},\dots,d_{k}}(\Rr)$ such that
$A^{\prime}= (L_{1},\dots,L_{k})\cdot A$. More strongly, we say that
$A, A^{\prime}$ are $\Oo$-equivalent if such a transformation~$L$
can be found in~$\Oo_{d_{1},\dots,d_{k}}(\Rr)$.
\end{definition}

For example, if $V_{1},\dots,V_{k}$ are vector spaces and $\dim(V_{i})=d_{i}$, then
$A,A^{\prime}\in\Rr^{d_{1}\times\dots\times d_{k}}$ represent the same
tensor in $V_{1}\otimes\dots\otimes V_{k}$ with respect to two different bases
iff $A, A'$ are $\GL$-equivalent.

We finish with some trivial properties of multilinear matrix multiplication: for $A,B\in\Rr%
^{d_{1}\times\dots\times d_{k}}$, and $\alpha,\beta\in\Rr$,%
\begin{equation}
(L_{1},\dots,L_{k})\cdot(\alpha A+\beta B)=\alpha(L_{1},\dots,L_{k})\cdot
A+\beta(L_{1},\dots,L_{k})\cdot B \label{mmult1}%
\end{equation}
and for $L_{i}\in\Rr^{c_{i}\times d_{i}}$, $M_{i}\in\Rr%
^{b_{i}\times c_{i}}$, $i=1,\dots,k$,%
\begin{equation}
(M_{1},\dots,M_{k})\cdot\lbrack(L_{1},\dots,L_{k})\cdot A]=(M_{1}L_{1}%
,\dots,M_{k}L_{k})\cdot A. \label{mmult2}%
\end{equation}
Lastly, the name \textit{multilinear} matrix multiplication is justified since
for any $M_{i},N_{i}\in\Rr^{c_{i}\times d_{i}}$, $\alpha,\beta
\in\Rr$,
\begin{align}
(L_{1},\dots,\alpha M_{i}+\beta N_{i},\dots,L_{k})\cdot A  &  =\alpha
(L_{1},\dots,M_{i},\dots,L_{k})\cdot A\\
&  \qquad+\beta(L_{1},\dots,N_{i},\dots,L_{k})\cdot A.\nonumber
\end{align}


\subsection{Outer-product rank and outer-product decomposition of a tensor}

Let $\Rr^{d_{1}}\otimes\dots\otimes\Rr^{d_{k}}$ be the tensor
product of the vector spaces $\Rr^{d_{1}},\dots,\Rr^{d_{k}}$.
Note that the \textit{Segre map}%
\begin{equation}
\Rr^{d_{1}}\times\dots\times\Rr^{d_{k}}\rightarrow
\Rr^{d_{1}\times\dots\times d_{k}},\qquad(\xx_{1},\dots
,\xx_{k})\mapsto\bigl\llbracket x_{j_{1}}^{(1)}\cdots x_{j_{k}}%
^{(k)}\bigr\rrbracket_{j_{1},\dots,j_{k}=1}^{d_{1},\dots,d_{k}} \label{segre}%
\end{equation}
is multilinear and so by the universal property of tensor product \cite{AW, B,
DF, Gr, Hun, L, M, N, R, Y}, we have a unique linear map $\varphi$ such that
the following diagram commutes:
\[
\begin{diagram}[labelstyle=\scriptstyle]
&  & \Rr^{n_{1}}\otimes\dots\otimes\Rr^{n_{k}}\\
& \ruTo & \dTo^{\theta}\\
\Rr^{n_{1}}\times\dots\times\Rr^{n_{k}} &
\rTo^{\varphi} & \Rr^{n_{1}\times\dots\times n_{k}}
\end{diagram}
\]
Clearly,%
\begin{equation}
\varphi(\xx_{1}\otimes\dots\otimes\xx_{k}%
)=\bigl\llbracket x_{j_{1}}^{(1)}\cdots x_{j_{k}}^{(k)}\bigr\rrbracket_{j_{1}%
,\dots,j_{k}=1}^{d_{1},\dots,d_{k}} \label{op0}%
\end{equation}
and $\varphi$ is a vector space isomorphism since $\dim(\Rr%
^{d_{1}\times\dots\times d_{k}})=\dim(\Rr^{d_{1}}\otimes\dots
\otimes\Rr^{d_{k}})=d_{1}\cdots d_{k}$. Henceforth we will not
distinguish between these two spaces. The elements of $\Rr^{d_{1}%
}\otimes\dots\otimes\Rr^{d_{k}}\cong\Rr^{d_{1}\times\dots\times
d_{k}}$ will be called a tensor and we will also drop $\varphi$ in \eqref{op0}
and write%
\begin{equation}
\xx_{1}\otimes\dots\otimes\xx_{k}=\bigl\llbracket x_{j_{1}%
}^{(1)}\cdots x_{j_{k}}^{(k)}\bigr\rrbracket_{j_{1},\dots,j_{k}=1}%
^{d_{1},\dots,d_{k}}. \label{op}%
\end{equation}
Note that the symbol $\otimes$ in \eqref{op0} denotes the formal tensor
product and by dropping $\varphi$, we are using the same symbol $\otimes$ to
define the \textit{outer product} of the vectors $\xx_{1}%
,\dots,\xx_{k}$ via the formula \eqref{op}. Hence, a tensor can be
represented either as a $k$-dimensional array or as a sum of formal tensor
products of $k$ vectors --- where the equivalence between these two objects is
established by taking the formal tensor product of $k$~vectors as defining a
$k$-way array via \eqref{op}.

It is clear that the map in \eqref{segre} is not surjective --- the image
consists precisely of the \textit{decomposable tensors}: a tensor
$A\in\Rr^{d_{1}\times\dots\times d_{k}}$ is said to be
\textit{decomposable} if it can be written in the form
\[
A=\xx_{1}\otimes\dots\otimes\xx_{k}%
\]
with $\xx_{i}\in\Rr^{d_{i}}$ for $i=1,\dots,k$. It is easy to
see that multilinear matrix multiplication of decomposable tensors obeys the
formula:%
\begin{equation}\label{mmdecomp}
(L_{1},\dots,L_{k}) \cdot (\xx_{1}\otimes\dots\otimes\xx_{k}%
)=L_{1}\xx_{1}\otimes\dots\otimes L_{k}\xx_{k}.
\end{equation}

{\it Remark.} The outer product can be viewed as a special case of multilinear matrix multiplication. For example, a linear combination of outer products of vectors may be expressed in terms of multilinear matrix multiplication:
\[
\sum\nolimits_{i=1}^{r}
\lambda_i \xx_{i}\otimes\yy_{i}\otimes\mathbf{z}_{i} = (X,Y,Z)\cdot \Lambda
\]
with matrices $X = [\xx_{1},\dots, \xx_{r}] \in \mathbb{R}^{l \times r}$, $Y = [\yy_{1},\dots, \yy_{r}] \in \mathbb{R}^{m \times r}$, $Z = [\zz_{1},\dots, \zz_{r}] \in \mathbb{R}^{n \times r}$ and a `diagonal tensor' $\Lambda = \operatorname{diag} [\lambda_1, \dots, \lambda_r] \in \mathbb{R}^{r \times r \times r}$.

We now come to the main concept of interest in this paper.

\begin{definition}
A tensor has \textbf{outer-product rank}~$r$ if it can be written as a sum of
$r$~decomposable tensors, but no fewer. We will write $\trank(A)$ for the outer-product rank of $A$. So%
\[
\trank(A):=\min\{r\mid A=
{\textstyle\sum\nolimits_{i=1}^{r}}
\uu_{i}\otimes\vv_{i}\otimes\dots\otimes\mathbf{z}_{i}\}.
\]
\end{definition}%
Note that a non-zero decomposable tensor has outer-product rank~$1$.
%

Despite several claims of originality as well as many misplaced attributions to these claims, the concepts of tensor rank and the decomposition of a tensor into a sum of outer-products of vectors was the product of much earlier work by Frank~L.~Hitchcock in 1927 \cite{Hi1, Hi2}. We call this the outer-product rank mainly to distinguish it from the \textit{multilinear rank} to be defined in Section~\ref{sectMRMD} (due to Hitchcock) but we will use the term \textit{rank} or \textit{tensor rank} most of the time when there is no danger of confusion.

\begin{lemma}[Invariance of tensor rank]
\label{lemRankMM}
{\rm (1)} If $A \in \Rr^{d_{1}\times\dots\times d_{k}}$ and
$(L_{1}, \dots, L_{k}) \in
\Rr^{c_1 \times d_1} \times \dots \times \Rr^{c_k \times d_k}$,
then
\begin{equation}\label{leqtrank}
\trank((L_{1},\dots,L_{k}) \cdot A) \leq \trank(A).
\end{equation}
{\rm (2)} If $A\in\Rr^{d_{1}\times\dots\times d_{k}}$ and
$(L_{1},\dots,L_{k})\in\GL_{d_{1},\dots,d_{k}}(\Rr%
):=\GL_{d_{1}}(\Rr)\times\dots\times\GL%
_{d_{k}}(\Rr)$, then
\begin{equation}\label{invtrank}
\trank((L_{1},\dots,L_{k}) \cdot A)=\trank(A).
\end{equation}
\end{lemma}

\begin{proof}
\eqref{leqtrank} follows from \eqref{mmdecomp} and \eqref{mmult1}. Indeed, if
$
A=\sum_{j=1}^{r} \xx_{1}^{j}\otimes\dots\otimes\xx_{k}^{j}
$
then
$
(L_1, \dots, L_k)\cdot A=
\sum_{j=1}^{r} L_1 \xx_{1}^{j}\otimes\dots\otimes L_k \xx_{k}^{j}
$.
%
Furthermore, if the~$L_i$ are invertible then by
\eqref{mmult2} we get%
\[
A=(L_{1}^{-1},\dots,L_{k}^{-1})\cdot\lbrack(L_{1},\dots,L_{k})\cdot A \rbrack
\]
and so $\trank(A)\leq \trank((L_{1},\dots,L_{k})\cdot A)$, hence~\eqref{invtrank}.
\end{proof}

\subsection{The outer product and direct sum operations on tensors}

The outer product of vectors defined earlier is a special case of the
\textit{outer product} of two tensors. Let $A\in\Rr^{d_{1}\times
\dots\times d_{k}}$ be a tensor of order $k$ and $B\in\Rr^{c_{1}%
\times\dots\times c_{\ell}}$ be a tensor of order $\ell$, then the outer
product of $A$ and $B$ is the tensor $C:=A\otimes B\in\Rr^{d_{1}%
\times\dots\times d_{k}\times c_{1}\times\dots\times c_{l}}$ of order $k+\ell$
defined by%
\[
c_{i_{1}\cdots i_{k}j_{1}\cdots j_{l}}=a_{i_{1}\cdots i_{k}}b_{j_{1}\cdots
j_{l}}.
\]

The \textit{direct sum} of two order-$k$ tensors $A\in\Rr^{d_{1}%
\times\dots\times d_{k}}$ and $B\in\Rr^{c_{1}\times\dots\times c_{k}}$
is the order-$k$ tensor $C:=A\oplus B\in\Rr^{(c_{1}+d_{1})\times
\dots\times(c_{k}+d_{k})}$ defined by%
\[
c_{i_{1},\dots,i_{k}}=%
\begin{cases}
a_{i_{1},\dots,i_{k}} & \text{if }1\leq i_{\alpha}\leq d_{\alpha}\text{,
}\alpha=1,\dots,k;\\
b_{i_{1}-d_{1},\dots,i_{k}-d_{k}} & \text{if }d_{\alpha}+1\leq i_{\alpha}\leq
c_{\alpha}+d_{\alpha}\text{, }\alpha=1,\dots,k;\\
0 & \text{otherwise.}%
\end{cases}
\]
For matrices, the direct sum of $A\in\Rr^{m_{1}\times n_{1}}$ and
$B\in\Rr^{m_{2}\times n_{2}}$ is simply the block-diagonal matrix%
\[
A\oplus B=%
\begin{bmatrix}
A & 0\\
0 & B
\end{bmatrix}
\in\Rr^{(m_{1}+m_{2})\times(n_{1}+n_{2})}.
\]
The direct sum of two order-$3$ tensors $A\in\Rr^{l_{1}\times
m_{1}\times n_{1}}$ and $B\in\Rr^{l_{2}\times m_{2}\times n_{2}}$ is a
`block tensor' with $A$ in the $(1,1,1)$-block and $B$ in the $(2,2,2)$-block
\[
A\oplus B=\left[
\begin{array}
[c]{rr}%
A & 0\\
0 & 0
\end{array}
\right\vert \!\left.
\begin{array}
[c]{rr}%
0 & 0\\
0 & B
\end{array}
\right]  \in\Rr^{(l_{1}+l_{2})\times(m_{1}+m_{2})\times(n_{1}+n_{2})}.
\]

In abstract terms, if $U_i, V_i, W_i$ are vector spaces such that
$W_i = U_i \oplus V_i$ for $i = 1, \dots, k$, then tensors $A \in
U_1 \otimes \dots \otimes U_k$ and $B \in V_1 \otimes \dots \otimes
V_k$ have direct sum $A \oplus B \in W_1 \otimes \dots \otimes W_k$.

\subsection{Tensor subspaces}

Whenever $c\leq d$ there is a canonical embedding $\Rr^{c}%
\subseteq\Rr^{d}$ given by identifying the $c$~coordinates of
$\Rr^{c}$ with the first~$c$ coordinates of $\Rr^{d}$.

Let $c_{i}\leq d_{i}$ for $i=1,\dots,k$. Then there is then a canonical embedding
%
$\Rr^{c_{1}\times\dots\times c_{k}}\subset\Rr^{d_{1}\times\dots\times d_{k}}$, defined as the tensor product of the embeddings
$\Rr^{c_{i}}\subseteq\Rr^{d_{i}}$. We say that $\Rr^{c_{1}\times\dots\times c_{k}}$ is a \textit{tensor subspace} of
$\Rr^{d_{1}\times\dots\times d_{k}}$.
More generally, if $U_i, V_i$ are vector spaces with $U_i \subset V_i$ for $i = 1, \dots, k$, then there is an inclusion $U_1 \otimes \dots \otimes U_k \subset V_1 \otimes \dots \otimes V_k$ defined as the tensor product of the inclusions $U_i \subset V_i$. Again we say that $U_1 \otimes \dots \otimes U_k$ is a tensor subspace of $V_1 \otimes \dots \otimes V_k$.

If $B \in \Rr^{c_{1}\times\dots\times c_{k}}$ then its image under the canonical embedding into $\Rr^{d_{1}\times\dots d_{k}}$ can be written in the form $B \oplus 0$, where $0 \in \Rr^{{(d_1-c_1)}\times\dots\times{(d_k-c_k)}}$ is the zero tensor.
%
A tensor $A \in\Rr^{d_{1}\times\dots\times d_{k}}$ is said to be \emph{$\GL$-equivalent} (or simply `equivalent')
to~$B$
%
%
if there exists $(L_{1},\dots,L_{k})\in\GL_{d_{1},\dots,d_{k}}(\Rr)$
%
such that $B \oplus 0= (L_{1},\dots,L_{k})\cdot A$.
%
More strongly, we say that $A$ is \emph{$\Oo$-equivalent} (`orthogonally equivalent')
to~$B$ if such a transformation can be found
in~$\Oo_{d_{1},\dots,d_{k}}(\Rr)$.

We note that $A$ is $\GL$-equivalent to~$B$ if and only if there exist {full-rank} matrices $M_i \in \Rr^{d_i\times c_i}$ such that $A = (M_1, \cdots, M_k) \cdot B$.
In one direction, $M_i$~can be obtained as the first~$c_i$ columns of~$L_i^{-1}$. In the other direction, $L_i^{-1}$ can be obtained from~$M_i$ by adjoining extra columns. There is a similar statement for $\Oo$-equivalence. Instead of full rank, the condition is that the matrices~$M_i$ have orthogonal columns.

An important simplifying principle in tensor algebra is that questions about a tensor --- such as `what is its rank?' --- can sometimes, as we shall see, be reduced to analogous questions about an equivalent tensor in a lower-dimensional tensor subspace.

\subsection{Multilinear rank and multilinear decomposition of a
tensor\label{sectMRMD}}

Although we focus on outer product rank in this paper, there is a simpler
notion of multilinear rank which directly generalizes the column and row ranks
of a matrix to higher order tensors.

For convenience, we will consider order-$3$ tensors only. Let $A =
\llbracket a_{ijk}\rrbracket \in\Rr^{d_{1}\times d_{2}\times d_{3}}$.
For fixed values of $j \in\{1, \dots, d_{2}\}$ and $k \in\{ 1, \dots, d_{3}%
\}$, consider the vector $A_{\bullet jk} := [a_{ijk}]_{i=1}^{d_{1}}
\in\Rr^{d_{1}}$. Likewise consider (column) vectors $A_{i \bullet k}
:= [a_{ijk}]_{j=1}^{d_{2}} \in\Rr^{d_{2}}$ for fixed values of
$i,k$, and (row) vectors $A_{ij \bullet} := [a_{ijk}]_{k=1}^{d_{3}}
\in\Rr^{d_{3}}$ for fixed values of $i,j$. In analogy with row rank
and column rank, define
\begin{align*}
r_{1}(A)  &  := \dim(\Span_{\Rr}\{ A_{\bullet
jk} \mid1 \leq j \leq d_{2}, 1 \leq k \leq d_{3} \}),\\
r_{2}(A)  &  := \dim(\Span_{\Rr}\{ A_{i \bullet
k} \mid1 \leq i \leq d_{1}, 1 \leq k \leq d_{3} \}),\\
r_{3}(A)  &  := \dim(\Span_{\Rr}\{ A_{ij \mathbb{\bullet
}} \mid1 \leq i \leq d_{1}, 1 \leq j \leq d_{2} \}).
\end{align*}

For another interpretation, note that $\Rr^{d_{1}\times d_{2}\times
d_{3}}$ can be viewed as $\Rr^{d_{1}\times d_{2}d_{3}}$ by ignoring the
multiplicative structure between the second and third factors. Then $r_{1}(A)$
is simply the rank of~$A$ regarded as $d_{1}\times d_{2}d_{3}$ matrix.
There are similar definitions for $r_{2}(A)$ and $r_{3}(A)$.

The \textit{multilinear rank} of $A$, denoted\footnote{The symbol
$\boxplus$ is meant to evoke an impression of the rows and columns in a
matrix.} $\multirank(A)$, is the $3$-tuple
$(r_{1}(A),r_{2}(A),r_{3}(A))$. Again, this concept is not new but
was first explored by Hitchcock under the name \textit{multiplex
rank} in the same papers where he defined tensor rank \cite{Hi1,
Hi2}. What we term multilinear rank will be equivalent to
Hitchcock's duplex rank. A point to note is that $r_{1}(A)$,
$r_{2}(A)$, $r_{3}(A)$, and $\trank(A)$ are in general all different
--- a departure from the case of matrices, where the row rank,
column rank and outer product rank are always equal. Observe that
we will always have
\begin{equation}
\label{mrtr}
r_{i}(A)\le \min\{\trank(A),d_{i}\}.
\end{equation}
Let us verify this for~$r_1$: if $A = \xx_1 \otimes \yy_1 \otimes \zz_1 + \dots + \xx_r \otimes \yy_r \otimes \zz_r$ then each vector $A_{\bullet j k}$ belongs to $\Span(\xx_1, \dots, \xx_r)$. This implies that $r_1 \leq \trank(A)$, and $r_1 \leq d_1$ is immediate from the definitions.
%
%
A simple but useful consequence of~\eqref{mrtr} is that
\begin{equation}
\label{trmr}
\trank(A) \ge \lVert \multirank(A) \rVert_\infty = \max \{r_i(A) \mid i = 1,\dots,k \}.
\end{equation}

If $A\in\Rr^{d_{1}\times d_{2}\times d_{3}}$ then and $\multirank(A) =
(r_{1},r_{2},r_{3})$, then there exist subspaces $U_i \subset
\Rr^{d_i}$ with $\dim(U_i) = r_i$, such that $A \in U_1 \otimes U_2 \otimes U_3$.
We call these the \emph{supporting subspaces} of~$A$. The supporting subspaces are minimal, in the sense that if $A \in V_1 \otimes V_2 \otimes V_3$ then $U_i \subset V_i$ for $i=1,2,3$. This observation leads to an alternate definition:
\[
{r_i(A) = \min \{ \dim(U_i)
\mid
U_1 \subset \Rr^{d_1}
,\,
U_2 \subset \Rr^{d_2}
,\,
U_3 \subset \Rr^{d_3}
,\,
A \in U_1 \otimes U_2 \otimes U_3
\}.}
\]

An immediate  consequence of this characterization is that
{$\multirank(A)$ is invariant under the action of $\GL_{d_1,d_2,d_3}(\Rr)$: if $A' =
(L,M,N) \cdot A$, where $(L,M,N) \in \GL_{d_1,d_2,d_3}(\Rr)$, then $\multirank(A) = \multirank((L,M,N)\cdot A)$.
Indeed, if $U_1, U_2, U_3$ are the supporting subspaces of~$A$, then $L(U_1)$, $M(U_2)$, $N(U_3)$ are the supporting subspaces of $(L,M,N)\cdot A$.

More generally, we have multilinear rank equivalents of \eqref{leqtrank} and~\eqref{invtrank}: if $A \in \mathbb{R}^{d_1 \times \dots \times d_k}$ and
$(L_1,\dots, L_k) \in \Rr^{c_1 \times d_1} \times \dots \times \Rr^{c_k \times d_k}$,
then
\begin{equation}\label{leqmrank}
\multirank((L_1,\dots,L_k)\cdot A) \leq \multirank(A),
\end{equation}
and if $A \in \mathbb{R}^{d_1 \times \dots \times d_k}$ and $(L_1,\dots, L_k) \in \GL_{d_1,\dots,d_k}(\Rr)$, then
\begin{equation}\label{invmrank}
\multirank((L_1,\dots,L_k)\cdot A)
=
\multirank(A).
\end{equation}

Suppose $\multirank(A) = (r_1, r_2, r_3)$. By applying transformations~$L_i \in \GL_{d_i}(\Rr)$ which carry $U_i$ to $\Rr^{r_i}$, it follows that $A$ is equivalent to some $B \in \Rr^{r_1 \times r_2 \times r_3}$.
Alternatively there exist $B\in\Rr^{r_{1}\times r_{2}\times r_{3}}$ and full-rank matrices
$L \in\Rr^{d_{1}\times r_{1}}$, $M \in\Rr^{d_{2}\times r_{2}}$,
$N \in\Rr^{d_{3}\times r_{3}}$, such that%
\[
A=(L,M,N)\cdot B.
\]
%
%
The matrices $L,M,N$ may be chosen to have orthonormal columns or to be unit
lower-triangular --- a fact easily deduced from applying the $QR$%
-decomposition or the $LU$-decomposition to the full-rank matrices $L,M,N$ and
using \eqref{mmult2}.


To a large extent, the study of tensors~$A\in\Rr^{d_{1}\times d_{2}\times d_{3}}$ with $\multirank(A) \leq(r_{1}, r_{2}, r_{3})$ reduces to the study of tensors in~$\Rr^{r_{1} \times r_{2} \times r_{3}}$. This is a useful reduction, but (unlike the matrix case) it does not even come close to giving us a full classification of tensor types.

\subsection{Multilinear orthogonal projection}
\label{subsec:multiprojection}

If $U$ is a subspace of an inner-product space~$V$ (for instance, $V = \Rr^n$ with the usual dot product), then there is an orthogonal projection from $V$ onto~$U$, which we denote $\pi_U$. We regard this as a map $V \to V$. As such, it is self-adjoint (i.e.\ has a symmetric matrix with respect to any orthonormal basis), and satisfies $\pi_U^2 = \pi_U$, $\operatorname{im}(\pi_U) = U$, $\ker(\pi_U) = U^\perp$.
We note Pythagoras' theorem for any $\vv \in V$:
\[
\| \vv \|^2 = \| \pi_U \vv \|^2 + \| (1 - \pi_U)\vv \|^2
\]

We now consider orthogonal projections for tensor spaces. If $U_1, U_2, U_3$ are subspaces of $V_1, V_2, V_3$, respectively, then $U_1 \otimes U_2 \otimes U_3$ is a tensor subspace of $V_1 \otimes V_2 \otimes V_3$, and the multilinear map $\Pi = (\pi_{U_1}, \pi_{U_2}, \pi_{U_3})$ is a projection onto that subspace. In fact, $\Pi$ is orthogonal with respect to the Frobenius norm. The easiest way to see this is to identify $U_i \subset V_i$ with $\Rr^{c_i} \subset \Rr^{d_i}$ by taking suitable orthonormal bases; then $\Pi$ acts by zeroing out all the entries of a $d_1 \times d_2 \times d_3$ array outside the initial $c_1 \times c_2 \times c_3$ block. In particular we have Pythagoras' theorem for any $A \in V_1 \otimes V_2 \otimes V_3$:
\begin{equation}
\label{eq:pythagoras}
\| A \|_F^2 = \| \Pi A \|_F^2 + \| (1 - \Pi)A \|_F^2
\end{equation}
Being a multilinear map, $\Pi$ is non-increasing for $\trank, \multirank$, as in \eqref{leqtrank}, \eqref{leqmrank}.

There is a useful orthogonal projection $\Pi_A$ associated with any tensor $A \in \Rr^{d_1 \times d_2 \times d_3}$. Let $U_1, U_2, U_3$ be the supporting subspaces of~$A$, so that $A \in U_1 \otimes U_2 \otimes U_3$, and $\dim(U_i) = r_i(A)$ for $i = 1,2,3$.
%
%
Define:
\[
\Pi_A = (\pi_1(A), \pi_2(A), \pi_3(A)) = (\pi_{U_1}, \pi_{U_2}, \pi_{U_3})
%
\]

\begin{proposition}
\label{prop:PiA}
$\Pi_A(A) = A$.
\end{proposition}

\begin{proof}
$A$ belongs to $U_1 \otimes U_2 \otimes U_3$, which is fixed by~$\Pi_A$.
\end{proof}

\begin{proposition}
\label{prop:PiAconst}
The function $A \mapsto \Pi_A$ is continuous over subsets of $\Rr^{d_1
\times d_2 \times d_3}$ on which $\multirank(A)$ is constant.
\end{proposition}

\begin{proof}
We show, for example, that $\pi_{1} = \pi_1(A)$ depends continuously on~$A$. For any $A \in \Rr^{d_1\times d_2\times d_3}$, select $r = r_1(A)$ index pairs $(j,k)$ such that the vectors $A_{\bullet j k}$ are linearly independent. For any $B$ near~$A$, assemble the marked vectors as a matrix $X = X(B) \in \Rr^{d_i \times r}$. Then $\pi_1 = X(X^{\transp}X)^{-1}X^{\transp} =: P(B)$ by a well-known formula in linear algebra. The function $P(B)$ is defined and continuous as long as the $r$ selected vectors remain independent, which is true on a neighborhood of~$A$. Finally, the orthogonal projection defined by $P(B)$ maps onto the span of the $r$~selected vectors. Thus, if $r_1(B) = r$ then $P(B) = \pi_1(B)$.
\end{proof}

It is clear that the results of this section apply to tensor spaces of all orders.


\section{The algebra of tensor rank}
\label{sect:rankalgebra}

We will state and prove a few basic results about the outer-product rank.

\begin{proposition}
\label{prop:rank}Let $A\in\Rr^{c_{1}\times\dots\times c_{k}}%
\subset\Rr^{d_{1}\times\dots\times d_{k}}$. The rank of~$A$ regarded as
a tensor in~$\Rr^{c_{1}\times\dots\times c_{k}}$ is the same as the
rank of~$A$ regarded as a tensor in~$\Rr^{d_{1}\times\dots\times d_{k}%
}$.
\end{proposition}

\begin{proof}
For each~$i$ the identity on~$\Rr^{c_{i}}$ factors as a pair of maps
$\Rr^{c_{i}}\overset{\iota_i}{\hookrightarrow}\Rr^{d_{i}}%
\overset{\pi_i}{\twoheadrightarrow}\Rr^{c_{i}}$, where $\iota_i$ is the
canonical inclusion and $\pi$~is the projection map given by deleting the last
$d_{i}-c_{i}$ coordinates.
Applying \eqref{leqtrank} twice,
we have
\begin{eqnarray*}
\trank(A)
\geq
 \trank((\iota_{1},\dots,\iota_{k})\cdot A)
&\geq&
 \trank((\pi_{1}, \dots, \pi_{k}) \cdot (\iota_{1},\dots,\iota_{k})\cdot A)
\\
&=&
 \trank((\pi_{1} \iota_{1}, \dots, \pi_{k} \iota_{k})\cdot A)
\\
&=& \trank(A)
\end{eqnarray*}
so $A \in \Rr^{c_1 \times \dots \times c_k}$ and its image $(\iota_{1}, \dots , \iota_{k})\cdot A \in \Rr^{d_1 \times \dots \times d_k}$ must have equal tensor ranks.
%
%
\end{proof}


\begin{corollary}
\label{cor:reduction}
Suppose $A \in \Rr^{d_1 \times \dots \times d_k}$ and $\multirank(A) \leq (c_1, \dots, c_k)$. Then $\trank(A) = \trank(B)$ for an equivalent tensor $B \in \Rr^{c_1\times\dots\times c_k}$.
\endproof
\end{corollary}

The next corollary asserts that tensor rank is consistent under a different
scenario: when order~$k$ tensors are regarded as order~$l$ tensors, for $l >
k$, by taking the tensor product with a non-zero monomial term.

\begin{corollary}
\label{cor:rank}Let $A\in\Rr^{d_{1}\times\dots\times d_{k}}$ be an
order-$k$ tensor and $\uu_{k+1}\in\Rr^{d_{k+1}}, \dots,$
$\uu_{k+\ell}\in\Rr^{d_{k+\ell}}$ be non-zero vectors. Then%
\[
\trank(A)=\trank(A\otimes
\uu_{k+1}\otimes\dots\otimes\uu_{k+\ell}).
\]

\end{corollary}

\begin{proof}
Let $c_{k+1}=\dots=c_{k+l}=1$ and apply Proposition~\ref{prop:rank} to
$A\in\Rr^{d_{1}\times\dots\times d_{k}}=\Rr^{d_{1}\times
\dots\times d_{k}\times c_{k+1}\times\dots\times c_{k+\ell}}\hookrightarrow
\Rr^{d_{1}\times\dots\times d_{k}\times d_{k+1}\times\dots\times
d_{k+\ell}}$. Note that the image of the inclusion is $A\otimes\ee%
_{1}^{(k+1)}\otimes\dots\otimes\ee_{1}^{(k+\ell)}$ where $\ee%
_{1}^{(i)}=(1,0,\dots,0)^{\transp}\in\Rr^{d_{i}}$. So we have%
\[
\trank(A\otimes\ee_{1}^{(k+1)}\otimes
\dots\otimes\ee_{1}^{(k+\ell)})=\trank(A).
\]
The general case for arbitrary non-zero $\uu_{i}\in\Rr^{d_{i}}$
follows from applying to $A\otimes\ee_{1}^{(k+1)}\otimes\dots
\otimes\ee_{1}^{(k+\ell)}$ a multilinear multiplication $(I_{d_{1}%
},\dots,I_{d_{k}},L_{1},\dots,L_{\ell})\in\GL_{d_{1}%
,\dots,d_{k+\ell}}(\Rr)$ where $I_{d}$ is the $d\times d$ identity
matrix and $L_{i}$ is a non-singular matrix with $L_{i}\ee%
_{i}=\uu_{i}$. It then follows from Lemma~\ref{lemRankMM} that%
\begin{multline*}
\trank(A\otimes\uu_{k+1}\otimes\dots
\otimes\uu_{k+\ell})\\
\begin{aligned} & = \trank[(I_{d_{1}},\dots,I_{d_{k}},L_{1},\dots,L_{\ell})\cdot(A\otimes\ee_{1}^{(k+1)}\otimes\dots\otimes\ee_{1}^{(k+\ell)})]\\ & = \trank(A\otimes\ee_{1}^{(k+1)}\otimes\dots\otimes\ee_{1}^{(k+\ell)}). \end{aligned}
\end{multline*}

\end{proof}

Let $E=\uu_{k+1}\otimes\uu_{k+2}\otimes\dots\otimes
\uu_{k+\ell}\in\Rr^{d_{k+1}\times\dots\times d_{k+\ell}}$. So
$\trank(E)=1$ and Corollary~\ref{cor:rank} says that
$\trank(A\otimes E)=\trank(A)\trank(E)$. Note that this last relation does not
generalize. If $\trank(A)>1$ and $\trank(B)>1$, then it is true that
\[
\trank(A \otimes B) \leq \trank(A) \trank(B),
\]
since one can multiply decompositions of $A,B$ term by term to obtain a decomposition of~$A \otimes B$, but it can happen (cf.\ \cite{BCS})} that%
\[
\trank(A\otimes B) < \trank(A) \trank(B).
\]

The corresponding statement for direct sum is still an open problem for tensors of order~$3$ or higher. It has been conjectured by Strassen \cite{S1} that
\begin{equation}
\trank(A\oplus B)\overset{?}{=}\trank(A)+\trank(B)
\label{DSC}
\end{equation}
for all order-$k$ tensors $A$ and~$B$. However J\'aJ\'a and Takche \cite{JT}
have shown that for the special case when $A$ and~$B$ are of order~$3$ and at
least one of them is a matrix pencil (i.e.\ a tensor of size $p\times q\times
2$, $p\times2\times q$, or $2\times p\times q$ that may be regarded as a pair
of $p\times q$ matrices), then the direct sum conjecture holds.

\begin{theorem}[J\'aJ\'a--Takche \cite{JT}]
\label{jaja}
Let $A\in\Rr^{c_{1}\times c_{2}\times c_{3}}$ and $B\in\Rr^{d_{1}\times d_{2}\times d_{3}}$. If $2\in\{c_{1},c_{2},c_{3},d_{1},d_{2},d_{3}\}$, then
\[
\trank(A\oplus B)=\trank(A)+\trank(B).
\]
\endproof
\end{theorem}

It is not hard to define tensors of arbitrarily high rank so long as we have sufficiently many linearly independent vectors in every factor.

\begin{lemma}\label{buildrank}
For $\ell=1,\dots,k$, let
$\xx_{1}^{(\ell)},\dots,\xx_{r}^{(\ell)}\in\Rr^{d_{i}}$ be linearly
independent.
Then the tensor defined by%
\[
A:=\sum_{j=1}^{r}\xx_{j}^{(1)}\otimes\xx_{j}^{(2)}\otimes
\dots\otimes\xx_{j}^{(k)}\in\Rr^{d_{1}\times d_{2}\times
\dots\times d_{k}}%
\]
has $\trank(A)=r$.
\end{lemma}

\begin{proof}
Note that $\multirank(A) = (r,r,...,r)$. By \eqref{trmr}, we get
\[
\trank(A) \ge \max\{r_i(A) \mid i =1,\dots,k \} = r.
\]
On the other hand, it is clear that $\trank(A) \le r$.
\end{proof}

Thus, in $\Rr^{d_{1}\times\dots\times d_{k}}$, it is easy to write
down tensors of any rank~$r$ in the range $0\leq
r\leq\min\{d_{1},\dots,d_{k}\}$. For matrices, this exhausts all
possibilities; the rank of $A\in \Rr^{d_{1}\times d_{2}}$ is at most
$\min\{d_{1},d_{2}\}$. In contrast, for $k\geq3$, there will always
be tensors in $\mathbb{R}^{d_1 \times d_k}$ that have rank
exceeding $\min\{d_{1},\dots,d_{k}\}$. We will see this in
Theorem~\ref{thm:rank+1}.

\section{The topology of tensor rank}
\label{sec:top} 

Let $A=\llbracket a_{i_{1}\cdots i_{k}}\rrbracket\in\Rr^{d_{1}%
\times\dots\times d_{k}}$. The \textit{Frobenius norm} of $A$
and its associated inner product are defined by
%
\[
\lVert A\rVert_{F}^{2}:=
\sum\nolimits_{i_{1},\dots,i_{k}=1}^{d_{1},\dots,d_{k}}
 \lvert a_{i_{1}\cdots i_{k}}\rvert^{2},
\quad
\langle A, B \rangle_F :=
\sum\nolimits_{i_{1},\dots,i_{k}=1}^{d_{1},\dots,d_{k}}
 a_{i_{1}\cdots i_{k}}  b_{i_{1}\cdots i_{k}}
.
\]

Note that for a decomposable tensor, the Frobenius norm satisfies%
\begin{equation}
\lVert\uu\otimes\vv\otimes\dots\otimes\mathbf{z}\rVert
_{F}=\lVert\uu\rVert_{2}\lVert\vv\rVert_{2}\cdots\lVert
\mathbf{z}\rVert_{2} \label{Fnorm2}%
\end{equation}
where $\lVert\cdot\rVert_{2}$ denotes the $l^{2}$-norm of a vector, and more generally
\begin{equation}
\label{eqn:Fnorm3}
\| A \otimes B \|_F = \|A\|_F \| B \|_F
\end{equation}
for arbitrary tensors $A,B$.
Another important property which follows from \eqref{mmdecomp} and \eqref{Fnorm2} is orthogonal invariance:
\begin{equation}
\label{eqn:Fnorm4}
\lVert (L_1, \dots, L_k) \cdot A \rVert_{F}
= \lVert A \rVert_{F}
\end{equation}
whenever $(L_1, \dots, L_k) \in \Oo_{d_1, \dots, d_k}(\Rr)$.
There are of course many other natural choices of norms on tensor product
spaces \cite{DF1, Gro}. The important thing to note is that $\Rr^{d_1
\times \dots \times d_k}$ being finite dimensional, all these norms will
induce the same topology.

We define the following (topological) subspaces of $\Rr^{d_{1}%
\times\dots\times d_{k}}$.
\begin{align*}
\Ss_{r}(d_{1},\dots,d_{k})  &  =\left\{  A\in\Rr^{d_{1}%
\times\dots\times d_{k}}\mid\trank(A)\leq r\right\} \\
\Ssbar_{r}(d_{1},\dots,d_{k})  &  =\text{closure of
$\Ss_{r}(d_{1},\dots,d_{k})\subset\Rr^{d_{1}\times\dots\times
d_{k}}$}%
\end{align*}
Clearly the only reason to define $\Ssbar_r$ is the
sad fact that $\Ss_r$ is not necessarily (or even usually)
closed --- the theme of this paper. See Section~\ref{sec:notsemic}.

We occasionally refer to elements of $\Ss_{r}$ as `rank-$r$ tensors'.
This is slightly inaccurate, since lower-rank tensors are included, but
convenient. However, the direct assertions `$A$ has
rank~$r$' and `$\operatorname*{rank}(A)=r$' are always meant in the precise sense.
The same remarks apply to `border rank', which is defined in Section~\ref{sect:weak}.
We refer to elements of ${\Ssbar}_{r}$ as `border-rank-$r$ tensors', and describe them as being `rank-$r$-approximable'.

Theorem~\ref{thm:b} asserts that $\Ssbar_{2}(d_1,d_2,d_3)\subset
\Ss_{3}(d_1,d_2,d_3)$ for all $d_1,d_2,d_3$, and that the
exceptional tensors $\overline{\Ss_{2}}(d_1,d_2,d_3) \setminus
{\Ss}_{2}(d_1,d_2,d_3)$ are all of a particular form.

\subsection{Upper semicontinuity}
\label{sec:uppersemic}

Discrete-valued rank functions on spaces of matrices or tensors cannot
usefully be continuous, because they would then be constant and would not have
any classifying power. As a sort of compromise, matrix rank is well known to
be an upper semicontinuous function; if $\operatorname*{rank}(A)=r$ then
$\operatorname*{rank}(B)\geq r$ for all matrices~$B$ in a neighborhood of~$A$.
This is not true for the outer-product rank of tensors (as we will see Section~\ref{sec:notsemic}). There are several equivalent ways of formulating this assertion.

\begin{proposition}
\label{propEquiv}Let $r\geq2$ and $k\geq3$. Given the norm-topology on
$\Rr^{d_{1}\times\dots\times d_{k}}$, the following statements are equivalent:

\begin{enumerate}
[\upshape (a)]

\item The set $\Ss_{r}(d_{1},\dots,d_{k}):=\{A\in\Rr%
^{d_{1}\times\dots\times d_{k}}\mid\trank(A)\leq r\}$
is not closed.

\item There exists a sequence $A_{n}\in\Rr^{d_{1}\times\dots\times
d_{k}}$, $\trank(A_{n})\leq r$, $n\in\mathbb{N}$,
converging to $B\in\Rr^{d_{1}\times\dots\times d_{k}}$ with
$\trank(B)>r$.

\item There exists $B\in\Rr^{d_{1}\times\dots\times d_{k}}$,
$\trank(B)>r$, that may be approximated arbitrarily
closely by tensors of strictly lower rank, i.e.%
\[
\inf\{\lVert B-A\rVert\mid\trank(A)\leq r\}=0.
\]

\item There exists $C\in\Rr^{d_{1}\times\dots\times d_{k}}$,
$\trank(C)>r$, that does not have a best rank-$r$
approximation, i.e.%
\[
\inf\{\lVert C-A\rVert\mid\trank(A)\leq r\}
\]
is not attained (by any $A$ with $\trank(A)\leq r$).
\end{enumerate}
\end{proposition}

\begin{proof}
It is obvious that (a) $\Rightarrow$ (b) $\Rightarrow$ (c) $\Rightarrow$ (d).
To complete the chain, we just need to show that (d) $\Rightarrow$ (a).
Suppose $\Ss:=\Ss_{r}(d_{1},\dots,d_{k})$ is closed. Since the
closed ball of radius $\lVert C\rVert$ centered at $C$, $\{A\in\Rr%
^{d_{1}\times\dots\times d_{k}}\mid\lVert C-A\rVert\leq\lVert C\rVert\}$,
intersects $\Ss$ non-trivially (e.g.\ $0$ is in both sets). Their
intersection $\mathcal{T}=\{A\in\Rr^{d_{1}\times\dots\times d_{k}}%
\mid\trank(A)\leq r,\lVert C-A\rVert\leq\lVert
C\rVert\}$ is a non-empty compact set. Now observe that%
\[
\delta:=\inf\{\lVert C-A\rVert\mid A\in\Ss\}=\inf\{\lVert
C-A\rVert\mid A\in\mathcal{T}\}
\]
since any $A^{\prime}\in\Ss\backslash\mathcal{T}$ must have $\lVert
C-A^{\prime}\rVert>\lVert C\rVert$ while we know that $\delta\leq\lVert
C\rVert$. By the compactness of $\mathcal{T}$, there exists $A_{\ast}%
\in\mathcal{T}$ such that $\lVert C-A_{\ast}\rVert=\delta$. So the required
infimum is attained by $A_{\ast}\in\mathcal{T}\subset\Ss$.
\end{proof}

We caution the reader that there exist tensors of rank~$>r$ that do
not have a best rank-$r$ approximation but \textit{cannot} be
approximated arbitrarily closely by rank-$r$ tensors, i.e.\
$\inf\{\lVert C-A\rVert\mid \trank(A)\leq r\}>0$. In other words,
statement~(d) applies to a strictly larger class of tensors than
statement~(c) (cf.\ Section~\ref{sec:volume}). The tensors in
statement~(d) are sometimes called `degenerate' in the psychometrics
and chemometrics literature (e.g.\ {\cite{KDS, KHL, Paa, Steg1, Steg2}})
but we prefer to
avoid this term since it is inconsistent (and often at odds) with
common usage in Mathematics. For example, in Table
\ref{table:orbits}, the tensors in the orbit classes of $D_{2},
D_{2}^{\prime}, D_{2}^{\prime\prime}$ are all degenerate but
statement~(d) does not apply to them; on the other hand, the tensors
in the orbit class of $G_3$ are non-degenerate but Theorem
\ref{thm:norank2} tells us that they are all of the form in
statement~(d).


We begin by getting three well-behaved cases out of the way. The proofs shed light on what can go wrong in all the other cases.

\begin{proposition}
\label{prop:goodranka}
For all $d_1, \dots, d_k$, we have $\Ssbar_1(d_1,\dots,d_k) = \Ss_1(d_1,\dots,d_k)$.
\end{proposition}

\begin{proof}
Suppose $A_n \to A$ where $\trank(A_n) \leq 1$. We can write
\[
A_n =
\lambda_n \uu_{1,n} \otimes \uu_{2,n} \otimes \dots \otimes \uu_{k,n}
\]
where $\lambda_n = \| A_n \|$ and the vectors $\uu_{i,n} \in
\Rr^{d_i}$ have unit norm. Certainly $\lambda_n = \|A_n\| \to \|A\|
=: \lambda$. Moreover, since the unit sphere in~$\Rr^{d_i}$ is
compact, each sequence $\uu_{i,n}$ has a convergent
subsequence, with limit $\uu_i$, say. It follows that there is
a subsequence of~$A_n$ which converges to $\lambda\uu_1\otimes
\dots \otimes\uu_k$. This must equal~$A$, and it has rank at
most~1.
\end{proof}

\begin{proposition}
\label{prop:goodrankb}
For all $r$ and $d_1, d_2$, we have $\Ssbar_r(d_1,d_2) = \Ss_r
(d_1,d_2)$. In other words, matrix rank is upper-semicontinuous.
\end{proposition}

\begin{proof}
Suppose $A_n \to A$ where $\rank(A_n) \leq r$, so we can write
\[
A_n = \lambda_{1,n} \uu_{1,n}\otimes \vv_{1,n} + \dots +
\lambda_{r,n} \uu_{r,n} \otimes \vv_{r,n}.
\]
Convergence of the sequence~$A_n$ does not imply convergence of the individual terms $\lambda_{i,n}$, $\uu_{i,n}$, $\vv_{i,n}$, even in a subsequence. However, if we take the singular value decomposition, then the $\uu_{i,n}$ and~$\vv_{i,n}$ are unit vectors and the $\lambda_{i,n}$ satisfy
\[
\lambda_{1,n}^2 + \dots + \lambda_{r,n}^2 = \| A_n \|
\]
Since $\|A_n\| \to \|A\|$ this implies that the $\lambda_{i,n}$ are uniformly bounded. Thus we can find a subsequence with convergence $\lambda_{i,n} \to \lambda_i$, $\uu_{i,n} \to \uu_i$, $\vv_{i,n} \to \vv_i$ for all~$i$. Then
\[
A = \lambda_1 \uu_1 \otimes \vv_1 + \dots +
\lambda_r \uu_r \otimes \vv_r
\]
which has rank at most~$r$.
\end{proof}

\begin{proposition}
\label{prop:goodrankc}
The multilinear rank function $\multirank(A) = (r_1(A), \dots, r_k(A))$ is upper semicontinuous.
\end{proposition}

\begin{proof}
Each $r_i$ is the rank of a matrix obtained by rearranging the
entries of~$A$, and is therefore upper semicontinuous in~$A$ by
Proposition \ref{prop:goodrankb}.
\end{proof}

\begin{corollary}
\label{cor:goodapprox} Every tensor has a best rank-$1$
approximation. Every matrix has a best rank-$r$ approximation. Every
order-$k$ tensor has a best approximation with $\multirank \leq
(r_1, \dots, r_k)$, for any specified $(r_1, \dots, r_k)$.
\end{corollary}

\begin{proof}
These statements follow from Proposition \ref{prop:goodranka}, \ref{prop:goodrankb} and~\ref{prop:goodrankc}, together with the implication (d)$\Rightarrow$(a) from
Proposition~\ref{propEquiv}.
\end{proof}

\subsection{Tensor rank is not upper semicontinuous}
\label{sec:notsemic}

Here is the simplest example of the failure of outer-product rank to be upper
semicontinuous.
This is the first example of a more general construction which we discuss in Section~\ref{subsec:derive}.
A formula similar to \eqref{eqdSL} appeared as Exercise 62 in Section
4.6.4 of Knuth's {\it The Art of Computer Programming}~\cite{Kn} (the original source is unknown to us but may well be
\cite{Kn}).
Other examples have appeared in \cite{BCLR} (the earliest known to us) and~\cite{Paa}, as well as in unpublished work of Kruskal.

\begin{proposition}
\label{prop:dSLtensor} Let $\xx_{1},\yy_{1}\in\Rr^{d_{1}%
}$, $\xx_{2},\yy_{2}\in\Rr^{d_{2}}$ and $\xx%
_{3},\yy_{3}\in\Rr^{d_{3}}$ be vectors such that each pair
$\xx_{i},\yy_{i}$ is linearly independent. Then the tensor%
\begin{equation}
A:=\xx_{1}\otimes\xx_{2}\otimes\yy_{3}+\xx%
_{1}\otimes\yy_{2}\otimes\xx_{3}+\yy_{1}\otimes
\xx_{2}\otimes\xx_{3}\in\Rr^{d_{1}\times d_{2}\times
d_{3}} \label{eqdSL}%
\end{equation}
has rank~$3$ but can be approximated arbitrarily closely by tensors of
rank~$2$. In particular, $A$ does not have a best rank-$2$ approximation.
\end{proposition}

\begin{proof}
For each $n\in\mathbb{N}$, define
\begin{equation}
A_{n}:= n \left(  \xx_{1} + \frac{1}{n}\yy_{1}\right)  \otimes\left(
\xx_{2}+\frac{1}{n}\yy_{2}\right)  \otimes \left(  \xx_{3} + \frac{1}{n}\yy_{3}\right)
- n \xx_{1}\otimes\xx_{2}\otimes \xx_{3}
\label{dSLseq}
\end{equation}
Clearly, $\trank(A_{n})\leq2$, and since, as $n\rightarrow
\infty$,
\begin{multline*}
\lVert A_{n}-A\rVert_{F} \le
\frac{1}{n}\lVert \yy_{1}\otimes\yy_{2}\otimes\xx_{3}+\yy_{1}\otimes\xx_{2}\otimes\yy_{3}+\xx_{1}\otimes\yy_{2}\otimes\yy_{3}\rVert_{F}\\
+ \frac{1}{n^2}\lVert \yy_{1}\otimes\yy_{2}\otimes\yy_{3}\rVert_{F}
\rightarrow 0,
\end{multline*}
we see that $A$ is approximated arbitrary closely by tensors $A_{n}$.
%

It remains to establish that $\trank(A) = 3$. From the
three-term format of~$A$, we deduce only that $\trank(A) \leq3$. A clean proof that $\trank(A)>2$ is
included in the proof of Theorem~\ref{thm:orbits}, but this depends on the
properties of the polynomial~$\Delta$ defined in Section~\ref{subsec:delta}. A
more direct argument is given in the next lemma.
\end{proof}

\begin{lemma}
\label{lem:dSLtensor-rank}
Let $\xx_{1},\yy_{1}\in\Rr^{d_{1}}$, $\xx_{2},\yy_{2}\in\Rr^{d_{2}}$,
$\xx_{3},\yy_{3}\in\Rr^{d_{3}}$ and
\[
A =
\xx_{1}\otimes\xx_{2}\otimes\yy_{3} +
\xx_{1}\otimes\yy_{2}\otimes\xx_{3} +
\yy_{1}\otimes \xx_{2}\otimes\xx_{3}.
\]
Then $\trank(A)=3$ if and only if $\xx_{i},\yy_{i}$ are linearly independent for $i=1,2,3$.
\end{lemma}

\begin{proof}
Only two distinct vectors are involved in each factor of the tensor product, so $\multirank(A) \leq (2,2,2)$ and we can work in~$\Rr^{2\times 2\times 2}$ (Corollary~\ref{cor:reduction}). More strongly, if any of the pairs $\{ \xx_{i}, \yy_{i}\}$ is
linearly dependent, then $A$ is $\GL$-equivalent to a tensor in $\Rr^{1 \times2 \times2}$, $\Rr^{2\times1 \times2}$ or $\Rr^{2 \times2 \times1}$. These spaces are
isomorphic to $\Rr^{2 \times2}$, so the maximum possible rank of~$A$ is~2.

Conversely, suppose each pair $\{\xx_{i}, \yy_{i}\}$ is linearly
independent. We may as well assume that
\begin{equation}
A =
\left[
\begin{array}[c]{rr} 0 & 1\\ 1 & 0 \end{array}
\right\vert \!\left.
\begin{array}[c]{rr} 1 & 0\\ 0 & 0 \end{array}
\right]
\label{eq:dSL2}
\end{equation}
since we can transform~$A$ to that form using a multilinear transformation
$(L_{1}, L_{2}, L_{3})$ where $L_{i}(\xx_{i}) = \ee_{1}$ and
$L_{i}(\yy_{i}) = \ee_{2}$ for $i=1,2,3$.

Suppose, for a contradiction, that $\trank(A)\leq2$; then we can write
\begin{equation}
A= \uu_{1}\otimes\uu_{2}\otimes\uu_{3}
+ \vv_{1}\otimes\vv_{2}\otimes\vv_{3}
\label{eq:dSL3}
\end{equation}
for some $\uu_{i},\vv_{i}\in\Rr^{d_{i}}$.

\textrm{Claim~1}: The vectors $\uu_{1},\vv_{1}$ are independent.
If they are not, then let $\varphi:\Rr^{2}\rightarrow\Rr$ be a
nonzero linear map such that $\varphi(\uu_{1})=\varphi(\vv%
_{1})=0$. Using the expressions in \eqref{eq:dSL3} and~\eqref{eq:dSL2}, we
find that
\[
\mathbf{0}=(\varphi,I,I)\cdot A=\left[
\begin{array}
[c]{cc}%
\varphi(\ee_{2}) & \varphi(\ee_{1})\\
\varphi(\ee_{1}) & 0
\end{array}
\right]
\]
in $\Rr^{1\times2\times2}\cong\Rr^{2\times2}$, which is a
contradiction because $\varphi(\ee_{1})$ and $\varphi(\ee_{2})$
cannot both be zero.

\textrm{Claim~2}: The vectors $\uu_{1},\ee_{1}$ are dependent.
Indeed, let $\varphi_{u}:\Rr^{2}\rightarrow\Rr$ be a linear map
whose kernel is spanned by~$\uu_{1}$. Then
\[
\varphi_{u}(\vv_{1})(\vv_{2}\otimes\vv_{3})=(\varphi
_{u},I,I)\cdot A=\left[
\begin{array}
[c]{cc}%
\varphi_{u}(\ee_{2}) & \varphi_{u}(\ee_{1})\\
\varphi_{u}(\ee_{1}) & 0
\end{array}
\right]
\]
in $\Rr^{1\times2\times2}\cong\Rr^{2\times2}$. The
\textsc{lhs}\ has rank at most~1, which implies on the \textsc{rhs}\ that
$\varphi_{u}(\ee_{1})=0$, and hence $\ee_{1}\in
\Span\{\uu_{1}\}$.

\textrm{Claim~3}: The vectors $\vv_{1},\ee_{1}$ are dependent.
Indeed, let $\varphi_{v}:\Rr^{2}\rightarrow\Rr$ be a linear map
whose kernel is spanned by~$\vv_{1}$. Then
\[
\varphi_{v}(\uu_{1})(\uu_{2}\otimes\uu_{3})=(\varphi
_{v},I,I)\cdot A=\left[
\begin{array}
[c]{cc}%
\varphi_{v}(\ee_{2}) & \varphi_{v}(\ee_{1})\\
\varphi_{v}(\ee_{1}) & 0
\end{array}
\right]
\]
in $\Rr^{1\times2\times2}\cong\Rr^{2\times2}$. The
\textsc{lhs}\ has rank at most~1, which implies on the \textsc{rhs}\ that
$\varphi_{v}(\ee_{1})=0$, and hence $\ee_{1}\in
\Span\{\vv_{1}\}$.

Taken together, the three claims are inconsistent. This is the desired
contradiction. Thus $\trank(A) > 2$ and therefore
$\trank(A) = 3$.
\end{proof}

{\it Remark.} Note that if we take $d_{1}=d_{2}=d_{3}=2$, then \eqref{eqdSL} is an
example of a tensor whose outer product rank exceeds $\min\{d_{1},d_{2},d_{3}\}$.

\subsection{Diverging coefficients}
\label{subsec:divergence}

What goes wrong in the example of Proposition~\ref{prop:dSLtensor}? Why do the rank-2 decompositions of the~$A_n$ fail to converge to a rank-2 decomposition of~$A$? We can attempt to
mimic the proofs of Propositions \ref{prop:goodranka} and~\ref{prop:goodrankb} by seeking convergent subsequences for the rank-$2$ decompositions of the~$A_n$. We fail because we cannot simultaneously keep all the variables bounded. For example, in the decomposition
\[
A_{n} =
n \left(  \xx_{1} + \frac{1}{n}\yy_{1}\right)  \otimes\left(
\xx_{2}+\frac{1}{n}\yy_{2}\right)  \otimes \left(  \xx_{3} + \frac{1}{n}\yy_{3}\right)
- n \xx_{1}\otimes\xx_{2}\otimes \xx_{3}
\]
the vector terms converge but the coefficients $\lambda_1 = \lambda_2 = n$ tend to infinity. In spite of this, the sequence~$A_n$ itself remains bounded.

In fact, rank-jumping always occurs like this (see also \cite{KDS}).

\begin{proposition}
\label{prop:unbounded}
Suppose $A_n \to A$, where $\trank(A) \geq r+1$ and $\trank(A_n) \leq r$ for all~$n$. If we write
\[
A_n = \lambda_{1,n} \uu_{1,n}\otimes \vv_{1,n}\otimes\ww_{1,n} + \dots +
\lambda_{r,n} \uu_{r,n} \otimes \vv_{r,n}\otimes \ww_{r,n},
\]
where the vectors $\uu_{i,n}$, $\vv_{i,n}$, $\ww_{i,n}$ are unit vectors, then $\max_i \{|\lambda_{i,n}| \} \to \infty$ as $n \to \infty$. Moreover, at least two of the coefficient sequences $\{ \lambda_{i,n} \mid n = 1, 2, \dots \}$ are unbounded.
\end{proposition}

\begin{proof}
If the sequence $\max_i \{|\lambda_{i,n}| \}$ does not diverge to $\infty$, then it has a bounded subsequence. In this subsequence, the coefficients and vectors are all bounded, so we can pass to a further subsequence in which each of the coefficient sequences and vector sequences is convergent:
\[
\lambda_{i,n} \to \lambda_i, \quad
\uu_{i,n} \to \uu_i, \quad
\vv_{i,n} \to \vv_i, \quad
\ww_{i,n} \to \ww_i
\]
It follows that $A = \lambda_1 \uu_1 \otimes \vv_1 \otimes \ww_1 + \dots +
\lambda_r \uu_r \otimes \vv_r \otimes \ww_r$, so it has rank at most~$r$, which is a contradiction.

Thus $\max_i\{|\lambda_{i,n}|\}$ diverges to~$\infty$. It follows that at least one of the coefficient sequences has a divergent subsequence. If there were only one such coefficient sequence, all the others being bounded, then (on the subsequence) $A_n$ would be dominated by this term and consequently $\|A_n\|$ would be unbounded. Since $A_n \to A$, this cannot happen. Thus there are at least two unbounded coefficient sequences.
\end{proof}

For a \emph{minimal} rank-jumping example, \emph{all} the coefficients must diverge to~$\infty$.
\begin{proposition}
\label{prop:unbounded2}
Suppose $A_n \to A$, where $\trank(A) = r+s$ and $\trank(A_n) \leq r$ for all~$n$. If we write
\[
A_n = \lambda_{1,n} \uu_{1,n}\otimes \vv_{1,n}\otimes\ww_{1,n} + \dots +
\lambda_{r,n} \uu_{r,n} \otimes \vv_{r,n}\otimes \ww_{r,n},
\]
where the vectors $\uu_{i,n}$, $\vv_{i,n}$, $\ww_{i,n}$ are unit vectors, then there are two possibilities: {either} (i) all of the sequences $|\lambda_{i,n}|$ diverge to~$\infty$ as $n \to \infty$; {or} (ii) in the same tensor space there exists $B_n \to B$, where $\trank(B) \geq r'+s$ and $\trank(B_n) \leq r'$ for all~$n$, for some $r' < r$.
\end{proposition}

\begin{proof}
Suppose one of the coefficient sequences, say $|\lambda_{i,n}|$, fails to diverge as $n\to \infty$; so it has a bounded subsequence. In a further subsequence, the $i$th term ${R}_n = \lambda_{i,n}\uu_{i,n}\otimes\vv_{i,n}\otimes\ww_{i,n}$ converges to a tensor ${R}$ of rank (at most)~1. Writing $B_n = A_n - R_n$, we find that $B_n \to B = A - R$ on this subsequence, with $\trank(B_n) \leq r-1$. Moreover, $r+s \leq \trank(A) \leq \trank(B) + \trank(R)$, so $\trank(B) \geq (r-1)+s$.
\end{proof}

{\it Remark.} Clearly the arguments in Propositions \ref{prop:unbounded} and~\ref{prop:unbounded2} apply to tensors of all orders, not just order~3. We also note that the vectors ($\uu_{i,n}$ etc.) need not be unit vectors; they just have to be uniformly bounded.

One interpretation of Proposition~\ref{prop:unbounded} is that if one attempts to minimize
\[
\lVert A - \lambda_1 \uu_1\otimes \vv_1 \otimes\ww_1 - \dots - \lambda_r
\uu_r \otimes \vv_r \otimes \ww_r \rVert
\]
for a tensor $A$ which does not have a best rank-$r$ approximation, then (at least some of) the coefficients $\lambda_i$ become unbounded.
This phenomenon of diverging summands has been observed in practical
applications of multilinear models in psychometrics and chemometrics and
is commonly referred to in those circles as
`\textsc{candecomp}/\textsc{parafac} degeneracy' or `diverging
\textsc{candecomp}/\textsc{parafac} components' \cite{KDS, KHL, Paa,
Steg1, Steg2}. More precisely, these are called `$k$-factor degeneracies'
when there are $k$ diverging summands whose sum stays bounded. $2$- and
$3$-factor degeneracies were exhibited in~\cite{Paa} and $4$-
and $5$-factor degeneracies were exhibited in~\cite{Steg1}. There are uninteresting (see Section~\ref{subsec:higher}) and interesting (see Section~\ref{subsec:derive}) ways of generating $k$-factor degeneracies for arbitrarily large~$k$.

\subsection{Higher orders, higher ranks, arbitrary norms}
\label{subsec:higher}

We will now show that the rank-jumping phenomenon --- that is, the failure of $\Ss_{r}(d_{1},\dots,d_{k})$ to be closed --- is independent of the choice of norms and can be extended to arbitrary order.
The norm independence is a trivial consequence of a basic fact in
functional analysis: all norms on finite dimensional
vector spaces are equivalent; in particular, any norm will induce
the same unique topology on a finite dimensional vector space.

\begin{theorem}\label{thm:rank+1}
For $k \ge 3$, and $d_1 ,\dots,d_k \ge 2$, the problem of
determining a best rank-$r$ approximation for an order-$k$ tensor in
$\Rr^{d_{1}\times\dots\times d_{k}}$ has no solution in general for
any $r=2,\dots,\min\{d_{1},\dots,d_{k}\}$. In particular, there
exists $A\in\Rr^{d_{1}\times\dots\times d_{k}}$ with
\[
\trank(A)=r+1
\]
that has no best rank-$r$ approximation. The result is independent of the choice of norms.
\end{theorem}

\begin{proof}
We begin by assuming $k=3$.

\textsc{Higher rank}.
Let $2 \leq r \leq \min\{d_1, d_2, d_3\}$. By Lemma~\ref{buildrank}, we can construct a tensor $B \in\Rr^{(d_1-2)\times (d_2-2) \times (d_3-2)}$ with rank~$r-2$. By Proposition~\ref{prop:dSLtensor}, we can construct a convergent sequence of tensors $C_n \to C$ in~$\Rr^{2\times 2\times 2}$ with $\trank(C_n) \leq 2$, and $\trank(C) = 3$. Let $A_n = B \oplus C_n \in \Rr^{d_1 \times d_2 \times d_3}$. Then $A_n \to A := B \oplus C$ and $\trank(A_n) \leq \trank(B)  + \trank(C_n) \leq r$. The result of J\'aJ\'a--Takche (Theorem~\ref{jaja}) implies that $\trank(A) = \trank(B) + \trank(C) = r+1$.

\textsc{Arbitrary order}. Let $\uu_{4}\in\Rr^{d_{4}},\dots,\uu_{k}\in\Rr^{d_{k}}$ be unit vectors and set
\[
\tilde{A}_{n}:=A_{n}\otimes\uu_{4}\otimes\dots\otimes\uu%
_{k},\qquad\tilde{A}:=A\otimes\uu_{4}\otimes\dots\otimes\uu%
_{k}.
\]
By \eqref{eqn:Fnorm3},
\[
\lVert \tilde{A}_{n}-\tilde{A} \rVert_{F}=
\lVert A_n - A \rVert =
\lVert B \oplus C_n - B \oplus C \rVert =
\lVert C_n - C \rVert
\to 0,
\mbox{ as $n\to \infty$.}
\]
Moreover, Corollary~\ref{cor:rank} ensures that $\trank(\tilde{A})=
r+1$ and $\trank(\tilde{A_n})\le r$.

\textsc{Norm independence}. Whether the sequence $\tilde{A}_{n}$
converges to~$\tilde{A}$ is entirely dependent on the norm-induced
topology on $\Rr^{d_{1}\times\dots\times d_{k}}$. Since it has a
unique topology induced by any of its equivalent norms as a
finite-dimensional vector space, the convergence is independent of
the choice of norms.
\end{proof}

We note that the proof above exhibits an order-$k$ tensor, namely
$\tilde{A}$, that has rank strictly larger than $\min \{d_1, \dots, d_k
\}$.

\subsection{Tensor rank can leap an arbitrarily large gap}
\label{subsec:leap}

How can we construct a sequence of tensors of rank~$r$ that converge
to a tensor of rank~$r+2$? An easy trick is to take the direct sum
of two sequences of rank-$2$ tensors of the form shown in~\eqref{dSLseq}.
%
The resulting sequence converges to a limiting tensor that
is the direct sum of two rank-$3$ tensors, each of form
shown in~\eqref{eqdSL}. To show that the limiting tensor has rank~6 (and
does not have some miraculous lower-rank decomposition), we once
again turn to the theorem of J\'aJ\'a--Takche, which contains just
enough of the direct sum conjecture~\eqref{DSC} for our purposes.

\begin{proposition}
\label{prop:leap}Given any $s\in\mathbb{N}$ and $r \geq 2s$, there
exists a sequence of order-$3$ tensors~$B_{n}$ such that
$\trank(B_{n})\leq r$ and $\lim_{n\rightarrow\infty}B_{n}=B$ with
{$\trank(B) = r+s$}.
\end{proposition}

\begin{proof}
Let $d = r-2s$. By Lemma~\ref{buildrank}, there exists a rank-$d$
tensor $C \in \Rr^{d \times d \times d}$. Let $A_n \to A$ be a
convergent sequence in $\Rr^{2\times2\times2}$ with $\trank(A) \leq
2$ and $\trank(A) = 3$. Define
\[
B_n = C \oplus {A_n \oplus \dots \oplus A_n},
\quad
B = C \oplus {A\oplus \dots \oplus A}
\]
where there are $s$ terms~$A_n$ and~$A$. Then $B_n \to B$, and $\trank(B_n) \leq r-2s + 2s = r$. By applying the J\'{a}J\'{a}--Takche sequentially $s$~times, once for each summand~$A$, we deduce that $\trank(B) = r-2s + 3s = r+s$.
\end{proof}

As usual the construction can be extended to order-$k$ tensors, by taking an outer product with a suitable number of non-zero vectors in the new factors.
\begin{corollary}
\label{cor:leap}Given any $s \ge 1$, $r \ge 2$, and $k \ge 3$, with $r \geq 2s$, there exists $A \in \mathbb{R}^{d_1 \times \dots \times d_k}$ such that $\trank(A) = r + s$ and $A$ has no best rank-$r$ approximation.
\end{corollary}

\begin{proof}
This follows from Proposition~\ref{prop:leap} and the previous remark.
\end{proof}

\subsection{Br\`{e}gman divergences and other continuous measures of proximity}
\label{subsec:bregman}

In data analytic applications, one frequently encounters low-rank
approximations with respect to `distances' that are more general than norms.
Such a `distance' may not even be a metric, an example being the Br\`{e}gman
divergence \cite{Bre, DT} (sometimes also known as Br\`{e}gman distance). The
definition here is based on the definition given in \cite{DT}. Recall first that if $S \subset \mathbb{R}^n$, the \textit{relative interior} of $S$ is simply the interior of $S$ considered as a subset of its affine hull, and is denoted by $\operatorname{ri}(S)$.

\begin{definition}
Let $S\subseteq\Rr^{d_{1}\times\dots\times d_{k}}$ be a {convex set}.
Let $\varphi:S\rightarrow\Rr$ be a lower semicontinuous, {convex} function that is {continuously differentiable and strictly convex} in $\operatorname{ri}(S)$.
Let $\varphi$ have the property that for any sequence $\{C_{n}\}\subset
\operatorname{ri}(S)$ that converges to $C\in S\setminus\operatorname{ri}%
(S)$, we have:
\[
\lim_{n\rightarrow\infty}\lVert\nabla\varphi(C_{n})\rVert=+\infty.
\]
The \textbf{Br\`{e}gman divergence} $D_{\varphi}:S\times\operatorname{ri}%
(S)\rightarrow\Rr$ is defined by%
\[
D_{\varphi}(A,B)=\varphi(A)-\varphi(B)-\langle\nabla\varphi(B),A-B\rangle.
\]

\end{definition}

It is natural to ask if the analogous problem {\brap} for Br\`{e}gman divergence will always have a solution. Note that a Br\`{e}gman divergence, unlike a metric, is not necessarily symmetric in its two arguments and thus there are \textit{two} possible problems:
\[
\operatorname*{argmin}\nolimits_{\operatorname*{rank}_{\otimes}(B)\leq
r}D_{\varphi}(A,B) \qquad \text{and} \qquad \operatorname*{argmin}\nolimits_{\operatorname*{rank}_{\otimes}(B)\leq
r}D_{\varphi}(B,A).
\]
As the following proposition shows, the answer is no in both cases.

\begin{proposition}
\label{prop:bregman}Let $D_{\varphi}$ be a Br\`{e}gman divergence. Let $A$ and
$A_{n}$ be defined as in \eqref{eqdSL} and \eqref{dSLseq} respectively. Then
\[
\lim_{n\rightarrow\infty}D_{\varphi}(A,A_{n})=0 = \lim_{n\rightarrow\infty}D_{\varphi}(A_n,A).
\]
\end{proposition}

\begin{proof}
The Br\`{e}gman divergence is jointly continuous in both arguments with respect to the norm topology, and $A_{n} \to A$ in norm, so $D_{\varphi}(A, A_n) \to D_\varphi(A,A) = 0 $ and $D_{\varphi}(A_n, A) \to D_\varphi(A,A) = 0 $.
\end{proof}
%

Proposition~\ref{prop:bregman} extends trivially to any other
measure of nearness that is continuous with respect to the norm topology in at least one argument.

\subsection{Difference quotients}
\label{subsec:derive}
We thank Landsberg~\cite{Land1} for the insight that the expression in~\eqref{eqdSL} is best regarded as a derivative. Indeed, if
\[
f(t) = (\xx + t\yy)^{\otimes 3}
= (\xx + t\yy) \otimes (\xx + t\yy) \otimes (\xx + t\yy)
\]
%
then
\[
\left. \frac{df}{dt} \right|_{t=0}
 =
 \yy \otimes \xx \otimes \xx
+ \xx \otimes \yy \otimes \xx
+ \xx \otimes \xx \otimes \yy
\]
by the Leibniz rule. On the other hand
\[
\left. \frac{df}{dt} \right|_{t=0}
=
\lim_{t \to 0}
\left[
\frac{(\xx+ t \yy)\otimes(\xx+ t \yy)\otimes(\xx+ t \yy)
- \xx \otimes \xx \otimes \xx}{t}
\right]
\]
and the difference quotient on the right-hand side has rank~2. The expression in~\eqref{dSLseq} can be obtained from this by taking $t = 1/N$.

%

We can extend Landsberg's idea to more general partial differential operators. It will be helpful to use the degree-$k$ Veronese map~\cite{H1}, which is
$
V_k(\xx) = \xx^{\otimes k} =
{\xx \otimes \dots \otimes \xx}
$ ($k$-fold product).
Then, for example, the 6-term symmetric tensor
\[
\xx \otimes \yy \otimes \zz
+\xx \otimes \zz \otimes \yy
+\yy \otimes \zz \otimes \xx
+\yy \otimes \xx \otimes \zz
+\zz \otimes \xx \otimes \yy
+\xx \otimes \yy \otimes \xx
\]
can be written as a partial derivative
\[
\left. \frac{\partial^2}{\partial s \, \partial t} \right|_{s=t=0}
(\xx + s\yy + t\zz)^{\otimes 3}
\]
which is a limit of a 4-term difference quotient:
\[
\lim_{s, t\to 0}
\left[
\frac{V_3(\xx+s\yy+t\zz) - V_3(\xx+s\yy) - V_3(\xx+t\zz)+V_3(\xx)}
{st}
\right]
\]
This example lies naturally in $\Rr^{3\times3\times3}$, taking $\xx, \yy, \zz$ to be linearly independent.
Another example, in $\Rr^{2\times2\times2\times2}$, is the 6-term symmetric order-4 tensor
\begin{eqnarray*}
&&\xx \otimes \xx \otimes \yy \otimes \yy
+\xx \otimes \yy \otimes \xx \otimes \yy
+\xx \otimes \yy \otimes \yy \otimes \xx
\\
&&\qquad
{} + \yy \otimes \xx \otimes \xx \otimes \yy
+\yy \otimes \xx \otimes \yy \otimes \xx
+\yy \otimes \yy \otimes \xx \otimes \xx.
\end{eqnarray*}
This can be written as the second-order derivative
\[
\left. \frac{\partial^2}{\partial t^2} \right|_{t=0}
\frac{(\xx + t\yy)^{\otimes 4}}{2!}
\]
which is a limit of a 3-term difference quotient:
\[
\lim_{t\to 0}
\left[
\frac{V_4(\xx+2t\yy) - 2V_4(\xx+t\yy) +V_4(\xx)}
{2!\,t^2}
\right]
\]

%
%

%
We call these examples {\bf symmetric Leibniz tensors} for the differential operators ${\partial^2}/{\partial s\,\partial t}$ and ${\partial^2}/{\partial t^2}$, of orders 3 and~4, respectively.
More generally, given positive integers $k$ and $a_1, \dots, a_j$ with $a_1 + \dots + a_j = a \leq k$, the symmetric tensor
\begin{eqnarray*}
L_k(a_1, \dots, a_j)
&:=&
\sum_{\mathrm{Sym}} \xx^{\otimes(k-a)} \otimes \yy_1^{\otimes a_1} \otimes
\dots \otimes \yy_j^{\otimes a_j}
\end{eqnarray*}
can be written as a partial derivative
\[
\left.
\frac{\partial^a}{\partial {t_1}^{a_1} \dots \partial {t_j}^{a_j}}
\right|_{t_1 = \,\dots\, = t_j = 0}
\frac{V_k(\xx + t_1 \yy_1 + \dots + t_j \yy_j)}
{(a_1!) \cdots (a_j!)}
\]
which is a limit of a difference quotient with $(a_1+1)\cdots(a_j+1)$ terms. On the other hand, the number of terms in the limit $L_k(a_1, \dots, a_j)$ is given by a multinomial coefficient, and that is usually much bigger.

This construction gives us a ready supply of candidates for rank-jumping.
However, we do not know --- even for the two explicit 6-term examples above
--- whether the limiting tensors actually have the ranks suggested by their
formulas. We can show that $\trank(L_k(1)) = k$, for all~$k$ and over any field,
generalizing Lemma~\ref{lem:dSLtensor-rank}. Beyond that it is not clear to us what is
likely to be true. The optimistic conjecture is:
\begin{equation}
\trank(L_k(a_1, \dots, a_j)) \stackrel{?}{=}
\binom{k}{k-a, a_1, \dots, a_j}
= \frac{k!}{(k-a)!\,a_i!\,\cdots\,a_j!}
\label{eq:dqconj}
\end{equation}
Comon et al.~\cite{CGLM2} show that the \emph{symmetric} rank of $L_k(1)$
over the complex numbers is~$k$, so that is another possible context which \eqref{eq:dqconj} may be true.

\section{Characterizing the limit points of order-$3$ rank-$2$ tensors\label{sect:char32}}

If an order-$3$ tensor can be expressed as a limit of a sequence of
rank-$2$ tensors, but itself has rank greater than~$2$, then we show
in this section that it takes a particular form. This kind of result
may make it possible to overcome the ill-posedness of {\brap}, by
defining weak solutions.

\begin{theorem}
\label{thm:b}
Let $d_{1},d_{2},d_{3}\geq2$. Let $A_{n}\in\Rr^{d_{1}\times d_{2}\times
d_{3}}$ be a sequence of tensors with $\trank(A_{n})\leq2$ and
\[
\lim_{n\rightarrow\infty}A_{n}=A,
\]
where the limit is taken in any norm topology. If the limiting tensor $A$ has
rank higher than $2$, then $\trank(A)$ must be exactly
$3$ and there exist pairs of linearly independent vectors $\xx%
_{1},\yy_{1}\in\Rr^{d_{1}}$, $\xx_{2},\yy_{2}%
\in\Rr^{d_{2}}$, $\xx_{3},\yy_{3}\in\Rr^{d_{3}}$
such that
\begin{equation}
A=\xx_{1}\otimes\xx_{2}\otimes\yy_{3}+\xx%
_{1}\otimes\yy_{2}\otimes\xx_{3}+\yy_{1}\otimes
\xx_{2}\otimes\xx_{3}.
\end{equation}
\end{theorem}

The proof of this theorem will occupy the next few subsections.


\subsection{Reduction\label{sect:red}}

Our first step is to show that we can limit our attention to the
particular tensor space~$\Rr^{2 \times 2 \times 2}$. Here the
orthogonal group action is
important. Recall that the actions of $\Oo_{d_{1},\dots,d_{k}%
}(\Rr)$ and $\GL_{d_{1},\dots,d_{k}}(\Rr)$ on
$\Rr^{d_{1} \times\dots\times d_{k}}$ are continuous and carry
decomposable tensors to decomposable tensors. It follows that the subspaces
$\Ss_{r}$ and $\Ssbar_{r}$ are preserved. The next
theorem provides a general mechanism for passing to a tensor subspace.

\begin{theorem}
\label{thm:reduction:k} Let $r_{i}=\min(r,d_{i})$ for all~$i$. The restricted
maps
\begin{gather*}
\Oo_{d_{1},\dots,d_{k}}(\Rr) \times \Ss_{r}(r_{1},\dots,r_{k}) \rightarrow \Ss_{r}(d_{1},\dots,d_{k})\\
\Oo_{d_{1},\dots,d_{k}}(\Rr) \times
\overline{\Ss}_{r}(r_{1},\dots,r_{k}) \rightarrow
\overline{\Ss}_{r}(d_{1},\dots,d_{k})
\end{gather*}
given by $((L_{1},\dots,L_{k}),A)\mapsto(L_{1},\dots,L_{k})\cdot A$ are both surjective.
\end{theorem}

In other words, every rank-$r$ tensor in $\Rr^{d_{1}\times\dots\times
d_{k}}$ is equivalent by an orthogonal transformation to a rank-$r$ tensor in
the smaller space $\Rr^{r_{1}\times\dots\times r_{k}}$. Similarly every
rank-$r$-approximable tensor in $\Rr^{d_{1}\times\dots\times d_{k}}$ is
equivalent to a rank-$r$-approximable tensor in $\Rr^{r_{1}\times
\dots\times r_{k}}$.

\begin{proof}
If $A\in\Ss_{r}(d_{1},\dots,d_{k})$ is any rank-$r$ tensor then we can
write $A=\sum_{j=1}^{r}\xx_{1}^{j}\otimes\dots\otimes\xx_{k}%
^{j}$ for vectors $\xx_{i}^{j}\in\Rr^{d_{i}}$. For each~$i$, the
vectors $\xx_{i}^{1},\dots,\xx_{i}^{r}$ span a subspace
$V_{i}\subset\Rr^{d_{i}}$ of rank at most~$r_{i}$. Choose $L_{i}%
\in\Oo_{d_{i}}(\Rr)$ so that $L_{i}(\Rr^{d_{i}%
})\supseteq V_{i}$. Let $B=(L_{1}^{-1},\dots,L_{k}^{-1})\cdot A$. Then
$A=(L_{1},\dots,L_{k})\cdot B$ and $B\in\Ss_{r}(d_{1},\dots,d_{k})$.
This argument shows that the first of the maps is surjective.

Now let $A\in\Ssbar_{r}(d_{1},\dots,d_{k})$ be any
rank-$r$-approximable tensor. Let $(A^{(n)})_{n=1}^{\infty}$ be any sequence
of rank-$r$ tensors converging to~$A$. For each~$n$, by the preceding result,
we can find $B^{(n)}\in\Ssbar_{r}(d_{1},\dots,d_{k})$ and
$(L_{1}^{(n)},\dots,L_{k}^{(n)})\in\Oo_{d_{1},\dots,d_{k}%
}(\Rr)$ with $(L_{1}^{(n)},\dots,L_{k}^{(n)})\cdot B^{(n)}=A^{(n)}$.
Since $\Oo_{d_{1},\dots,d_{k}}(\Rr)$ is compact, there is
a convergent subsequence $(L_{1}^{(n_{j})},\dots,L_{k}^{(n_{j})}%
)\rightarrow(L_{1},\dots,L_{k})$. Let $B=(L_{1},\dots,L_{k})^{-1}\cdot A$.
Then $A=(L_{1},\dots,L_{k})\cdot B$; and $B^{(n_{j})}=(L_{1}^{(n_{j})}%
,\dots,L_{k}^{(n_{j})})^{-1}\cdot A^{(n_{j})}\rightarrow(L_{1},\dots
,L_{k})^{-1}\cdot A=B$, so $B\in\Ssbar_{r}(d_{1},\dots,d_{k}%
)$. Thus the second map is also surjective.
\end{proof}

\begin{corollary}
\label{cor:reduceto222} If Theorem~\ref{thm:b} is true for the tensor
space~$\Rr^{2\times2\times2}$ then it is true in general.
\end{corollary}

\begin{proof}
The general case is $V_{1}\otimes V_{2}\otimes V_{3}\cong\Rr^{d_1\times
d_2\times d_3}$. Suppose $A\in\Ssbar_{2}(d_1,d_2,d_3)$ and
$\trank(A)\geq3$. By Theorem~\ref{thm:reduction:k},
there exists $(L_{1},L_{2},L_{3})\in\Oo_{d_1,d_2,d_3}(\Rr)$ and
$B\in\Ssbar_{2}(2,2,2)$ with $(L_{1},L_{2},L_{3})\cdot B=A$.
Moreover $\trank(B)=\trank(A)\geq3$ in $\Rr^{l\times m\times n}$ and hence $\trank(B)\geq3$ in $\Rr^{2\times2\times2}$ by Proposition~\ref{prop:rank}. Since the theorem is assumed true for~$\Rr^{2\times2\times2}$ and $B$~satisfies the hypotheses, it can be written in the specified form in terms of vectors $\xx_{1},\xx_{2}%
,\xx_{3}$ and $\yy_{1},\yy_{2},\yy_{3}$. It
follows that $A$ takes the same form with respect to the vectors
$L_{1}\xx_{1},L_{2}\xx_{2},L_{3}\xx_{3}$ and
$L_{1}\yy_{1},L_{2}\yy_{2},L_{3}\yy_{3}$.
\end{proof}


\subsection{Tensors of rank~$1$ and~$2$}

We establish two simple facts, for later use.

\begin{proposition}
\label{prop:reduction:1} If $A \in\Rr^{d_{1} \times\cdots\times
d_{k}}$ has rank~$1$, then we can write $A =
(L_{1},\dots,L_{k})\cdot B$, where $(L_{1}, \dots, L_{k})
\in\GL_{d_{1}, \dots, d_{k}}(\Rr)$ and $B = \ee_{1}
\otimes\dots\otimes\ee_{k}$.
\end{proposition}

\begin{proof}
Write $A = \xx_{1} \otimes\dots\otimes\xx_{k}$ and choose
the~$L_{i}$ so that $L_{i}(\ee_{i}) = \xx_{i}$.
\end{proof}

\begin{proposition}
\label{thm:reduction:2} Assume $d_{i}\geq2$ for all~$i$. If $A\in
\Rr^{d_{1}\times\dots\times d_{k}}$ has rank~$2$, then we can write
$A =
(L_{1},\dots,L_{k})\cdot B$, where $(L_{1},\dots,L_{k})\in\GL_{d_{1}%
,\dots,d_{k}}(\Rr)$ and $B\in\Rr^{2\times\dots\times2}$ is of
the form $B=\ee_{1}\otimes\dots\otimes\ee_{1}+\ff%
_{1}\otimes\dots\otimes\ff_{k}$. Here $\ee_{1}$ denotes the standard
basis vector $(1,0)^{\transp}$; each $\ff_{i}$ is equal either
to~$\ee_{1}$ or to $\ee_{2}=(0,1)^{\transp}$; and at least two of
the~$\ff_{i}$ are equal to~$\ee_{2}$.
\end{proposition}

\begin{proof}
We can write $A=\xx_{1}\otimes\dots\otimes\xx_{k}+\yy%
_{1}\otimes\dots\otimes\yy_{k}$. Since $\trank(A)=2$ all of the $\xx_{i}$ and $\yy_{i}$ must be nonzero. We
claim that $\yy_{i},\xx_{i}$ must be linearly independent for at
least two different indices~$i$. Otherwise, suppose $\yy_{i}%
=\lambda_{i}\xx_{i}$ for $k-1$ different indices, say $i=1,\dots,k-1$.
It would follow that
\[
A=\xx_{1}\otimes\dots\otimes\xx_{k-1}\otimes\ (\xx%
_{k}+(\lambda_{1}\cdots\lambda_{k-1})\yy_{k})
\]
contradicting $\trank(A)=2$.

For each~$i$ choose $L_{i}:\Rr^{2}\rightarrow\Rr^{d_{i}}$ such
that $L_{i}\ee_{1}=\xx_{i}$, and such that $L_{i}\ee%
_{2}=\yy_{i}$ if $\yy_{i}$ is linearly independent
of~$\xx_{i}$; otherwise $L_{i}\ee_{2}$ may be arbitrary. It is
easy to check that $(L_{1},\dots,L_{k})^{-1}\cdot A=\ee_{1}\otimes
\dots\otimes\ee_{1}+\lambda\ff_{1}\otimes\dots\otimes
\ff_{k}$ where the $\ff_{i}$ are as specified in the theorem,
and $\lambda$ is the product of the~$\lambda_{i}$ over those indices where
$\yy_{i}=\lambda_{i}\xx_{i}$. This is almost in the correct form. To get rid
of the~$\lambda$, replace $L_{i}\ee_{2}=\yy_{i}$ with
$L_{i}\ee_{2}=\lambda\yy_{i}$ at one of the indices~$i$ for
which $\xx_{i},\yy_{i}$ are linearly independent. This completes
the construction.
\end{proof}


\subsection{The discriminant polynomial~$\Delta$\label{subsec:delta}}

The structure of tensors in~$\Rr^{2\times2\times2}$ is largely
governed by a quartic polynomial~$\Delta$ which we define and
discuss here. This same polynomial was discovered by Cayley in 1845
\cite{Cay}.
More generally, $\Delta$ is the $2\times2\times2$ special case of an
object called the \textit{hyperdeterminant} revived in its modern
form by Gelfand, Kapranov, and Zelevinsky \cite{GKZ1, GKZ2}. We give
an elementary treatment of the properties we need.

As in our discussion in Section~\ref{sect:MMM}, we identify a tensor
$\mathbf{A}\in\Rr^{2}\otimes\Rr^{2}\otimes\Rr^{2}$ with
the array $A\in\Rr^{2\times2\times2}$ of its eight coefficients with
respect to the standard basis $\{\ee_{i}\otimes\ee_{j}%
\otimes\ee_{k}:i,j,k=1,2\}$. Pictorially, we can represent it as a pair
of side-by-side $2\times2$~slabs:
\[
\mathbf{A}=\sum_{i=1}^{2}\sum_{j=1}^{2}\sum_{k=1}^{2}a_{ijk}\ee%
_{i}\otimes\ee_{j}\otimes\ee_{k}=\left[
\begin{array}
[c]{rr}%
a_{111} & a_{112}\\
a_{121} & a_{122}%
\end{array}
\right\vert \!\left.
\begin{array}
[c]{rr}%
a_{211} & a_{212}\\
a_{221} & a_{222}%
\end{array}
\right]  =A
\]
The general strategy is to find ways of simplifying the representation of~$A$
by applying transformations in~$\GL_{2,2,2}(\Rr%
)=\GL_{2}(\Rr)\times\GL_{2}(\Rr%
)\times\GL_{2}(\Rr)$. This group is generated by the
following operations: decomposable row operations applied to both slabs
simultaneously; decomposable column operations applied to both slabs
simultaneously; decomposable slab operations (for example, adding a multiple
of one slab to the other).

Slab operations on a tensor $A=[A_{1}\,|\,A_{2}]$ generate new $2\times
2$~slabs of the form $S=\lambda_{1}A_{1}+\lambda_{2}A_{2}$. One can check
that:
\begin{equation}
\det(S)=\lambda_{1}^{2}\det(A_{1})+\lambda_{1}\lambda_{2}\frac{\det
(A_{1}+A_{2})-\det(A_{1}-A_{2})}{2}+\lambda_{2}^{2}\det(A_{2})
\label{eqn:quadratic}%
\end{equation}
We define $\Delta$ to be the discriminant of this quadratic polynomial:
\begin{equation}
\Delta([A_{1}\,|\,A_{2}])=\left[  \frac{\det(A_{1}+A_{2})-\det(A_{1}-A_{2}%
)}{2}\right]  ^{2}-4\det(A_{1})\det(A_{2}) \label{eqn:DeltaT}%
\end{equation}
Explicitly, if $A = \llbracket a_{ijk} \rrbracket_{i,j,k=1,2} \in
\mathbb{R}^{2\times 2\times 2}$, then
%
\begin{multline*}
\Delta(A)=(a_{111}^{2}a_{222}^{2}+a_{112}^{2}a_{221}^{2}+a_{121}^{2}%
a_{212}^{2}+a_{122}^{2}a_{211}^{2})\\
-2(a_{111}a_{112}a_{221}a_{222}+a_{111}a_{121}a_{212}a_{222}+a_{111}%
a_{122}a_{211}a_{222}\\
+a_{112}a_{121}a_{212}a_{221}+a_{112}a_{122}a_{221}a_{211}+a_{121}%
a_{122}a_{212}a_{211})\\
+4(a_{111}a_{122}a_{212}a_{221}+a_{112}a_{121}a_{211}a_{222}).
\end{multline*}

\begin{proposition}
Let $A\in\Rr^{2\times2\times2}$, let $A^{\prime}$ be obtained from $A$
by permuting the three factors in the tensor product, and let $(L_{1}%
,L_{2},L_{3})\in\GL_{2,2,2}(\Rr)$. Then $\Delta
(A^{\prime})=\Delta(A)$ and $\Delta((L_{1},L_{2},L_{3})\cdot A)=\det
(L_{1})^{2}\det(L_{2})^{2}\det(L_{3})^{2}\Delta(A)$.
\end{proposition}

\begin{proof}
To show that $\Delta$ is invariant under all permutations of the
factors of~$\Rr^{2\times2\times2}$, it is enough to check invariance
in the cases of two distinct transpositions. It is clear from
equation~(\ref{eqn:DeltaT}) that $\Delta$~is invariant under the
transposition of the second and third factors, since this amounts to
replacing $A_{1},A_{2}$ with their transposes
$A_{1}^{\transp},A_{2}^{\transp}$. To show that $\Delta$~is
invariant under transposition of the first and third factors, write
$A=[\uu_{11},\,\uu_{12}\mid\uu_{21},\,\uu_{22}]$, where the
$\uu_{ij}$ are column vectors. One can verify that
\begin{align*}
\Delta(A)  &  =\det[\uu_{11},\uu_{22}]^{2}+\det[\uu%
_{21},\uu_{12}]^{2}\\
&  \quad-2\det[\uu_{11},\uu_{12}]\det[\uu_{21}%
,\uu_{22}]-2\det[\uu_{11},\uu_{21}]\det[\uu%
_{12},\uu_{22}]
\end{align*}
which has the desired symmetry.

In view of the permutation invariance of~$\Delta$, it is enough to verify the
second claim in the case $(L_{1},L_{2},L_{3})=(I,L_{2},I)$. Then $(L_{1}%
,L_{2},L_{3})\cdot A=[L_{2}A_{1}\,|\,L_{2}A_{2}]$ and an extra factor~$\det
(L_{2})^{2}$ appears in all terms of equation~(\ref{eqn:DeltaT}), exactly as required.
\end{proof}

\begin{corollary}
\label{cor:sign} The sign of~$\Delta$ is invariant under the action of
$\GL_{2,2,2}(\Rr)$.
\end{corollary}

\begin{corollary}
\label{cor:signO} The value of~$\Delta$ is invariant under the action of
$\Oo_{2,2,2}(\Rr)$.
\end{corollary}

Using the properties of $\Delta$, we can easily prove, in a slightly different
way, a result due originally to Kruskal (unpublished work) and ten Berge
\cite{tB}.

\begin{proposition}
\label{prop:rank2:b} If $\Delta(A) > 0$ then $\trank(A)\leq 2$.
\end{proposition}

\begin{proposition}
\label{prop:rank2:a} If $\trank(A)\leq2$ then
$\Delta(A)\geq0$.
\end{proposition}

{\it Proof of Proposition~\ref{prop:rank2:b}.}
If the discriminant $\Delta(A)$ is positive then the homogeneous quadratic
equation~(\ref{eqn:quadratic}) has two linearly independent root pairs
$(\lambda_{11},\lambda_{12})$ and $(\lambda_{21},\lambda_{22})$. It follows
that we can use slab operations to transform $[A_{1}\,|\,A_{2}]\rightarrow
\lbrack B_{1}\,|\,B_{2}]$, where $B_{i}=\lambda_{i1}A_{1}+\lambda_{i2}A_{2}$.
By construction $\det(B_{i})=0$ so we can write $B_{i}=\ff_{i}%
\otimes\mathbf{g}_{i}$ for some $\ff_{i},\mathbf{g}_{i}\in
\Rr^{2}$ (possibly zero). It follows that $[B_{1}\,|\,B_{2}%
]=\ee_{1}\otimes\ff_{1}\otimes\mathbf{g}_{1}+\ee%
_{2}\otimes\ff_{2}\otimes\mathbf{g}_{2}$; so $\trank(A)=\trank([B_{1}\,|\,B_{2}])\leq2$.
\endproof

{\it Proof of Proposition~\ref{prop:rank2:a}.}
It is easy to check that
$\Delta(A)=0$ if $\trank(A)\leq1$, since we can write
$A=(L_{1},L_{2},L_{3})\cdot(\ee_{1}\otimes\ee_{1}\otimes
\ee_{1})$ or else $A=0$.

It remains to be shown that $\Delta(A)$ is not negative when
$\trank(A)=2$. Proposition~\ref{thm:reduction:2}
implies that $A$ can be transformed by an element of~$\GL%
_{2,2,2}(\Rr)$ (and a permutation of factors, if necessary) into one of
the following tensors:
\[
I_{1}=\left[
\begin{array}
[c]{rr}%
1 & 0\\
0 & 0
\end{array}
\right\vert \!\left.
\begin{array}
[c]{rr}%
0 & 0\\
0 & 1
\end{array}
\right]  \quad\text{or}\quad I_{2}=\left[
\begin{array}
[c]{rr}%
1 & 0\\
0 & 1
\end{array}
\right\vert \!\left.
\begin{array}
[c]{rr}%
0 & 0\\
0 & 0
\end{array}
\right]
\]
Since $\Delta(I_{1})=1$ and $\Delta(I_{2})=0$ it follows that $\Delta(A)\geq0$.
\endproof

Kruskal and also ten Berge deserve complete credit for {discovering} the above
result. In fact, the hyperdeterminant for $2\times2\times2$ tensor
$\Delta$ is known by the name \textit{Kruskal polynomial} in the
psychometrics community {\cite{tB}}. Our goal is not so much to provide
alternative proofs for Propositions \ref{prop:rank2:b}
and~\ref{prop:rank2:a} but to include them so that our proof of
Theorem~\ref{thm:b} can be self-contained. We are now ready to
give that proof, thereby characterizing all limit points of order-$3$ rank-$2$ tensors.


{\it Proof of Theorem~\ref{thm:b}.}
Note that the theorem is stated for order-$3$ tensors of any size
$d_{1}\times d_{2}\times d_{3}$. We begin with the case $A \in \Rr^{2\times2\times2}$. Suppose $A\in\Ssbar_{2}(2,2,2)\setminus\Ss_{2}(2,2,2)$. Then we claim that $\Delta(A)=0$. Indeed, since $A \not\in\Ss_{2}$,  Proposition~\ref{prop:rank2:b} implies that $\Delta(A)\leq0$. On the other hand, since $A \in \Ssbar_{2}$, it follows from Proposition~\ref{prop:rank2:a} and the continuity of~$\Delta$ that $\Delta(A)\geq0$.

Since $\Delta(A)=0$, the homogeneous quadratic
equation~(\ref{eqn:quadratic}) has a nontrivial root pair $(\lambda_{1},\lambda_{2})$.
It follows that $A$ can be transformed by slab operations into the form
$[A_{i}\,|\,S]$ where $S=\lambda_{1}A_{1}+\lambda_{2}A_{2}$ and $i=1$ or~2. By
construction $\det(S)=0$, but $S\neq0$ for otherwise $\trank(A)=\operatorname*{rank}(A_{i})\leq2$. Hence $\rank(S)=1$ and by a
further transformation we can reduce~$A$ to the form:
\[
B=\left[
\begin{array}
[c]{rr}%
p & q\\
r & s
\end{array}
\right\vert \!\left.
\begin{array}
[c]{rr}%
1 & 0\\
0 & 0
\end{array}
\right]
\]
In fact we may assume $p=0$ (the operation `subtract $p$
times the second slab from the first slab' will achieve
this), and moreover $s^2=\Delta(B)=0$. Both $q$ and~$r$ must be non-zero, otherwise $\trank(A)=\trank(B)\leq2$. If we rescale
the bottom rows by~$1/r$ and the right-hand columns by~$1/q$ we are finally
reduced to:
\[
B^{\prime}=\left[
\begin{array}
[c]{rr}%
0 & 1\\
1 & 0
\end{array}
\right\vert \!\left.
\begin{array}
[c]{rr}%
1 & 0\\
0 & 0
\end{array}
\right]  =\ee_{2}\otimes\ee_{1}\otimes\ee_{1}%
+\ee_{1}\otimes\ee_{2}\otimes\ee_{1}+\ee%
_{1}\otimes\ee_{1}\otimes\ee_{2}%
\]
By reversing all the row, column and slab operations we can obtain a
transformation $(L_{1},L_{2},L_{3})\in\GL_{2,2,2}(\Rr)$
such that $A=(L_{1},L_{2},L_{3})\cdot B^{\prime}$. Then $A$~can be written in
the required form, with $\xx_{i}=L_{i}\ee_{1}$, $\yy%
_{i}=L_{i}\ee_{2}$ for $i=1,2,3$.

This completes the proof of Theorem~\ref{thm:b} in the case of the tensor
space~$\Rr^{2\times2\times2}$. By Corollary~\ref{cor:reduceto222} this
implies the theorem in general.
\endproof

\subsection{Ill-posedness and ill-conditioning of the best rank-$r$ approximation
problem}
\label{sect:condition}

Recall that a problem is called \textit{well-posed} if a solution
exists, is unique, and is stable (i.e.\ depends continuously on the
input data). If one or more of these three criteria are not
satisfied, the problem\footnote{Normally, existence is taken for
granted and an ill-posed problem often means one whose solution
lacks either uniqueness or stability. In this paper, the
ill-posedness is of a more serious kind --- the existence of a
solution is itself in question.} is called \textit{ill-posed}.

From Sections \ref{sec:top} and \ref{sec:volume}, we see that tensors will often fail to have a best rank-$r$ approximation. In all applications that rely on {\brap} or a variant of it as the underlying mathematical model, we should fully expect the ill-posedness of {\brap} to pose a serious difficulty. Even if it is known \textit{a priori} that a tensor $A$ has a best rank-$r$ approximation, we should remember that in applications, the data array $\hat{A}$ available at our disposal is almost always one that is corrupted by noise, i.e.\ $\hat{A} = A + E$ where $E$ denotes the collective contributions of various errors --- limitations in measurements, background noise, rounding off, etc. Clearly there is no guarantee that $\hat{A}$ will also have a best rank-$r$ approximation.

In many situations, one only needs a `good' rank-$r$ approximation rather than the best rank-$r$ approximation. It is tempting to argue, then, that the non-existence of the best solution does not matter --- it is enough to seek an `approximate solution'. We discourage this point of view, for two main reasons. First, there is a serious conceptual difficulty: if there is no solution, then what is the `approximate solution' an approximation of? Second, even if one disregards this, and ploughs ahead regardless to compute an `approximate solution', we argue below that this task is ill-conditioned and the computation is unstable.

For notational simplicity and since there is no loss of generality (cf.\ Theorem~\ref{thm:rank+1} and Corollary~\ref{cor:leap}), we will use the problem of finding a best rank-$2$ approximation to a rank-$3$ tensor to make our point. Let $A\in\Rr^{d_{1}\times d_{2}\times d_{3}}$ be an instance where
\begin{equation}
\operatorname*{argmin}\nolimits_{\xx_{i},\yy_{i}\in\Rr^{d_{i}}%
}\lVert A-\xx_{1}\otimes\xx_{2}\otimes\xx_{3}%
-\yy_{1}\otimes\yy_{2}\otimes\yy_{3}\rVert\label{nosoln}%
\end{equation}
does not have a solution (such examples abound, cf.\ Section \ref{sec:volume}). If we disregard the fact that a solution does not exist and plug the problem into a computer program%
\footnote{While there is no known globally convergent algorithm for {\brap}, we will ignore this difficulty for a moment and assume that the ubiquitous alternating least squares algorithm would yield the required solution.}%
, we will still get some sort of `approximate solution' because of the finite-precision error inherent in the computer. What really happens here \cite{W} is that we are effectively solving a problem perturbed by some small $\varepsilon>0$; the `approximate solution' $\xx_{i}^{\ast}(\varepsilon),\yy_{i}^{\ast}(\varepsilon)\in\Rr^{d_{i}}$ ($i=1,2,3$) is really a solution to the perturbed problem:%
\begin{multline}
\lVert A
-\xx_{1}^{\ast}(\varepsilon)\otimes\xx_{2}^{\ast}(\varepsilon)\otimes\xx_{3}^{\ast}(\varepsilon)
-\yy_{1}^{\ast}(\varepsilon)\otimes\yy_{2}^{\ast}(\varepsilon)\otimes\yy_{3}^{\ast}(\varepsilon)\rVert             \label{perturb}\\
=\varepsilon+\inf\nolimits_{\xx_{i},\yy_{i}\in\Rr^{d_{i}}}\lVert
A-\xx_{1}\otimes\xx_{2}\otimes\xx_{3}-\yy%
_{1}\otimes\yy_{2}\otimes\yy_{3}\rVert.
\end{multline}
Since we are attempting to find a solution of \eqref{nosoln} that does not exist, in exact arithmetic the algorithm will never terminate, but in reality the computer is limited by its finite precision and so the algorithm terminates at an `approximate solution', which may be viewed as a solution to a perturbed problem \eqref{perturb}. This process of forcing a solution to an ill-posed problem is almost always guaranteed to be ill-conditioned because of the infamous rule of thumb in numerical analysis \cite{Dem1, Dem2,Dem3}:
\medskip
\begin{center}
\framebox{
\textit{A well-posed problem near to an ill-posed one is ill-conditioned.}
}
\end{center}
\medskip
The root of the ill-conditioning lies in the fact that we are solving the (well-posed but ill-conditioned) problem \eqref{perturb} that is a slight perturbation of the ill-posed problem \eqref{nosoln}. The ill-conditioning manifests itself as the phenomenon described in Proposition~\ref{prop:unbounded}, namely,
\[
\lVert \xx_{1}^{\ast}(\varepsilon)\otimes\xx_{2}^{\ast}(\varepsilon)\otimes\xx_{3}^{\ast}(\varepsilon) \rVert \to \infty \qquad \text{and} \qquad
\lVert \yy_{1}^{\ast}(\varepsilon)\otimes\yy_{2}^{\ast}(\varepsilon)\otimes\yy_{3}^{\ast}(\varepsilon) \rVert \to \infty
\]
as $\varepsilon \to 0$. The ill-conditioning described here was originally observed in numerical experiments by psychometricians and chemometricians, who named the phenomenon  `diverging \textsc{candecomp}/\textsc{parafac} components' or `\textsc{candecomp}/\textsc{parafac} degeneracy' \cite{KDS, KHL, Paa, Steg1, Steg2}.

To fix the ill-conditioning, we should first fix the ill-posedness, i.e.\ find a well-posed problem. This leads us to the subject of the next section.

\subsection{Weak solutions}
\label{sect:weak}
In the study of partial differential equations~\cite{E}, there often
arise systems of \textsc{pde}s that have no solutions in the
traditional sense. A standard way around this is to define a
so-called \textit{weak solution}, which may not be a continuous
function or even a function (which is a tad odd since one would
expect a \textit{solution} to a \textsc{pde} to be at least
differentiable). Without going into the details, we will just say
that weak solution turns out to be an extremely useful concept and
is indispensable in modern studies of \textsc{pde}s. Under the
proper context, a weak solution to an ill-posed \textsc{pde} may be
viewed as the limit of \textit{strong} or \textit{classical
solutions} to a sequence of well-posed \textsc{pde} that are
slightly perturbed versions of the ill-posed one in question.
Motivated by the \textsc{pde} analogies, we will define weak
solutions to {\brap}.

We let $\Ss_{r}(d_{1},\dots,d_{k}):=\left\{ A\in\Rr^{d_{1}
\times\dots\times d_{k}}\mid\trank(A)\leq r\right\}$ and let
$\Ssbar_{r}(d_{1},\dots,d_{k})$ denote its closure in the (unique)
norm topology.

\begin{definition}
\label{def:brank}
An order-$k$ tensor $A\in\Rr^{d_{1} \times\dots\times d_{k}}$ has
\textbf{border rank} $r$ if
\[
A\in\overline {\Ss}_{r}(d_{1},\dots,d_{k}) \qquad \text{and} \qquad
A\not \in \overline{\Ss }_{r-1}(d_{1},\dots,d_{k}).
\]
This is denoted by $\brank(A)$. Note that
\[
\Ssbar_{r}(d_{1},\dots,d_{k})=\{A\in\Rr ^{d_{1}\times\dots\times
d_{k}}\mid \brank (S)\leq r\}.
\]
\end{definition}
\textit{Remark.} Clearly $\brank(A) \leq \trank(A)$ for any tensor~$A$. Since $\Ssbar_{0} = \Ss_{0}$ (trivially) and $\Ssbar_{1} = \Ss_{1}$ (by Proposition~\ref{prop:goodranka}), it follows that
$\brank(A) = \trank(A)$ whenever $\trank(A) \leq 2$. Moreover, $\brank(A) \geq 2$ if $\trank(A) \geq 2$.

Our definition differs slightly from the usual definition of border
rank in the algebraic computational complexity literature \cite{Bi1,
Bi2, BCS, Kn, Land}, which uses the Zariski topology (and is
normally defined for tensors over $\mathbb{C}$).

Let $A\in\Rr^{d_{1}\times\dots\times d_{k}}$ with $d_{i}\geq2$ and
$k\geq3$. Then the way to ensure that {\brap}, the optimal rank-$r$
approximation problem
\begin{equation}\label{eq:opt}
\operatorname*{argmin}\nolimits_{\trank(B)\leq
r}\lVert A-B\rVert
\end{equation}
always has a meaningful solution for any $A\in\Rr^{d_{1}\times
\dots\times d_{k}}$ is to instead consider the optimal
border-rank-$r$
approximation problem%
\begin{equation}\label{eq:optbar}
\operatorname*{argmin}\nolimits_{\brank(B)\leq
r}\lVert A-B\rVert.
\end{equation}

It is an obvious move to propose to fix the ill-posedness of {\brap} by
taking the closure. However, without a characterization of the limit points
such a proposal will at best be academic --- it is not enough to simply say that
weak solutions are limits of rank-$2$ tensors, without
giving an explicit expression (or a number of
expressions) for them that may be plugged into the objective
function to be minimized.

Theorem~\ref{thm:b} solves this problem in the order-$3$
rank-$2$ case --- it gives a complete description of these limit
points with an explicit formula and, in turn, a constructive
solution to the border-rank approximation problem.
In case this is not obvious, we will spell out the implication of
Theorem~\ref{thm:b}:

\begin{corollary}\label{cor:b}
Let $d_{1},d_{2},d_{3}\geq2$. Let $A \in\Rr^{d_{1}\times d_{2}\times
d_{3}}$ with $\trank(A)=3$. $A$ is the limit of a
sequence $A_{n}\in\Rr^{d_{1}\times d_{2}\times d_{3}}$ with
$\trank(A_{n})\leq2$ if and only if%
\[
A =
\yy_{1}\otimes \xx_{2}\otimes \xx_{3} +
\xx_{1} \otimes \yy_{2}\otimes \xx_{3} +
\xx_{1} \otimes \xx_{2}\otimes \yy_{3}
\]
for some $\xx_{i},\yy_{i}$ linearly independent vectors in
$\Rr^{d_{i}}$, $i=1,2,3$.
\end{corollary}

This implies that every tensor in $\Ssbar_{2}(d_{1},\dots,d_{k})$
can be written in one of two forms:
\begin{equation}\label{rank2form}
\yy_{1}\otimes \xx_{2}\otimes \xx_{3} +
\xx_{1} \otimes \yy_{2}\otimes \xx_{3} +
\xx_{1} \otimes \xx_{2}\otimes \yy_{3}
\end{equation}
or
\begin{equation}\label{dSLform}
\xx_{1}\otimes\xx_{2}\otimes\xx_{3}+\yy_{1}\otimes\yy_{2}\otimes\yy_{3}.
\end{equation}
These expressions may then be used to define the relevant objective
function(s) in the minimization of \eqref{eq:optbar}. As in the case
of \textsc{pde}, every classical (strong) solution is also a weak
solution to {\brap}.
\begin{proposition}
If $B$ is a solution to~\eqref{eq:opt} then $B$ is a solution to~\eqref{eq:optbar}.
\end{proposition}
\begin{proof}
If $\| A-B \| \leq \| A - B' \|$ for all $B' \in \Ss_r$, then
$\|A-B\| \leq \|A-B'\|$ for all $B' \in \Ssbar_r$ by continuity.
\end{proof}

\section{Semialgebraic description of tensor rank\label{sect:SA}}

One may wonder whether the result in Propositions~\ref{prop:rank2:b} and
\ref{prop:rank2:a} extends to more general \textit{hyperdeterminants}. We know
from \cite{GKZ1, GKZ2} that a hyperdeterminant may be uniquely defined (up to
a constant scaling) in $\Rr^{d_{1}\times\dots\times d_{k}}$ whenever
$d_{1},\dots,d_{k}$ satisfy%
\begin{equation}
d_{i}-1\leq\sum_{j\neq i}(d_{j}-1)\qquad\text{for }i=1,\dots,k.
\label{DetCond}%
\end{equation}
(Note that for matrices, \eqref{DetCond} translates to $d_{1}=d_{2}$, which
may be viewed as one reason why the determinant is defined only for square
matrices). Let $\Det_{d_{1},\dots,d_{k}}:\Rr^{d_{1}%
\times\dots\times d_{k}}\rightarrow\Rr$ be the polynomial function
defined by the hyperdeterminant, whenever \eqref{DetCond} is satisfied.
Propositions \ref{prop:rank2:b} and~\ref{prop:rank2:a} tell us that the rank
of a tensor is~$2$ on the set $\{A\mid\Det_{2,2,2}(A)>0\}$ and
$3$ on the set $\{A\mid\Det_{2,2,2}(A)<0\}$. One may start by
asking whether the tensor rank in $\Rr^{d_{1}\times\dots\times d_{k}}$
is constant-valued on the sets
\[
\{A\mid\Det_{d_{1},\dots,d_{k}}(A)<0\}\qquad\text{and}%
\qquad\{A\mid\Det_{d_{1},\dots,d_{k}}(A)>0\}.
\]
The answer, as Sturmfels has kindly communicated to us \cite{St}, is
\textit{no} with explicit counterexamples in cases
$2\times2\times2\times2$ and $3\times3\times3$. We will not
reproduce Sturmfels' examples here (one reason is that
$\Det_{2,2,2,2}$ already contains close to $3$ million monomial
terms \cite{GHSY}) but instead refer our readers to his forthcoming
paper.

We will prove that although there is no single polynomial $\Delta$
that will separate $\mathbb{R}^{d_1 \times \dots \times d_k}$ into
regions of constant rank as in the case of $\mathbb{R}^{2 \times 2
\times 2}$, there is always a finite number of polynomials
$\Delta_1,\dots, \Delta_m$ that will achieve this.

Before we state and prove the result, we will introduce a few
notions and notations. We will write $\mathbb{R}[X_1,\dots,X_m]$ for
the ring of polynomials in $m$ variables $X_1,\dots,X_m$ with real
coefficients. Subsequently, we will be considering polynomial
functions on tensor spaces and will index our variables in a
consistent way (for example, when discussing polynomial functions on
$\mathbb{R}^{l\times m \times n}$, the polynomial ring in question
will be denoted $\mathbb{R}[X_{111}, X_{112}, \dots X_{lmn}]$).
Given $A = \llbracket a_{ijk} \rrbracket \in \mathbb{R}^{l\times m
\times n}$ and $p(X_{111}, X_{112}, \dots X_{lmn}) \in
\mathbb{R}[X_{111}, X_{112}, \dots X_{lmn}]$, $p(A)$ will mean the
obvious thing, namely, $p(A) = p(a_{111}, a_{112},\dots, a_{lmn})
\in \mathbb{R}$.

A \textit{polynomial map} is a function $F : \mathbb{R}^n \to
\mathbb{R}^m$, defined for each $\mathbf{a} =
[a_1,\dots,a_n]^{\transp} \in \mathbb{R}^n$, by $F(\mathbf{a}) =
[f_1(\mathbf{a}),\dots, f_m(\mathbf{a})]^{\transp}$ where $f_i
\in \mathbb{R}[X_1,\dots, X_n]$ for all $i=1,\dots,m$.

A \textit{semialgebraic set} in $\mathbb{R}^n$ is a union of
finitely many sets of the form\footnote{%
Only one $p$ is necessary, because multiple equality constraints $p_1(\mathbf{a})=0$, \dots, $p_k(\mathbf{a})=0$ can always be amalgamated into a single equation $p(\mathbf{a}) = 0$ by setting {$p = p_1^2 + \dots + p_k^2$}.
}
\[
\{ \mathbf{a} \in \mathbb{R}^n \mid p(\mathbf{a}) = 0,\,
q_1(\mathbf{a})> 0,\,\dots,\, q_\ell(\mathbf{a}) > 0 \}
\]
where $\ell \in \mathbb{N}$ and $p,q_1,\dots,q_\ell \in
\mathbb{R}[X_1,\dots, X_n]$.
%
%
Note that we do not exclude the
possibility of $p$ or any of the $q_i$ being constant (degree-$0$)
polynomials. For example, if $p$ is the zero polynomial, then the
first relation $0 = 0$ is trivially satisfied and the semialgebraic
set will be an open set in $\mathbb{R}^n$.

It is easy to see that the class of all semialgebraic sets in
$\mathbb{R}^n$ is closed under finite unions, finite intersections,
and taking complement. Moreover, if $\mathcal{S} \subseteq
\mathbb{R}^{n+1}$ is a semialgebraic set and $\pi : \mathbb{R}^{n+1}
\to \mathbb{R}^n$ is the projection onto the first $n$ coordinates,
then $\pi (\mathcal{S})$ is also a semialgebraic set --- this
seemingly innocuous statement is in fact the Tarski--Seidenberg
theorem~\cite{Sei, Tar}, possibly the most celebrated result about
semialgebraic sets. We will restate it in a (somewhat less common)
form that better suits our purpose.

\begin{theorem}[Tarski--Seidenberg]\label{thm:tarski}
If $\mathcal{S} \subseteq \mathbb{R}^n$ is a semialgebraic set and
$F : \mathbb{R}^n \to \mathbb{R}^m$ is a polynomial map, then the
image $F(\mathcal{S}) \subseteq \mathbb{R}^m$ is also a
semialgebraic set.
\end{theorem}

These and other results about semialgebraic sets may be found in
\cite[Chapter~2]{Co}, which, in addition, is a very readable
introduction to semialgebraic geometry.

\begin{theorem}\label{thm:sa}
The set $\mathcal{R}_{r}(d_{1},\dots,d_{k}):=\{
A\in\Rr^{d_{1}\times\dots\times d_{k}}\mid\trank(A) = r\}$ is a
semialgebraic set.
\end{theorem}
\begin{proof}
Let $\psi_{r}:(\Rr^{d_{1}}\times\Rr^{d_{1}}\times\dots
\times\Rr^{d_{k}})^{r}\rightarrow\Rr^{d_{1}\times d_{2}%
\times\dots\times d_{k}}$ be defined by
\[
\psi_{r}(\uu_{1},\vv_{1},\dots,\mathbf{z}_{1};\dots
;\uu_{r},\vv_{r},\dots,\mathbf{z}_{r})=\uu_{1}%
\otimes\vv_{1}\otimes\dots\otimes\mathbf{z}_{1}+\dots+\uu%
_{r}\otimes\vv_{r}\otimes\dots\otimes\mathbf{z}_{r}.
\]
It is clear that the image of $\psi_{r}$ is exactly
$\mathcal{S}_{r}(d_{1},\dots,d_{k}) =
\{A\mid\operatorname*{rank}_{\otimes }(A)\leq r\}$. It is also clear
that $\psi_{r}$ is a polynomial map.

It follows from Theorem~\ref{thm:tarski} that
$\mathcal{S}_{r}(d_{1},\dots,d_{k})$ is semialgebraic. This holds
for arbitrary $r$. So
$\mathcal{R}_{r}(d_{1},\dots,d_{k})=\mathcal{S}_{r}(d_{1},\dots,d_{k})\setminus\mathcal{S}_{r-1}(d_{1},\dots,d_{k})$
is also semialgebraic.
\end{proof}

\begin{corollary}
\label{cor:sa}
There exist $\Delta_0, \dots ,\Delta_m \in \mathbb{R}[X_{1\cdots
1},\dots, X_{d_1\cdots d_k}]$ from which the rank of a tensor $A \in
\mathbb{R}^{d_1 \times \dots \times d_k}$ can be determined purely
from the signs (i.e.\ $+$ or $-$ or $0$) of $\Delta_0(A), \dots,
\Delta_m(A)$.
\end{corollary}

In the next section, we will see examples of such polynomials for
the tensor space $\Rr^{2\times 2 \times 2}$. We will stop short of
giving an explicit semialgebraic characterization of rank, but it
should be clear to the reader how to get one.

\section{Orbits of real $2\times2\times2$ tensors}
\label{sect:orbits}

In this section, we study the equivalence of tensors in $\Rr^{2\times2\times2}$ under multilinear matrix multiplication.
We will use the results and techniques of this section later on in
Section~\ref{sec:volume} where we determine \emph{which} tensors
in~$\Rr^{2 \times 2 \times 2}$ have an optimal rank-$2$
approximation.

Recall that $A$ and $B\in\Rr^{2\times2\times2}$ are said to be
($\GL_{2,2,2}(\Rr)$-)equivalent iff there exists a transformation
$(L,M,N)\in\GL_{2,2,2}(\Rr)$ such that $A=(L,M,N)\cdot B$. The
question is whether there is a finite list of `canonical tensors' so
that every $A\in\Rr^{2\times 2\times 2}$ is equivalent to one of
them. For matrices, $A\in\Rr^{m\times n}$,
$\operatorname*{rank}(A)=r$ if and only if there exists $M\in\GL
_{m}(\Rr),N\in\GL_{n}(\Rr)$ such that
\[
(M,N)\cdot A=MAN^{\transp}=
\begin{bmatrix}
I_{r} & 0\\
0 & 0
\end{bmatrix}
.
\]
So every matrix of rank~$r$ is equivalent to one that takes the
canonical form $ \left[ \begin{smallmatrix}
I_{r} & 0\\
0 & 0
\end{smallmatrix}\right]
$. Note that this is the same as saying that the matrix $A$ can be transformed
into $%
\left[ \begin{smallmatrix}
I_{r} & 0\\
0 & 0
\end{smallmatrix}\right]
$ using elementary row- and column-operations: adding a scalar
multiple of a row/column to another, scaling a row/column by a
non-zero scalar, interchanging two rows/columns --- since every
$(L_{1},L_{2})\in \GL_{m,n}(\Rr)$ is a sequence of such operations.

We will see that there is indeed a finite number of canonical forms for
tensors in~$\Rr^{2\times2\times2}$; although the
classification is somewhat more intricate than the case of matrices --- two
tensors in $\Rr^{2\times2\times2}$ can have the same rank but be
inequivalent (i.e.\ reduce to different canonical forms).

In fancier language, what we are doing is classifying the
\textit{orbits} of the \textit{group action} $\GL_{2,2,2}(\Rr)$ on
$\Rr^{2\times2\times2}$. We are doing for $\Rr^{2\times2\times 2}$
what Gelfand, Kapranov, and Zelevinsky did for $\mathbb{C}^{2\times
2\times2}$ in the last sections of \cite{GKZ1, GKZ2}. Not
surprisingly, the results that we obtained are similar but not
identical --- there are \textit{eight} distinct orbits for the
action of $\GL_{2,2,2}(\Rr)$ on $\Rr^{2\times2\times2}$ as opposed
to \textit{seven} distinct orbits for the
action of $\GL_{2,2,2}(\mathbb{C})$ on $\mathbb{C}%
^{2\times2\times2}$ --- a further reminder of the dependence of such results
on the choice of field.


\begin{theorem}
\label{thm:orbits}
Every tensor in $\Rr^{2\times2\times2}$ is equivalent
via a transformation in $\GL_{2,2,2}(\Rr)$ to precisely one of the canonical forms indicated in Table~\ref{table:orbits}, with its invariants taking the values shown.
\begin{table}
\[%
\begin{array}
[c]{ccccccc}%
\text{\rm tensor}
& \quad & \operatorname{sign}(\Delta) & \multirank
& \quad & \trank & \brank
\\
\hline
&  &  &  &  &  & \\
D_{0}=\left[
\begin{array}[c]{rr} 0 & 0\\ 0 & 0 \end{array}
\right\vert \!\left.
\begin{array}[c]{rr} 0 & 0\\ 0 & 0 \end{array}
\right]  &  & 0 & (0,0,0) & & 0 & 0\\
&  &  &  &  &  & \\
D_{1}=\left[
\begin{array}[c]{rr} 1 & 0\\ 0 & 0 \end{array}
\right\vert \!\left.
\begin{array}[c]{rr} 0 & 0\\ 0 & 0 \end{array}
\right]  &  & 0 & (1,1,1) &  & 1 & 1\\
&  &  &  &  &  & \\
D_{2}=\left[
\begin{array}[c]{rr} 1 & 0\\ 0 & 1 \end{array}
\right\vert \!\left.
\begin{array}[c]{rr} 0 & 0\\ 0 & 0 \end{array}
\right]  &  & 0 & (1,2,2) & & 2 & 2\\
&  &  &  &  &  & \\
D_{2}^{\prime}=\left[
\begin{array}[c]{rr} 1 & 0\\ 0 & 0 \end{array}
\right\vert \!\left.
\begin{array}[c]{rr} 0 & 1\\ 0 & 0 \end{array}
\right]  &  & 0 & (2,1,2) & & 2 & 2\\
&  &  &  &  &  & \\
D_{2}^{\prime\prime}=\left[
\begin{array}[c]{rr} 1 & 0\\ 0 & 0 \end{array}
\right\vert \!\left.
\begin{array}[c]{rr} 0 & 0\\ 1 & 0 \end{array}
\right]  &  & 0 & (2,2,1) & & 2 & 2\\
&  &  &  &  &  & \\
G_{2}=\left[
\begin{array}[c]{rr} 1 & 0\\ 0 & 0 \end{array}
\right\vert \!\left.
\begin{array}[c]{rr} 0 & 0\\ 0 & 1 \end{array}
\right]  &  & + & (2,2,2) & & 2 & 2\\
&  &  &  &  &  & \\
D_{3}=\left[
\begin{array}[c]{rr} 1 & 0\\ 0 & 0 \end{array}
\right\vert \!\left.
\begin{array}[c]{rr} 0 & 1\\ 1 & 0 \end{array}
\right]  &  & 0 & (2,2,2) & & 3 & 2\\
&  &  &  &  &  & \\
G_{3}=\left[
\begin{array}[c]{rr} 1 & 0\\ 0 & 1 \end{array}
\right\vert \!\left.
\begin{array} [c]{rr} 0 & -1\\ 1 & 0 \end{array}
\right]  &  & - & (2,2,2) & & 3 & 3\\
&  &  &  &  &  & \\\hline
\end{array}
\]
\caption{$\GL$-orbits of~$\Rr^{2\times 2\times 2}$.
The letters $D, G$ stand for `degenerate' and `generic' respectively.
\label{table:orbits}
}
\end{table}
%
\end{theorem}

\begin{proof}
Write $A=[A_{1}\mid A_{2}]$, $A_{i}\in\Rr^{2\times2}$, for
$\llbracket a_{ijk}\rrbracket \in\Rr^{2\times2\times2}$. If
$\rank(A_{1})=0$, then {
\[
A=\left[
\begin{array}[c]{rr} 0 & 0\\ 0 & 0 \end{array}
\right\vert \!\left.
\begin{array}[c]{rr} \times & \times\\ \times & \times \end{array}
\right]  .
\]
} Using matrix operations, $A$ must then be equivalent to one of the following
forms (depending on $\rank(A_{2})$) {
\[
\left[
\begin{array}[c]{rr} 0 & 0\\ 0 & 0 \end{array}
\right\vert \!\left.
\begin{array}[c]{rr} 0 & 0\\ 0 & 0 \end{array}
\right]  ,\qquad\left[
\begin{array}[c]{rr} 0 & 0\\ 0 & 0 \end{array}
\right\vert \!\left.
\begin{array}[c]{rr} 1 & 0\\ 0 & 0 \end{array}
\right]  ,\qquad\left[
\begin{array}[c]{rr} 0 & 0\\ 0 & 0 \end{array}
\right\vert \!\left.
\begin{array}[c]{rr} 1 & 0\\ 0 & 1 \end{array}
\right]  ,
\]
} which correspond to $D_{0}$, $D_{1}$, $D_{2}$ respectively (after reordering
the slabs).

If $\rank(A_{1})=1$, then we may assume that {
\[
A=\left[
\begin{array}[c]{rr} 1 & 0\\ 0 & 0 \end{array}
\right\vert \!\left.
\begin{array}[c]{rr} a & b\\ c & d \end{array}
\right]  .
\]
} If $d\not =0$ then we may transform this to~$G_{2}$ as follows:
\[
\left[
\begin{array}[c]{rr} 1 & 0\\ 0 & 0 \end{array}
\right\vert \!\left.
\begin{array}[c]{rr} a & b\\ c & d \end{array}
\right]  \leadsto\left[
\begin{array}[c]{rr} 1 & 0\\ 0 & 0 \end{array}
\right\vert \!\left.
\begin{array}[c]{rr} \times & 0\\ 0 & d \end{array}
\right]  \leadsto\left[
\begin{array}[c]{rr} 1 & 0\\ 0 & 0 \end{array}
\right\vert \!\left.
\begin{array}[c]{rr} 0 & 0\\ 0 & 1 \end{array}
\right]
\]
If $d=0$ then:
\[
\left[
\begin{array}
[c]{rr} 1 & 0\\ 0 & 0 \end{array}
\right\vert \!\left.
\begin{array}[c]{rr} a & b\\ c & 0 \end{array}
\right]  \leadsto\left[
\begin{array}
[c]{rr} 1 & 0\\ 0 & 0 \end{array}
\right\vert \!\left.
\begin{array}[c]{rr} 0 & b\\ c & 0 \end{array}
\right]
\]
In this situation we can normalize $b,c$ separately, reducing these matrices
to one of the following four cases (according to whether $b,c$ are zero):
\[
\left[
\begin{array}[c]{rr} 1 & 0\\ 0 & 0 \end{array}
\right\vert \!\left.
\begin{array}[c]{rr} 0 & 0\\ 0 & 0 \end{array}
\right]  ,\quad\left[
\begin{array}[c]{rr} 1 & 0\\ 0 & 0 \end{array}
\right\vert \!\left.
\begin{array}[c]{rr} 0 & 1\\ 0 & 0 \end{array}
\right]  ,\quad\left[
\begin{array}[c]{rr} 1 & 0\\ 0 & 0 \end{array}
\right\vert \!\left.
\begin{array}[c]{rr} 0 & 0\\ 1 & 0 \end{array}
\right]  ,\quad\left[
\begin{array} [c]{rr} 1 & 0\\ 0 & 0 \end{array}
\right\vert \!\left.
\begin{array}[c]{rr} 0 & 1\\ 1 & 0 \end{array}
\right]  ,
\]
which are $D_{1}$, $D_{2}^{\prime}$, $D_{2}^{\prime\prime}$, $D_{3}$ respectively.

Finally, if $\rank(A_{1})=2$, then we may assume that {
\[
A=[A_{1}\mid A_{2}]=\left[
\begin{array}
[c]{rr}%
1 & 0\\
0 & 1
\end{array}
\right\vert \!\left.
\begin{array}
[c]{rr}%
\times & \times\\
\times & \times
\end{array}
\right]  .
\]
} By applying a transformation of the form $(I,L,L^{-1})$, we can keep $A_{1}$
fixed while conjugating~$A_{2}$ into (real) Jordan canonical form. There are
four cases.

If $A_{2}$ has repeated real eigenvalues and is diagonalizable, then we
get~$D_{2}$:
\[
\left[
\begin{array}
[c]{rr}%
1 & 0\\
0 & 1
\end{array}
\right\vert \!\left.
\begin{array}
[c]{rr}%
\lambda & 0\\
0 & \lambda
\end{array}
\right]  \leadsto\left[
\begin{array}
[c]{rr}%
1 & 0\\
0 & 1
\end{array}
\right\vert \!\left.
\begin{array}
[c]{rr}%
0 & 0\\
0 & 0
\end{array}
\right]
\]

If $A_{2}$ has repeated real eigenvalues and is not diagonalizable, then we
have
\[
\left[
\begin{array}
[c]{rr}%
1 & 0\\
0 & 1
\end{array}
\right\vert \!\left.
\begin{array}
[c]{rr}%
\lambda & 1\\
0 & \lambda
\end{array}
\right]  \leadsto\left[
\begin{array}
[c]{rr}%
1 & 0\\
0 & 1
\end{array}
\right\vert \!\left.
\begin{array}
[c]{rr}%
0 & 1\\
0 & 0
\end{array}
\right]  ,
\]
which is equivalent (after swapping columns and swapping slabs) to~$D_{3}$.

If $A_{2}$ has distinct real eigenvalues, then $A$ reduces to~$G_{2}$:
\[
\left[
\begin{array}
[c]{rr}%
1 & 0\\
0 & 1
\end{array}
\right\vert \!\left.
\begin{array}
[c]{rr}%
\lambda & 0\\
0 & \mu
\end{array}
\right]  \leadsto\left[
\begin{array}
[c]{cc}%
1 & 0\\
0 & 1
\end{array}
\right\vert \!\left.
\begin{array}
[c]{cc}%
0 & 0\\
0 & \mu-\lambda
\end{array}
\right]  \leadsto\left[
\begin{array}
[c]{rr}%
1 & 0\\
0 & 0
\end{array}
\right\vert \!\left.
\begin{array}
[c]{rr}%
0 & 0\\
0 & 1
\end{array}
\right]
\]

If $A_{2}$ has complex eigenvalues, then we can reduce $A$ to~$G_{3}$:
\[
\left[
\begin{array}
[c]{rr}%
1 & 0\\
0 & 1
\end{array}
\right\vert \!\left.
\begin{array}
[c]{rr}%
a & -b\\
b & a
\end{array}
\right]  \leadsto\left[
\begin{array}
[c]{rr}%
1 & 0\\
0 & 1
\end{array}
\right\vert \!\left.
\begin{array}
[c]{rr}%
0 & -b\\
b & 0
\end{array}
\right]  \leadsto\left[
\begin{array}
[c]{rr}%
1 & 0\\
0 & 1
\end{array}
\right\vert \!\left.
\begin{array}
[c]{rr}%
0 & -1\\
1 & 0
\end{array}
\right]
\]

Thus, every $2\times2\times2$ tensor can be transformed to one of the
canonical forms listed in the statement of the theorem. Moreover, the
invariants $\operatorname{sign}(\Delta)$ and~$\multirank$
are easily computed for the canonical forms, and suffice to distinguish them.
It follows that the listed forms are pairwise inequivalent.

We confirm the given values of $\trank$. It is clear that $\trank(D_{0}) = 0$ and $\trank(D_{1}) = 1$. By Proposition~\ref{prop:reduction:1}, any tensor of
rank~$1$ must be equivalent to~$D_{1}$. Thus $D_{2}$, $D_{2}^{\prime}$,
$D_{2}^{\prime\prime}$ and $G_{2}$ are all of rank~$2$. By
Proposition~\ref{thm:reduction:2}, every tensor of rank~$2$ must be equivalent
to one of these. In particular, $D_{3}$ and $G_{3}$ must have rank at least~3.
Evidently $\trank(D_{3})=3$ from its definition; and
the same is true for~$G_{3}$ by virtue of the less obvious relation
\[
G_{3} = (\ee_{1}+\ee_{2})\otimes\ee_{2} \otimes
\ee_{2} + (\ee_{1}-\ee_{2})\otimes\ee_{1}
\otimes\ee_{1} + \ee_{2} \otimes(\ee_{1}+\ee_{2})
\otimes(\ee_{1}-\ee_{2}) .
\]

Finally, we confirm the tabulated values of $\brank$. By virtue of the remark after Definition~\ref{def:brank}, it is enough to verify that $\brank(D_3) \leq 2$ and that $\brank(G_3) = 3$. The first of these assertions follows from Proposition~\ref{prop:dSLtensor}. The set of tensors of type~$G_3$ is an open set, which implies the second assertion.
\end{proof}

\textit{Remark.} We note that $D_3$ is equivalent to any of the
tensors obtained from it by permutations of the three factors.
Indeed, all of these tensors have $\multirank = (2,2,2)$ and $\Delta
= 0$. Similar remarks apply to $G_2$, $G_3$.

\textit{Remark.} The classification of $\GL_{2,2,2}(\mathbb{C}%
)$-orbits in $\mathbb{C}^{2\times2 \times2}$ differs only in the treatment
of~$G_{3}$, since there is no longer any distinction between real and complex eigenvalues.

We caution the reader that the \emph{finite} classification in
Theorem \ref{thm:orbits} is, in general, not possible for tensors of
arbitrary size and order simply because the dimension or `degrees of freedom' of $\Rr%
^{d_{1}\times\dots\times d_{k}}$ exceeds that of $\GL%
_{d_{1},\dots,d_{k}}(\Rr)$ as soon as $d_{1}\cdots d_{k}>d_{1}%
^{2}+\cdots+d_{k}^{2}$ (which is almost always the case). Any
attempt at an explicit classification must necessarily include
continuous parameters. For the case of $\Rr^{2\times2\times2}$ this
argument is not in conflict with our finite classification, since
$2\cdot2\cdot2<2^{2}+2^{2}+2^{2}$.

\subsection{{Generic rank}}
We called the tensors in the orbit classes of $G_2$ and $G_3$
\textit{generic} in the sense that the property of being in either
one of these classes is an open condition. One should note that
there is often no one single \textit{generic {outer product} rank} for tensors over
$\mathbb{R}$ \cite{Krus1,tBK}. (For tensors over $\mathbb{C}$ such a generic rank always exists \cite{CGLM2}.) The `generic {outer product} rank' for tensors over $\mathbb{R}$
should be regarded as set-valued:
\[
\operatorname*{generic-rank}\nolimits_{\otimes}(\mathbb{R}^{d_1 \times \dots \times d_k})
= \{r \in \mathbb{N} \mid \Ss_r(d_1,\dots,d_k) \text{ has non-empty interior} \}.
\]
So the generic outer product rank in $\mathbb{R}^{2\times 2 \times 2}$ is
$\{2,3\}$. Another term, preferred by some and coined originally by
ten Berge, is \textit{typical rank} \cite{tBK}.

Given $d_1,\dots d_k$, the determination of the generic outer product rank is an open problem in general and a nontrivial problem even in simple cases --- see \cite{CGG1, CGG2} for results over $\mathbb{C}$ and \cite{tB, tBK} for results over $\mathbb{R}$. Fortunately, the difficulty does not extend to multilinear rank --- a single unique \textit{generic multilinear rank} always exist and depends only on $d_1,\dots d_k$ (and not on the base field, cf.\ Proposition \ref{prop:indmrank}).
\begin{proposition}\label{prop:gmrank} Let $A \in \mathbb{R}^{d_1 \times \dots \times d_k}$. If $\multirank(A) = (r_1 (A),\dots, r_k(A))$, then
\[
  r_i(A) = \min \Bigl( d_i, \prod\nolimits_{j\ne i} d_j \Bigr),\quad i = 1,\dots,k,
\]
generically.
\end{proposition}
\begin{proof}
Let $\mu_i : \mathbb{R}^{d_1 \times \dots \times d_k} \to \mathbb{R}^{d_i \times \prod_{j \ne i} d_j}$ be the forgetful map that `flattens' or `unfolds' a tensor into a matrix in the $i$th mode. It is easy to see that
\begin{equation}\label{eq:matrank}
r_i(A) = \rank (\mu_i (A))
\end{equation}
where `$\rank$' here denotes matrix rank. The results then follow from the fact that the generic rank of a matrix in $\mathbb{R}^{d_i \times \prod_{j \ne i} d_j}$ is $\min (d_i, \prod_{j\ne i} d_j)$.
\end{proof}

For example, for order-$3$ tensors,
\[
\operatorname*{generic-rank}\nolimits_{\boxplus}(\mathbb{R}^{l \times m \times n}) = (\min(l, mn), \min(m, ln), \min(n, lm)).
\]

\subsection{Semialgebraic description of orbit classes}
For a general tensor $A \in\Rr^{2\times2\times2}$, its orbit class is
readily determined by computing the invariants $\operatorname{sign}(\Delta(A))$ and $\multirank(A)$, and comparing with the canonical forms. The ranks
$r_{i}(A)$ which constitute $\multirank(A)$ can be evaluated
algebraically as follows. If $A \not = 0$ then each~$r_{i}(A)$ is
either 1 or~2. For example, note that $r_{1}(A) < 2$ if and only if
the vectors $A_{\bullet 1 1}$, $A_{\bullet 1 2}$, $A_{\bullet 2 1}$,
$A_{\bullet 2 2}$ are linearly dependent, which happens if and only
if all the 2-by-2 minors of the matrix
\[%
\begin{bmatrix}
a_{111} & a_{112} & a_{121} & a_{122}\\
a_{211} & a_{212} & a_{221} & a_{222}%
\end{bmatrix}
\]
are zero. Explicitly, the following six equations must be satisfied:
\begin{align}
\label{eq:minor1}a_{111}a_{212}  &  =a_{211}a_{112},\quad a_{111}%
a_{221}=a_{211}a_{121},\quad a_{111}a_{222}=a_{211}a_{122},\\
a_{112}a_{221}  &  =a_{212}a_{121},\quad a_{112}a_{222}=a_{212}a_{122},\quad
a_{121}a_{222}=a_{221}a_{122}.\nonumber
\end{align}
Similarly, $r_{2}(A) < 2$ if and only if
\begin{align}
\label{eq:minor2}a_{111}a_{122}  &  =a_{121}a_{112},\quad a_{111}%
a_{221}=a_{121}a_{211},\quad a_{111}a_{222}=a_{121}a_{212},\\
a_{112}a_{122}  &  =a_{122}a_{211},\quad a_{112}a_{222}=a_{122}a_{212},\quad
a_{211}a_{222}=a_{221}a_{212};\nonumber
\end{align}
and $r_{3}(A) < 2$ if and only if
\begin{align}
\label{eq:minor3}a_{111}a_{122}  &  =a_{112}a_{121},\quad a_{111}%
a_{212}=a_{112}a_{211},\quad a_{111}a_{222}=a_{112}a_{221},\\
a_{121}a_{212}  &  =a_{122}a_{211},\quad a_{121}a_{222}=a_{122}a_{221},\quad
a_{211}a_{222}=a_{212}a_{221}.\nonumber
\end{align}
The equations \eqref{eq:minor1}--\eqref{eq:minor3} lead to twelve distinct polynomials (beginning with $\Delta_1 = a_{111}a_{212} - a_{211}a_{112}$) which, together with $\Delta_0 := \Delta$, make up the collection $\Delta_0,\dots,\Delta_{12}$ of polynomials used in the semialgebraic description of the orbit structure of $\mathbb{R}^{2 \times 2 \times 2}$, as in Corollary~\ref{cor:sa}. Indeed,
we note that in Table \ref{table:orbits} the information in the
fourth and fifth columns ($\trank(A)$, $\brank(A)$) is determined by the information in the second and third columns ($\operatorname{sign}(\Delta)$, $\multirank(A)$).

\subsection{Generic rank on $\Delta = 0$}

The notion of generic rank also makes sense on subvarieties of~$\Rr^{2\times2\times2}$, for instance on the $\Delta=0$ hypersurface.

\begin{proposition}
The tensors on the hypersurface $\Dd_3 = \{A\in\Rr^{2\times2\times2}\mid
\Delta(A)=0\}$ are all of the form
\[
\xx_{1}\otimes\xx_{2}\otimes\yy_{3}+\xx_{1}%
\otimes\yy_{2}\otimes\xx_{3}+\yy_{1}\otimes\xx%
_{2}\otimes\xx_{3}%
\]
and they have rank~$3$ generically.
\end{proposition}


\begin{proof}
From the canonical forms in Table~\ref{table:orbits}, we see that if $\Delta(A)=0$, then%
\[
A=\xx_{1}\otimes\xx_{2}\otimes\yy_{3}+\xx%
_{1}\otimes\yy_{2}\otimes\xx_{3}+\yy_{1}\otimes
\xx_{2}\otimes\xx_{3}%
\]
for some $\xx_{i},\yy_{i}\in\Rr^{2}$,
not necessarily linearly independent. It
remains to be shown that $\trank(A)=3$ generically.

From Theorem~\ref{thm:orbits} and the subsequent discussion, if $\Delta(A)=0$
then $\trank(A) \leq2$ if and only if at least one of
the equation sets \eqref{eq:minor1}, \eqref{eq:minor2}, \eqref{eq:minor3} is
satisfied. Hence $\Dd_{2} := \{A\mid\Delta(A)=0,\,\trank(A)\leq2\}$ is an algebraic subset of $\Dd_{3}$.

On the other hand, $\Dd_{3}\setminus \Dd_{2}$ is dense
in~$\Dd_{3}$ with respect to the Euclidean, and hence the Zariski,
topology. Indeed, each of the tensors $D_{0}$, $D_{1}$, $D_{2}$,
$D_{2}^{\prime}$, $D_{2}^{\prime\prime}$ can be approximated by tensors of
type~$D_{3}$; for instance
\[
\left[
\begin{array}[c]{cc} 1 & 0\\ 0 & 1 \end{array}
\right\vert \!\left.
\begin{array}[c]{cc} 0 & \epsilon\\ 0 & 0 \end{array}
\right]
\to
\left[
\begin{array}[c]{cc} 1 & 0\\ 0 & 1 \end{array}
\right\vert \!\left.
\begin{array}[c]{cc} 0 & 0\\ 0 & 0 \end{array}
\right]  =D_{2},\quad\text{as }\epsilon\to0.
\]
Multiplying by an arbitrary $(L,M,N)\in\GL_{2,2,2}(\Rr)$, it
follows that any tensor in~$\Dd_{2}$ can be approximated by tensors of
type~$D_{3}$.

It follows that the rank-$3$ tensors $\Dd_{3} \setminus \Dd_{2}$ in~$\Dd_{3}$ constitute a generic subset of~$\Dd_{3}$, in the Zariski sense (and hence in all the other usual senses).
\end{proof}

\textit{Remark.} In fact, it can be shown that $\Dd_{3}$ is an irreducible variety. If we accept that, then the fact that $\Dd_{2}$ is a proper subvariety of~$\Dd_{3}$ immediately implies that the rank-$3$ tensors form a generic subset of~$\Dd_{3}$. The denseness argument becomes unnecessary.

\subsection{Base field dependence}

It is interesting to observe that the $\GL_{2,2,2}(\Rr)$-orbit classes of $G_2$ and $G_3$ merge into a single orbit class over $\mathbb{C}$ (under the action of $\GL_{2,2,2}(\Cc)$). Explicitly, if we write $\mathbf{z}_k = \mathbf{x}_k + i\mathbf{y}_k$ and  $\bar{\mathbf{z}}_k = \mathbf{x}_k - i\mathbf{y}_k$, then
\begin{multline}\label{eq:fielddep}
\mathbf{x}_1\otimes\mathbf{x}_2\otimes\mathbf{x}_3
+\mathbf{x}_1\otimes\mathbf{y}_2\otimes\mathbf{y}_3
-\mathbf{y}_1\otimes\mathbf{x}_2\otimes\mathbf{y}_3
+\mathbf{y}_1\otimes\mathbf{y}_2\otimes\mathbf{x}_3\\
=\frac{1}{2}(
\bar{\mathbf{z}}_1\otimes\mathbf{z}_2\otimes\bar{\mathbf{z}}_3
+\mathbf{z}_1\otimes\bar{\mathbf{z}}_2\otimes\mathbf{z}_3
).
\end{multline}
The \textsc{lhs} is in the $\GL_{2,2,2}(\Rr)$-orbit class of $G_3$ and has outer product rank $3$ over $\mathbb{R}$. The \textsc{rhs} is in the $\GL_{2,2,2}(\Cc)$-orbit class of $G_2$ and has outer product rank $2$ over $\mathbb{C}$. To see why this is unexpected, recall that an $m \times n$ matrix with real entries has the same rank whether we regard it as an element of $\mathbb{R}^{m \times n}$ or of $\mathbb{C}^{m \times n}$. Note however that $G_2$ and $G_3$ have the same multilinear rank --- this is not coincidental but is a manifestation of the following result.
\begin{proposition}\label{prop:indmrank}
The multilinear rank of a tensor is independent of the choice of base field. If $\mathbb{K}$ is an extension field of $\Bbbk$, the value $\multirank(A)$ is the same whether $A$ is regarded as an element of $\Bbbk^{d_1 \times \dots \times d_k}$ or of $\mathbb{K}^{d_1 \times \dots \times d_k}$.
\end{proposition}
\begin{proof}
This follows immediately from~\eqref{eq:matrank} and the base field independence of matrix rank.
\end{proof}

In 1968, Bergman~\cite{Berg} considered linear subspaces of matrix spaces, and showed that the minimum rank on a subspace can become strictly smaller upon taking a field extension. He gave a class of examples, the simplest instance being the 2-dimensional subspace
\[
s
\begin{bmatrix}1 & 0\\ 0 & 1\end{bmatrix}
+ t
\begin{bmatrix}0 & 1\\-1 & 0\end{bmatrix}
\]
of~$\Rr^{2\times2}$. Every (nonzero) matrix in this subspace has rank~2, but the complexified subspace contains a matrix of rank~1.
Intriguingly, this example is precisely the subspace spanned by the slabs of~$G_3$. We suspect a deeper connection.

\subsection{Injectivity of orbits\label{subsec:inject}}

The tensor rank has the property of being invariant under the general
multilinear group (cf.\ \eqref{invtrank}). Indeed, much of its relevance comes from this fact.
Moreover, from Proposition~\ref{prop:rank} we know that tensor rank is
preserved when a tensor space is included in a larger tensor space. Similar
assertions are true for the multilinear rank (cf.\ \eqref{invmrank}).

The situation is more complicated for the function~$\Delta$ defined
on~$\Rr^{2\times2\times2}$. The sign of~$\Delta$ is $\GL_{2,2,2}(\Rr)$-invariant,
and $\Delta$ itself is invariant under $\Oo_{2,2,2}(\Rr)$. For general $d_1,d_2,d_3 \ge2$, we
do not have an obvious candidate function~$\Delta$ defined on~$\Rr%
^{d_{1}\times d_{2}\times d_{3}}$. However, there is a natural definition
of~$\Delta$ restricted to the subset of tensors~$A$ for which
$\multirank(A)\leq(2,2,2)$. Such a tensor can be expressed
as
\[
A=(L,M,N)\cdot (B \oplus 0)
\]
where $B\in\Rr^{2\times2\times2}$, $0 \in \Rr^{(d_1 - 2)\times (d_2 - 2) \times (d_3 - 2)}$
and $(L,M,N)\in\Oo_{d_{1},d_{2},d_{3}}(\Rr)$. We provisionally define $\Delta
(A)=\Delta(B)$, subject to a check that this is independent of the choices
involved. Given an alternative expression $A=(L^{\prime},M^{\prime},N^{\prime
})\cdot (B^{\prime}\oplus 0)$, it follows that $B \oplus 0$ and $B^{\prime} \oplus 0$ are in the same
$\Oo_{d_{1},d_{2},d_{3}}(\Rr)$-orbit. Indeed:
\[
B\oplus 0 =(L^{-1}L^{\prime},M^{-1}M^{\prime},N^{-1}N^{\prime})\cdot
(B^{\prime}\oplus 0).
\]
If we can show, more strongly, that $B,B^{\prime}$ belong to the same
$\Oo_{2,2,2}(\Rr)$-orbit, then the desired equality
$\Delta(B)=\Delta(B^{\prime})$ follows from the orthogonal invariance
of~$\Delta$.

The missing step is supplied by the next theorem, which we state in a
basis-free form for abstract vector spaces. If $V$ is a vector space, we write
$\GL(V)$ for the group of invertible linear maps from
$V\rightarrow V$. If, in addition, $V$ is an inner-product space, we write
$\Oo(V)$ for the group of norm-preserving linear maps $V\rightarrow V$.
In particular, $\GL(\Rr^{d})\cong\GL_{d}(\Rr)$ and $\Oo(\Rr^{d})\cong\Oo_{d}(\Rr)$.

\begin{theorem}
[injectivity of orbits]\label{thm:injectivity} Let $\Bbbk=\Rr$ or
$\mathbb{C}$ and $V_{1},\dots,V_{k}$ be $\mathbb{\Bbbk}$-vector spaces. Let
$U_{1}\leq V_{1},\dots,U_{k}\leq V_{k}$. (1)~Suppose $B,B^{\prime}\in
U_{1}\otimes\dots\otimes U_{k}$ are in distinct $\GL%
(U_{1})\times\dots\times\GL(U_{k})$-orbits of $U_{1}\otimes
\dots\otimes U_{k}$, then $B$ and~$B^{\prime}$ are in distinct
$\GL(V_{1})\times\dots\times\GL(V_{k})$-orbits of
$V_{1}\otimes\dots\otimes V_{k}$. (2)~Suppose $B,B^{\prime}\in U_{1}%
\otimes\dots\otimes U_{k}$ are in distinct $\Oo(U_{1})\times
\dots\times\Oo(U_{k})$-orbits of $U_{1}\otimes\dots\otimes U_{k}%
$, then $B$ and~$B^{\prime}$ are in distinct $\Oo(V_{1}%
)\times\dots\times\Oo(V_{k})$-orbits of $V_{1}\otimes\dots\otimes
V_{k}$.
\end{theorem}

\begin{lemma}
\label{lem:tuv} Let $W\leq U\leq V$ be vector spaces and $L\in
\GL(V)$. Suppose $L(W)\leq U$. Then there exists $\tilde{L}%
\in\GL(U)$ such that $L|_{W}=\tilde{L}|_{W}$. Moreover, if
$L\in\Oo(V)$ then we can take $\tilde{L}\in\Oo(U)$.
\end{lemma}

\begin{proof}
Extend $L|_{W}$ to~$U$ by mapping the orthogonal complement of~$W$ in~$U$
by a norm-preserving map to the orthogonal complement of~$L(W)$ in~$U$. The resulting linear map~$\tilde{L}$ has the desired properties and is orthogonal if $L$ is orthogonal.
\end{proof}

{\it Proof of Theorem~\ref{thm:injectivity}.}
We prove the contrapositive form of the theorem. Suppose $B^{\prime}=(L_{1},\dots ,L_{k})\cdot B$, where $L_i \in \GL(V_i)$.
Let $W_{i}\leq U_{i}$ be minimal subspaces such that $B$ is in the image of $W_{1}\otimes\dots\otimes W_{k}\hookrightarrow U_{1} \otimes\dots\otimes U_{k}$. It follows that $L_{i}(W_{i})\leq U_{i}$, for otherwise we could replace $W_{i}$ by $L_{i}^{-1}(L_{i}(W_{i})\cap U_{i})$. We
can now use Lemma~\ref{lem:tuv} to find $\tilde{L}_{i}\in\GL%
(U_{i})$ which agree with~$L_{i}$ on~$W_{i}$. By construction, $(\tilde{L}%
_{1},\dots,\tilde{L}_{k})\cdot B=(L_{1},\dots,L_{k})\cdot B=B^{\prime}$.
In the orthogonal case, where $L_i \in \Oo(V_i)$, we may choose $\tilde{L}_{i} \in \Oo(U_i)$.
\endproof

\begin{corollary}
Let $\varphi$ be a $\GL_{d_{1},\dots,d_{k}}(\Rr)$-invariant (respectively $\Oo_{d_{1},\dots,d_{k}}(\Rr)$-invariant) function on $\Rr^{d_{1}\times\dots\times d_{k}}$. Then $\varphi$ naturally extends to a $\GL_{d_{1},\dots,d_{k}}(\Rr)$-invariant (respectively $\Oo_{d_{1},\dots,d_{k}}(\Rr)$-invariant) function on the subset
\[
\{
A \in \Rr^{(d_{1}+e_{1})\times\dots\times(d_{k}+e_{k})}
\mid
\mbox{$r_i(A) \leq d_i$ for $i = 1, \dots, k$}
\}
\]
of $\Rr^{(d_{1}+e_{1})\times\dots\times(d_{k}+e_{k})}$.
\end{corollary}

\begin{proof}
As with $\Delta$ above, write $A=(L_{1},\dots,L_{k})\cdot B$ for
$B\in\Rr^{d_{1}\times\dots\times d_{k}}$ and define $\varphi
(A)=\varphi(B)$. By Theorem~\ref{thm:injectivity} this is independent of the
choices involved.
\end{proof}

The problem of classification is closely related to finding invariant
functions. We end this section with a strengthening of
Theorem~\ref{thm:orbits}.

\begin{corollary}
The eight orbits in Theorem~\ref{thm:orbits} remain distinct under the
embedding $\Rr^{2\times2\times2}\hookrightarrow\Rr^{d_{1}\times
d_{2}\times d_{3}}$ for any $d_{1},d_{2},d_{3}\geq2$. Thus,
Theorem~\ref{thm:orbits} immediately gives a classification of tensors
$A\in\Rr^{d_{1}\times d_{2}\times d_{3}}$ with
$\multirank(A)\leq(2,2,2)$, into eight classes under
$\GL_{d_{1},d_{2},d_{3}}(\Rr)$-equivalence.
\end{corollary}

The corollary allows us to extend the notion of tensor-type to
$\Rr^{d_1 \times d_2 \times d_3}$. For instance, we will say that $A
\in \Rr^{d_1\times d_2 \times d_3}$ has type~$G_3$ iff $A$ is
$\GL$-equivalent to $G_3 \in \Rr^{2\times 2 \times 2} \subset
\Rr^{d_1 \times d_2 \times d_3}$.

Note that order-$k$ tensors can be embedded in order-$(k+1)$ tensors by taking the tensor product with a 1-dimensional factor. Distinct orbits remain distinct, so the results of this subsection extend to inclusions into tensor spaces of higher order.

\section{Volume of tensors with no optimal low-rank
approximation\label{sec:volume}}

At this point, it is clear that there exist tensors that can fail to
have optimal low-rank approximations. However it is our experience
that practitioners have sometimes expressed optimism that such
failures might be rare abnormalities that are not encountered in
practice. In truth, such optimism is misplaced: the set of tensors
with no optimal low-rank approximation has positive volume. In other
words, a randomly chosen tensor will have a non-zero chance of
failing to have a optimal low-rank approximation.

We begin this section with a particularly striking instance of this.

\begin{theorem}\label{thm:norank2}
No tensor of rank~$3$ in~$\Rr^{2 \times 2 \times 2}$ has an optimal
rank-$2$ approximation (with respect to the Frobenius norm). In
particular, \app$(A,2)$ has no solution for tensors of type~$G_3$,
which comprise a set that is open and therefore of positive volume.
\end{theorem}

\begin{lemma}
\label{lem:r-1}
Let $A \in \Rr^{d_1 \times \dots \times d_k}$ with $\trank(A) \geq r$. Suppose $B \in \Ss_r(d_1, \dots, d_k)$ is an optimal rank-$r$ approximation for~$A$. Then $\trank(B) = r$.
\end{lemma}

\begin{proof}
Suppose $\trank(B) \leq r-1$. Then $B \not= A$, and so $B-A$ has at
least one nonzero entry in its array representation. Let $E \in
\Rr^{d_1 \times \dots \times d_k}$ be the rank-$1$ tensor which
agrees with $B-A$ at that entry and is zero everywhere else. Then
$\trank(B+E) \leq r$ but $\|A - (B+E)\|_F < \|A - B\|_F$, so $B$
is not optimal.
\end{proof}

{\it Proof of Theorem~\ref{thm:norank2}.}
Let $A \in \Rr^{2\times 2\times 2}$ have rank~$3$,
%
%
and suppose that $B$ is an optimal rank-$2$ approximation to~$A$.
Propositions \ref{prop:rank2:b} and~\ref{prop:rank2:a}, together with the continuity of~$\Delta$, imply that $\Delta(B)=0$. Lemma~\ref{lem:r-1} implies that $\trank(B) = 2$.
By Theorem~\ref{thm:orbits}, it follows that $B$ is of type $D_2$, $D_2'$
or~$D_2''$.

We may assume without loss of generality that $B$ is of type~$D_2$.
The next step is to put $B$ into a helpful form by making an orthogonal change of
coordinates. This gives an equivalent approximation problem, thanks to the $\Oo$-invariance of the Frobenius norm.
From Table~\ref{table:orbits}, we know that $\multirank(B) = (1,2,2)$. Such a~$B$ is
orthogonally equivalent to a tensor of the following form:
\begin{equation}
\label{eq:B122}
\left[
\begin{array}[c]{cc}\lambda & 0\\0 & \mu \end{array}
\right\vert \!\left.
\begin{array}[c]{cc}0 & 0\\0 & 0\end{array}
\right]
\end{equation}
Indeed, a rotation in the first tensor factor brings~$B$ entirely into the first slab, and further rotations in the second and third factors put the resulting matrix into diagonal form, with singular values $\lambda, \mu \not= 0$.

Henceforth we assume that $B$ is equal to the tensor in \eqref{eq:B122}.
We will consider perturbations of the form $B + \epsilon H$, which will be chosen so that
$\Delta(B + \epsilon H) = 0$ for all $\epsilon \in \Rr$. Then $B + \epsilon H \in \Ssbar_2(2,2,2)$, and we must have
\[
\lVert A - B \rVert_F \leq \lVert A - (B + \epsilon H) \rVert_F
\]
for all~$\epsilon$. In fact
\[
\lVert A - (B + \epsilon H) \rVert_F^2 - \lVert A - B \rVert_F^2
= -2 \epsilon \langle A-B, H \rangle_F + \epsilon^2 \lVert H \rVert_F^2
\]
so if this is to be nonnegative for all small values of~$\epsilon$, it is necessary that
\begin{equation}
\label{eqn:dot}
\langle A - B, H \rangle_F = 0.
\end{equation}

Tensors~$H$ which satisfy the condition $\Delta(B + \epsilon H) \equiv 0$ include the following:
\[
\left[
\begin{array}[c]{cc}\times & \times \\ \times & \times \end{array}
\right\vert \!\left.
\begin{array}[c]{cc}0 & 0\\0 & 0\end{array}
\right]
,
\quad
\left[
\begin{array}[c]{cc}0 & 0\\0 & 0 \end{array}
\right\vert \!\left.
\begin{array}[c]{cc}0 & 1\\0 & 0\end{array}
\right]
,
\quad
\left[
\begin{array}[c]{cc}0 & 0\\0 & 0 \end{array}
\right\vert \!\left.
\begin{array}[c]{cc}0 & 0\\1 & 0\end{array}
\right]
,
\quad
\left[
\begin{array}[c]{cc}0 & 0\\0 & 0 \end{array}
\right\vert \!\left.
\begin{array}[c]{cc}\lambda & 0\\0 & \mu \end{array}
\right]
\]
since the resulting tensors have types $D_2$, $D_3$, $D_3$, and~$D_2$ respectively.

Each of these gives a constraint on $A-B$, by virtue of~\eqref{eqn:dot}. Putting the constraints together, we find that
\[
A - B =
\left[
\begin{array}[c]{cc}0 & 0\\0 & 0 \end{array}
\right\vert \!\left.
\begin{array}[c]{cc}a\mu & 0\\0 & -a\lambda \end{array}
\right]
\qquad
\mbox{or}
\qquad
A =
\left[
\begin{array}[c]{cc}\lambda & 0\\0 & \mu \end{array}
\right\vert \!\left.
\begin{array}[c]{cc}a\mu & 0\\0 & -a\lambda \end{array}
\right]
\]
for some $a \in \Rr$. Thus $A = (\lambda \ee_1 + a\mu \ee_2)\otimes \ee_1
\otimes \ee_1 + (\mu \ee_1 - a\lambda \ee_2)\otimes \ee_2 \otimes \ee_2$
has rank~$2$, a contradiction. %
\endproof

\begin{corollary}
\label{cor:g3}
Let $d_1, d_2, d_3 \geq 2$. If $A \in \Rr^{d_1 \times d_2 \times d_3}$ is
of type~$G_3$, then $A$ does not have an optimal rank-$2$ approximation.
\end{corollary}

\begin{proof}
We use the projection $\Pi_A$ defined in subsection~\ref{subsec:multiprojection}. For any $B \in \Rr^{d_1 \times d_2 \times d_3}$, Pythagoras' theorem \eqref{eq:pythagoras}
gives:
\begin{eqnarray*}
\| B - A \|_F^2
&=& \| \Pi_A(B - A) \|_F^2 + \| (1 - \Pi_A) (B - A) \|_F^2
\\
&=& \| \Pi_A(B) - A \|_F^2 + \| B - \Pi_A(B) \|_F^2
\end{eqnarray*}
If $B$ is an optimal rank-$2$ approximation, then it follows that $B =
\Pi_A(B)$; for otherwise $\Pi_A(B)$ would be a better approximation. Thus $B \in U_1 \otimes U_2 \otimes U_3$, where $U_1, U_2, U_3$ are the supporting subspaces of~$A$. These are 2-dimensional, since $\multirank(A) = (2,2,2)$, so $U_1 \otimes U_2 \otimes U_3 \cong \Rr^{2 \times 2 \times 2}$. The optimality of~$B$ now contradicts Theorem~\ref{thm:norank2}.
\end{proof}


Our final result is that the set of tensors~$A$ for which
\app$(A,2)$ has no solution is a set of positive volume, for all
tensor spaces of order~3 except those isomorphic to a matrix space;
in other words, Theorem~\ref{thm:posvol}. Note that the $G_3$-tensors comprise a set of zero volume in all cases except $\Rr^{2 \times 2 \times 2}$.
Here is the precise statement.

\begin{theorem}
\label{thm:posvol2} Let $d_1, d_2, d_3 \geq 2$. The set of tensors
$A \in \Rr^{d_1 \times d_2 \times d_3}$ for which \app$(A,2)$ does
not have a solution (in the Frobenius norm) contains an open
neighborhood of the set of tensors of type~$G_3$. In particular,
this set is nonempty and has positive volume.
\end{theorem}

For $A\in \Rr^{d_1 \times d_2 \times d_3}$, let $\Bb(A)$ denote the set of
optimal border-rank-$2$ approximations for~$A$. Since $\Ssbar_2(d_1, d_2,
d_3)$ is nonempty and closed, it follows that $\Bb(A)$ is nonempty and compact.

We can restate the theorem as follows.
Let $A_0$ be an arbitrary $G_3$-tensor. We must show that if $A$ is close to~$A_0$, and $B \in \Bb(A)$, then $\trank(B) > 2$, i.e.\ $B$ is a $D_3$-tensor. Our proof strategy is contained in the steps of the following lemma.
%
\begin{lemma}
\label{lem:final}
Let $A_0 \in \Rr^{d_1 \times d_2 \times d_3}$ be a fixed tensor of type $G_3$. Then there exist positive numbers $\rho = \rho(A_0)$, $\delta = \delta(A_0)$ such that the following statements are true for all $A \in \Rr^{d_1 \times d_2 \times d_3}$.
\begin{enumerate}
\item[{\rm(1)}]
If $A$ is a $G_3$-tensor and $B \in \Bb(A)$, then $B$ is a $D_3$-tensor and $\Pi_B = \Pi_A$.
\item[{\rm(2)}]
If $\| A - A_0 \|_F < \rho$ and $\multirank(A) \leq(2,2,2)$, then $A$ is a $G_3$-tensor.
\item[{\rm(3)}]
If $\| A - A_0\|_F < \delta$ and $B \in \Bb(A)$, define $A' = \Pi_B(A)$. Then $\| A' - A_0\|_F < \rho$ and $B \in \Bb(A')$.
\end{enumerate}
\end{lemma}

%



%

{\it Proof of Theorem~\ref{thm:posvol2}, assuming Lemma~\ref{lem:final}.}
Fix $A_0 \in \Rr^{d_1 \times d_2 \times d_3}$ and suppose $\| A - A_0\|_F < \delta$. It is not generally true that $\multirank(A) \leq (2,2,2)$, so we cannot apply~(2) directly to~$A$. Let $B \in \Bb(A)$. Then $A' = \Pi_B(A)$ is close to~$A_0$, by~(3). Since $\multirank(B) \leq (2,2,2)$ and $\Pi_B$ is the projection onto the subspace spanned by~$B$, it follows that $\multirank(A') \leq (2,2,2)$. Now (2) implies that $A'$ is a $G_3$-tensor. Since $B \in \Bb(A')$, by~(3), it follows from~(1) that $B$ is a $D_3$-tensor.
\endproof

{\it Proof of Lemma~\ref{lem:final}, (1).}
This is essentially Corollary~\ref{cor:g3}: $B$ cannot have rank~$2$ or
less, but it has border-rank~$2$, so $B$ must be a $D_3$-tensor. Since $B
=
\Pi_A(B)$ it follows that the supporting subspaces of~$B$ are contained in the supporting subspaces of~$A$. However, $\multirank(B) = (2,2,2) = \multirank(A)$, so the two tensors must have the same supporting subspaces, and so $\Pi_B = \Pi_A$.
\endproof

{\it Proof of Lemma~\ref{lem:final}, (2).}
Let $\Ssbar_2^+(d_1,d_2,d_3)$ denote the set of non-$G_3$ tensors in $\Rr^{d_1 \times d_2 \times d_3}$ with $\multirank \leq (2,2,2)$. Since $A_0 \not\in \Ssbar_2^+(d_1, d_2, d_3)$, it is enough to show that $\Ssbar_2^+(d_1, d_2, d_3)$ is closed, for then it would be disjoint from the $\rho$-ball about~$A_0$, for some $\rho > 0$. Note that
\[
\Ssbar_2^+(d_1, d_2, d_3) = \Oo_{d_1,d_2,d_3}(\Rr) \cdot \Ssbar_2^+(2,2,2).
\]
Now $\Ssbar_2^+(2,2,2) = \{ A \in \Rr^{2 \times 2 \times 2} \mid \Delta(A) \geq 0\}$ is a closed subset of $\Rr^{2 \times 2 \times 2}$, and the action of the compact group $\Oo_{d_1, d_2, d_3}(\Rr)$ is proper. It follows that $\Ssbar_2^+(d_1, d_2, d_3)$ is closed, as required.
\endproof

{\it Proof of Lemma~\ref{lem:final}, (3).}
We begin with the easier part of the statement, which is that $B \in \Bb(A')$.
To prove this, we will show that $\|A' - B\|_F \leq \|A' - B'\|_F$ whenever $B' \in \Bb(A')$, establishing the optimality of~$B$ as an approximation to~$A'$. Accordingly, let $B' \in \Bb(A')$. Since $\Pi_B(A') = A'$, it follows from \eqref{eq:pythagoras} with $\Pi_B$ that
\[
\| A' - B' \|_F^2 = \| A' - \Pi_B(B')\|_F^2 + \| B' - \Pi_B(B') \|_F^2
\]
so, since $B'$ is optimal, we must have $\Pi_B(B') = B'$. We can now apply \eqref{eq:pythagoras} with $\Pi_B$ to both sides of the inequality
$\| A - B \|_F^2 \leq \| A - B' \|_F^2$
to get
\[
\| A' - B\|_F^2 + \| A - A' \|_F^2
\leq
\| A' - B'\|_F^2 + \| A - A' \|_F^2
\]
and hence $\|A' - B\|_F \leq \|A' - B'\|_F$, as claimed.

We now turn to the proof that $\Pi_B(A)$ is close to~$A_0$ if $A$ is close to $A_0$. This is required to be uniform in $A$ and $B$. In other words, there exists $\delta = \delta(A_0) > 0$ such that for all $A$ and all $B \in \Bb(A)$, if $\|A - A_0\|_F < \delta$ then $\|\Pi_B(A) - A_0\| < \rho$. Here $\rho = \rho(A_0)$ is fixed from part~(2) of this lemma.

We need control over the location of~$B$. Let $\Bbe(A_0)$ denote the $\epsilon$-neighborhood of~$\Bb(A_0)$ in~$\Ssbar_2(d_1, d_2, d_3)$.

\begin{proposition}
\label{prop:delta2} Given $\epsilon > 0$, there exists $\delta > 0$
such that if $\|A - A_0\|_F < \delta$ then $\Bb(A) \subset
\Bbe(A_0)$.
\end{proposition}

\begin{proof}
The set $\Ssbar_2(d_1,d_2,d_3) \setminus \Bbe(A_0)$ is closed, and so it attains
its minimum distance from~$A_0$. This must exceed the absolute
minimum $\|A_0 - B_0\|_F$ for $B_0 \in \Bb(A_0)$ by a positive
quantity~$2\delta$, say. If $\|A - A_0\|_F < \delta$ and $B' \in
\Ssbar_2(d_1,d_2,d_3) \setminus \Bbe(A_0)$ then
\begin{eqnarray*}
\| A - B' \|_F &\geq& \| B' - A_0 \|_F - \|A - A_0\|_F
\\
&>& \|A_0 - B_0\|_F + 2\delta - \delta
\\
&=& \|A_0 - B_0\|_F + \delta
\\
&>& \|A_0 - B_0\|_F + \| A - A_0 \|_F
\\
&\geq& \|A - B_0\|_F
\end{eqnarray*}
using the triangle inequality in the first and last line. Thus $B' \not\in \Bb(A)$.
\end{proof}

We claim that if $\epsilon$ is small enough, then $\multirank(B) = (2,2,2)$ for all $B \in \Bbe(A_0)$. Indeed, this is
already true on $\Bb(A_0)$, by part~(1). Since $\multirank$ is upper-semicontinuous and does not exceed $(2,2,2)$ on~$\Ssbar_2(d_1,d_2,d_3)$, it must be constant on a neighborhood of~$\Bb(A_0)$ in~$\Ssbar_2(d_1,d_2,d_3)$. Since $\Bb(A_0)$ is compact, the neighborhood can be taken to be an $\epsilon$-neighborhood.

Part~(1) implies that $\Pi_{B_0} = \Pi_{A_0}$ for all $B_0 \in \Bb(A_0)$. If $\epsilon$ is small enough that $\multirank(B) = (2,2,2)$ on $\Bbe(A_0)$, then $\Pi_B$ depends continuously on~$B \in \Bbe(A_0)$, by Proposition~\ref{prop:PiAconst}. Since $\Bb(A_0)$ is compact, we can choose $\epsilon$ small enough so that the operator norm of $\Pi_B - \Pi_{A_0}$ is as small as we like, uniformly over~$\Bbe(A_0)$.


We are now ready to confine $\Pi_B(A)$ to the $\rho$-neighborhood of~$A_0$.
Suppose, initially, that $\| A - A_0\|_F \leq \rho/2$ and $B \in \Bbe(A_0)$. Then
\begin{align*}
\| \Pi_B(A) - A_0 \|_F &\leq \| (\Pi_B - \Pi_{A_0})\cdot A \|_F + \|
\Pi_{A_0} \cdot A - A_0\|_F \notag
\\
&\leq \| \Pi_B - \Pi_{A_0}\| \|A\|_F + \| \Pi_{A_0}\cdot(A -
A_0)\|_F \notag
\\
&\leq \| \Pi_B - \Pi_{A_0}\| (\|A_0\|_F + \rho/2) + \|A - A_0\|_F
\notag
\\
&\leq \| \Pi_B - \Pi_{A_0}\| (\|A_0\|_F + \rho/2) + \rho/2
\end{align*}
Now choose $\epsilon > 0$ so that the
operator norm $\|\Pi_B - \Pi_{A_0}\|$ is kept small enough to guarantee that the
right-hand side is less than~$\rho$. For this~$\epsilon$, choose
$\delta$ as given by Proposition~\ref{prop:delta2}. Ensure also that $\delta < \rho/2$.

Then, if $\| A - A_0 \|_F < \delta$ and $B \in \Bb(A)$, we have $B \in \Bbe(A_0)$. By the preceding calculation, $\| A' - A_0 \|_F < \rho$. This completes the proof.
\endproof


\section{Closing remarks}

We refer interested readers to \cite{CGLM1, CGLM2, L1, LG} for a discussion of
similar issues for symmetric tensors and nonnegative tensors. In particular,
the reader will find in \cite{CGLM2} an example of a symmetric tensor of
symmetric rank~$r$ (may be chosen to be arbitrarily large) that does not
have a best symmetric-rank-$2$ approximation. In \cite{L1, LG}, we show that such
failures do not occur in the context of nonnegative tensors --- a nonnegative
tensor of nonnegative-rank~$r$ will always have a best
nonnegative-rank-$s$
approximation for any $s\leq r$.

In this paper we have focused our attention on the real case; the complex
case has been studied in great detail in algebraic computational
complexity theory and algebraic geometry. For the interested reader, we note that the rank-jumping phenomenon still occurs: Proposition~\ref{prop:dSLtensor} and its proof carry straight through to the complex case. On the other hand, there is no distinction between $G_3$ and $G_2$ tensors over the complex numbers; if $\Delta(A) \ne 0$ then $A$ has rank~$2$. The results of Section~\ref{sec:volume} have no direct analogue.

The major open question in tensor approximation is how to overcome the ill-posedness of {\brap}. In general this will conceivably require an equivalent of Theorem~\ref{thm:b} that characterizes the limit points of rank-$r$ order-$k$ tensors. It is our hope that some
of the tools developed in our study, such as Theorems~\ref{thm:reduction:k} and \ref{thm:injectivity} (both of which apply to general $r$ and $k$), may be used in future studies. The type of characterization in Corollary~\ref{cor:b}, for $r=2$ and
$k = 3$, is an example of what one might hope to achieve.

\subsection*{Acknowledgements}
We thank the anonymous reviewers for some exceptionally helpful comments. We also gratefully acknowledge Joseph Landsberg and Bernd Sturmfels for enlightening pointers that helped improved Sections \ref{subsec:derive} and~\ref{sect:SA}.
Lek-Heng Lim thanks Gene Golub for his encouragement and helpful discussions.
Both authors thank Gunnar Carlsson and the Department of Mathematics, Stanford University, where some of this work was done.
Lek-Heng~Lim is supported by the Gerald~J.~Lieberman Fellowship from Stanford University.  Both authors have been partially supported by the DARPA Grant 32905 and the NSF Grant DMS 01-01364.

\end{document}